\documentclass[aop,MSNbibl,citesort,dvips]{arximspdf}

\doi{10.1214/11-AOP728}
\volume{41}
\issue{2}
\pubyear{2013}
\firstpage{527}
\lastpage{618}

\makeatletter
\newcommand{\R} {{\mathbb R}}
\newcommand{\N} {{\mathbb N}}

\newtheorem{Theo}{Theorem}
\newtheorem{Lemma}[Theo]{Lemma}
\newtheorem{Cor}[Theo]{Corollary}
\newtheorem{Prop}[Theo]{Proposition}
\newproclaim{Rmk}[Theo]{Remark}
\def\eps{\varepsilon}

\newcommand{\cal}{\mathcal}
\newcommand{\eqref}[1]{(\ref{#1})}
\newcommand{\fracd}[2]{({#1}/{#2})}
\makeatother

\begin{document}
\begin{frontmatter}

\title{The genealogy of branching Brownian motion with~absorption}
\runtitle{Genealogy of branching Brownian motion}

\begin{aug}
\author[A]{\fnms{Julien} \snm{Berestycki}\thanksref{au1}\ead[label=e1]{julien.berestycki@upmc.fr}},
\author[B]{\fnms{Nathana\"el}~\snm{Berestycki}\corref{}\thanksref{au2}\ead[label=e2]{N.Berestycki@statslab.cam.ac.uk}}\\
\and
\author[C]{\fnms{Jason} \snm{Schweinsberg}\thanksref{au3}\ead[label=e3]{jschwein@math.ucsd.edu}}
\runauthor{J. Berestycki, N. Berestycki and J. Schweinsberg}
\affiliation{Universit\'e Paris VI, University of Cambridge and~University~of~California~at~San~Diego}
\address[A]{J. Berestycki\\
LPMA/UMR 7599\\
Universit\'{e} Pierre et Marie Curie (P6)---\\
\quad Boite courrier 188\\
75252 PARIS Cedex 05\\
France\\
\printead{e1}}
\address[B]{N. Berestycki\\
Statistical Laboratory, DPMMS\\
Cambridge University\\
Wilberforce Rd.\\
Cambridge CB3 0WB\\
United Kingdom\\
\printead{e2}}
\address[C]{J. Schweinsberg\\
Department of Mathematics, 0112\\
University of California, San Diego\\
9500 Gilman Drive\\
La Jolla, California 92093-0112\\
USA\\
\printead{e3}} %adresu isvedimo komanda gale!
\end{aug}
\thankstext{au1}{Supported in part by ANR Grant BLAN06-3-146282 MAEV.}
\thankstext{au2}{Supported in part by EPSRC Grant EP/GO55068/1.}
\thankstext{au3}{Supported in part by NSF Grant DMS-08-05472.}

% HISTORY:
\received{\smonth{2} \syear{2010}}
\revised{\smonth{6} \syear{2011}}

% ABSTRACT
%
\begin{abstract}
We consider a system of particles which perform branching Brownian
motion with negative drift and are killed upon reaching zero, in the
near-critical regime where the total population stays roughly constant
with approximately $N$ particles. We show that the characteristic time
scale for the evolution of this population is of order $(\log N)^3$, in
the sense that when time is measured in these units, the scaled number
of particles converges to a variant of Neveu's continuous-state
branching process. Furthermore, the genealogy of the particles is then
governed by a coalescent process known as the Bolthausen--Sznitman
coalescent. This validates the nonrigorous predictions by Brunet,
Derrida, Muller and Munier for a closely related model.
\end{abstract}

% KEYWORDS
%

\end{frontmatter}

\tableofcontents

%s1 ###
%s1 #&#
\section{Introduction}\label{sec1}

Branching Brownian motion is a stochastic process in which, at time
zero, there is a single particle at the origin. Each particle moves
according to a standard Brownian motion for an exponentially
distributed time with mean one, at which point it splits into two
particles. Early work on branching Brownian motion, going back to
McKean~\cite{mckean}, focused on the position $M(t)$ of the right-most
particle. Bramson~\cite{bram1,bram2} obtained asymptotics for the
median of the distribution of $M(t)$, and Lalley and Sellke~\cite{ls1}
found the asymptotic distribution of $M(t)$.

In 1978, Kesten~\cite{kesten} introduced branching Brownian motion with
absorption. This process follows the same dynamics as branching
Brownian motion except that the initial particle is located at $x > 0$,
the Brownian particles have a drift of $-\mu$, where $\mu> 0$, and
particles are killed when they reach the origin. Kesten showed that
there exists a critical value $\mu_c = \sqrt{2}$ such that if $\mu
\geq
\mu_c$, then the process dies out almost surely, while if $\mu< \mu
_c$, the process survives with positive probability. More recent work
on this process can be found in~\cite{hh07} and~\cite{hhk06}.

Our interest in branching Brownian motion with absorption comes from
its possible interpretation as a model of a population undergoing
selection. To see this connection, imagine that each individual in a
population is represented by a position on the real line, which
measures her fitness. The fitness of an individual evolves according to
Brownian motion due to mutations, and initially the fitness of a child
is identical to the fitness of the parent. Selection progressively
eliminates all individuals whose fitness becomes too low; we
effectively imagine selection as a moving wall with constant speed $\mu
$. Every individual whose fitness falls beyond the current threshold is
instantly removed from the population.

To obtain asymptotic results as the population size tends to infinity,
we consider a sequence of branching Brownian motions with absorption.
For each positive integer $N$, we have a branching Brownian motion with
absorption $(X_N(t), t \geq0)$. We consider the near-critical case,
where the drift $\mu$ depends on $N$ and for $N \geq2$,
%
%e1 ###
%
%e1 #&#
\begin{equation}\label{mu}
\mu= \sqrt{2- \frac{2\pi^2}{(\log N + 3 \log\log N)^2}}.
\end{equation}
We also start the process with many particles, rather than just one, at
time zero, and we make some rather technical assumptions on the initial
conditions, which are given later in Proposition~\ref{ZNCSBP}. While
(\ref{mu}) and the initial conditions may seem unnatural, they are
necessary to ensure that the number of particles in the system stays of
order $N$ on the time scale of interest, so that the process can be
viewed as a model of a population of size approximately~$N$.

We focus on understanding the genealogy of a sample from the population
after a large time. We show that the time to the most recent common
ancestor of a sample behaves like $(\log N)^3$. Moreover we identify
the limiting geometry of the coalescence tree of a sample, which we
show is governed by a coalescent process $(\Pi(t), t \geq0)$ known as
the Bolthausen--Sznitman coalescent. The Bolthausen--Sznitman
coalescent, which is defined precisely in  Section~\ref{boszsec}, is
a coalescent process that allows many ancestral lines to merge at once.
This result is in sharp contrast with the standard case of the Moran
model or the Wright--Fisher model, where random genetic drift leads to
a characteristic genealogical time of $N$ generations and a
genealogical tree given by Kingman's coalescent, which permits only
pairwise mergers of ancestral lines.

The main result of this paper can thus be stated as follows. Fix $t >
0$. Choose $n$ particles uniformly at random from the population at
time $(\log N)^3 t$, and label these particles at random by the
integers $1,\ldots, n$. For $0 \leq s \leq2 \pi t$, define $\Pi_N(s)$
to be the partition of $\{1,\ldots, n\}$ such that $i$ and $j$ are in
the same block of $\Pi_N(s)$ if and only if the particles labeled $i$
and $j$ are descended from the same ancestor at time
$(t - s/2 \pi) (\log N)^3$. This is the standard ``ancestral
partition''
of the sample. Then, with our initial conditions, we have the following
result, which is stated precisely later as Theorem~\ref{bosz}.

\begin{result*} The sequence of processes $(\Pi_N(s), 0
\le s \le2 \pi t)$ converges in the sense of finite-dimensional
distributions as $N \to\infty$ to the Bolthausen--Sznitman coalescent
$(\Pi(s), 0 \leq s \leq2 \pi t)$.
\end{result*}

The reason for the multiple mergers is that when a particle gets very
far to the right [in fact, at position $\frac1{\sqrt{2}}(\log N + 3
\log\log N + O(1))$], many descendants of this particle survive for a
long time, as they are able to avoid being killed at zero. They quickly
generate a positive fraction of the population. As a result, when a
sample of particles is taken far into the future, many of their
ancestral lines get traced back to this particle and coalesce at nearly
the same time. Our result is accompanied by Theorem~\ref{MNCSBP}, which
gives the evolution of the total number of particles $M_N(t)$ in the
system. Under the same assumptions, $M_N((\log N)^3t)/(2 \pi N)$
converges in the sense of finite-dimensional distributions toward a
continuous-state branching process with branching mechanism $\Psi(u) =
au + 2\pi^2 u \log u$ for some constant $a \in\R$.

%s1.1 ###
%s1.1 #&#
\subsection{Related models and conjectures}\label{sec11}

Our inspiration for this model comes from the work of Brunet et al.
\cite{bdmm1,bdmm2} concerning the effect of natural selection on the
genealogy of a population. They considered a model of a population with
fixed size $N$ in which each individual has a fitness. They assumed
that each individual has $k \geq2$ offspring in the next generation,
and that the fitness of each offspring is the parent's fitness plus an
independent random variable with some distribution $\mu$. Of the $kN$
offspring, the $N$ with the highest fitness survive to form the next
generation. This process repeats itself in each generation. Brunet et
al.~\cite{bdmm1,bdmm2} gave a detailed and intricate, but not
mathematically rigorous, analysis of this model and arrived at the
following three conjectures:
\begin{longlist}[(3)]
\item[(1)] If $L_m$ is the maximum of the fitnesses of the $N$
individuals in generation $m$, then $L_m/m$ converges almost surely to
some limiting velocity $v_N$. Furthermore, the limit $v_{\infty} =
\lim
_{N \rightarrow\infty} v_N$ exists, and there is a constant $C$ such
that
%
%e2 ###
%
%e2 #&#
\begin{equation}\label{vcorrec}
v_{\infty} - v_N \sim\frac{C}{(\log N)^2}.
\end{equation}

\item[(2)] If two individuals are sampled from the population at random
in some generation, then the number of generations that we need to look
back to find their most recent common ancestor is of order $(\log N)^3$.

\item[(3)] If $n$ individuals are sampled from the population at random
in some generation, and their ancestral lines are traced backwards in
time, the coalescence of these lineages can be described by the
Bolthausen--Sznitman coalescent.
\end{longlist}

This model is similar to a branching random walk in which the positions
of the particles correspond to the fitnesses of the individuals.
Indeed, this model would be precisely a branching random walk if all
individuals were permitted to survive. The limiting velocity $v_{\infty
}$ that appears in the first conjecture is the limiting velocity of the
right-most particle in branching random walk, which was studied in the
1970s by Kingman~\cite{king}, Hammersley~\cite{hamm} and Biggins
\cite
{bigg}. Interest in variations of the branching random walk in which
the number of particles stays fixed is more recent. B\'erard and Gou\'
er\'e~\cite{bego08} recently proved the first conjecture in the form
stated above, in the case $k = 2$, under suitable regularity conditions
on $\mu$. Their proof builds on previous work of Gantert, Hu and Shi
\cite{ghs} and Pemantle~\cite{pem}. See also the work of Durrett and
Mayberry~\cite{dumay}, who considered a model very similar to this one
while studying predator-prey systems, and Durrett and Remenik~\cite{durem}.

The analysis of Brunet et al. involves studying solutions $u(x,t)$ to
the noisy FKPP equation
%
%e3 ###
%
%e3 #&#
\begin{equation}\label{noiseKPP}
\frac{\partial u}{\partial t} = \frac{\partial^2 u}{\partial x^2} +
u -
u^2 + \sqrt{\frac{u(1-u)}{N}} W(x,t),
\end{equation}
where $W(x,t)$ is space--time white noise. If the noise term were
removed, this partial differential equation would be the well-known
FKPP equation, which was introduced in 1937 by Fisher~\cite{fisher} and
by Kolmogorov, Petrovskii and Piscunov~\cite{kpp} and is one of the
simplest nonlinear partial differential equations that admits traveling
wave solutions. The link between the FKPP equation and branching
Brownian motion has been known since the work of McKean~\cite{mckean},
who showed that if $M(t)$ denotes the position of the right-most
particle at time $t$ for branching Brownian motion with variance
parameter $2$ and $u(t,x) = P(M(t) > x)$, then $u$ is the unique
solution to the FKPP equation with the initial condition $u(0,x) =
\mathbf{1}_{\{x < 0\}}$. In~\cite{hhk06}, Harris, Harris and Kyprianou
use branching Brownian motion with absorption to give a probabilistic
analysis of solutions to the FKPP equation.

The first conjecture above can also be expressed as a conjecture about
the velocity of solutions to equations such as (\ref{noiseKPP}). This
form of the conjecture goes back to the work~\cite{bd1} of Brunet and
Derrida, who refined their analysis and simulations in~\cite{bd2,bd3}.
In this form, the conjecture states that the velocity of traveling wave
solutions to the original FKPP equation exceeds the velocity of
solutions to equation (\ref{noiseKPP}) with the noise term by a
quantity that is of the order $1/(\log N)^2$. Recently Mueller, Mytnik
and Quastel~\cite{mmq} proved this result.

With the first conjecture having been largely settled, the purpose of
the present paper is to provide rigorous versions of the second and
third conjectures. As explained above, the model that we work with is
not exactly the model studied in~\cite{bdmm1,bdmm2}. Instead, to
simplify the analysis, we replace branching random walk by branching
Brownian motion, and rather than keeping the population size exactly
fixed, we control the population size by killing particles that drift
too far to the left. Note in particular that with our choice (\ref{mu}),
%
%e4 ###
%
%e4 #&#
\begin{equation}\label{muasymp}
\mu_c - \mu\sim\frac{\pi^2}{ \sqrt{2}(\log N)^2}
\end{equation}
as $N \to\infty$, which matches precisely (\ref{vcorrec}) in
Conjecture 1 above. Models with nonconstant population size where
already discussed by Derrida and Simon in~\cite{ds1,ds2} using
nonrigorous methods. Although they do not study genealogies, their
analysis strongly suggests that a result similar to ours may be expected.

Another question of interest related to these nearly critical branching
particle systems (in the sense that the drift $\mu$ of particles is
slightly above the critical value $\mu_c= \sqrt{2}$), concerns
asymptotics for the survival probability. This is a topic that has
attracted a considerable amount of attention in recent years; see, for
example,~\cite{aija11,bego08,bego09,fz,ghs,hhk06,jaffuel}. In
\cite{bbs}, we use the techniques developed here to derive fairly sharp
estimates for the survival probability of nearly critical branching
Brownian motion.

We also emphasize that the Bolthausen--Sznitman coalescent describes
precisely the ultrametric structure that is expected to emerge in the
low-temperature regime of mean-field spin glass models such as the
well-known Sherrington--Kirkpatrick model. This is perhaps not a
coincidence, as the model which we study here may be seen as a
degenerate form of spin glass models, with the position of the
particles being approximately given by a Gaussian field with a
covariance structure which is closely related to their genealogy.

%s1.2 ###
%s1.2 #&#
\subsection{Continuous-state branching processes}\label{sec12}

A continuous-state branching process is a $[0, \infty]$-valued Markov
process $(Z(t), t \geq0)$ whose transition functions $p_t(x, \cdot)$ satisfy
\[
p_t(x + y, \cdot) = p_t(x, \cdot) * p_t(y, \cdot) \qquad\mbox{for all
}x, y \geq0.
\]
That is, the sum of independent copies of the process started from $x$
and $y$ has the same law as the process started from $x + y$.
Continuous-state branching processes were introduced by Jirina \cite
{jirina}. Lamperti~\cite{lamp1} showed that continuous-state branching
processes are precisely the processes that can be obtained by taking
scaling limits of Galton--Watson processes. Lamperti~\cite{lamp2} and
Silverstein~\cite{silv} observed a one-to-one correspondence between
continuous-state branching processes and L\'evy processes with no
negative jumps, by showing that it is possible to obtain any
continuous-state branching process through a time change of the
corresponding L\'evy process; see also~\cite{cablam} for a very
readable account of this theory and proofs.

If we exclude processes that can make an instantaneous jump to
infinity, continuous-state branching processes can be characterized by
a function $\Psi\dvtx  [0, \infty) \rightarrow\R$ of the form
\[
\Psi(u) = \alpha u + \beta u^2 + \int_0^{\infty} \bigl(e^{-ux} - 1 + ux
\mathbf{1}_{\{x \leq1\}}\bigr) \nu(dx),
\]
where $\alpha\in\R$, $\beta\geq0$ and $\nu$ is a measure on $(0,
\infty)$ satisfying $\int_0^{\infty} (1 \wedge x^2) \nu(dx) <
\infty
$. The function $\Psi$ is called the branching mechanism. If $(Z(t), t
\geq0)$ is a continuous-state branching process with branching
mechanism $\Psi$, then for $\lambda\geq0$,
%
%e5 ###
%
%e5 #&#
\begin{equation}\label{lapcsbp}
E\bigl[e^{-\lambda Z(t)}|Z(0) = a\bigr] = e^{-a u_t(\lambda)},
\end{equation}
where the function $t \mapsto u_t(\lambda)$ is a solution to the
differential equation
%
%e6 ###
%
%e6 #&#
\begin{equation}\label{csbpdiffeq}
\frac{\partial}{\partial t} u_t(\lambda) = -\Psi(u_t(\lambda)),
u_0(\lambda) = \lambda.
\end{equation}

Neveu~\cite{neveu} studied the continuous-state branching process with
$\Psi(u) = u \log u$. We will be interested, more generally, in a
continuous-state branching process $(Z(t), t \geq0)$ whose branching
mechanism is of the form $\Psi(u) = au + b u \log u$, where $a \in\R$
and $b > 0$. In this case (see, e.g., page 256 of~\cite{beleg00}),
there exists a real number $c$ such that
\[
\Psi(u) = -cu + b \int_0^{\infty} \bigl(e^{-ux} - 1 + ux \mathbf{1}_{\{x
\leq1\}}\bigr) x^{-2} \, dx.
\]
Also, it is not difficult to solve (\ref{csbpdiffeq}) to obtain
%
%e7 ###
%
%e7 #&#
\begin{equation}\label{utlambda}
u_t(\lambda) = \lambda^{e^{-bt}} e^{a(e^{-bt} - 1)/b}.
\end{equation}
Because $\int_0^{\delta} 1/\Psi(u) \, du = \infty$ for all $\delta>
0$, the process does not explode. That is, almost surely $Z(t) < \infty
$ for all $t$. Because $\int_{\delta}^{\infty} 1/\Psi(u) = \infty$ for
all $\delta> 0$, the process does not go extinct, that is, almost
surely $Z(t) > 0$ for all $t$. Proofs of these facts can be found in
\cite{grey}.

%s1.3 ###
%s1.3 #&#
\subsection{The Bolthausen--Sznitman coalescent}\label{sec13}
\label{boszsec}

In mathematical population genetics, it is standard to represent the
ancestral relationships among a sample of $n$ individuals using a
coalescent process $(\Pi(t), t \geq0)$, which is a continuous-time
Markov process taking its values in the set of partitions of $\{
1,\ldots, n\}$. Here $\Pi(0)$ is the partition of $\{1,\ldots, n\}$
into $n$
singletons, and blocks of the partition merge over time. The merging of
blocks of the partition corresponds to the merging of ancestral lines
when the ancestral lines of the $n$ sampled individuals are traced
backwards in time. The standard coalescent model is Kingman's
coalescent. Kingman's coalecent was introduced in~\cite{king82} and is
now the basis for much work in mathematical population genetics.
Kingman's coalescent has the property that only two blocks of the
partition ever merge at a time, and each transition that involves two
blocks merging into one happens at rate one.

Within the last decade, alternative models of coalescence, allowing for
multiple ancestral lines to merge at once, have been studied in some
depth. These coalescent processes, known as coalescents with multiple
mergers or $\Lambda$-coalescents, were introduced by Pitman \cite
{pit99} and Sagitov~\cite{sag99}. If $\Lambda$ is a finite measure on
$[0, 1]$, then the $\Lambda$-coalescent has the property that whenever
there are $b$ blocks, each transition that involves merging $k$ blocks
of the partition into one happens at rate
\[
\lambda_{b,k} = \int_0^1 x^{k-2} (1 - x)^{b-k} \Lambda(dx).
\]
Kingman's coalescent is the special case of the $\Lambda$-coalescent in
which $\Lambda$ is the unit mass at zero.

If $\Lambda$ is the uniform distribution on $[0, 1]$, then the
$\Lambda
$-coalescent is known as the Bolthausen--Sznitman coalescent. The
Bolthausen--Sznitman coalescent was introduced in~\cite{bosz98} in the
context of Ruelle's probability cascades. The Bolthausen--Sznitman
coalescent has been studied extensively, and has been found to be
related to stable subordinators~\cite{bepi00} and random recursive
trees~\cite{goma05}. It also shows up in Derrida's generalized random
energy model~\cite{bovkur}. Properties of the Bolthausen--Sznitman
coalescent have been worked out, for example, in~\cite{pit99,dimr07,bago07}.

Bertoin and Le Gall~\cite{beleg00} showed how to define precisely the
notion of the genealogy of a continuous-state branching process. They
found that the genealogy of Neveu's continuous-state branching process
is given by the Bolthausen--Sznitman coalescent. These results were
extended in~\cite{bbcemsw}, where it was shown that the genealogy of
any continuous-state branching process whose branching mechanism is of
the form $\Psi(u) = au + bu \log u$ can still be described by the
Bolthausen--Sznitman coalescent. This connection between the
Bolthausen--Sznitman coalescent and Neveu's continuous-state branching
process played a central role in Bovier and Kurkova's analysis of
Derrida's generalized random energy model~\cite{bovkur}. A survey of
this material can be found in~\cite{ensaios}.

%s1.4 ###
%s1.4 #&#
\subsection{Main results}\label{sec14}

Recall that for each positive integer $N$, we have a branching Brownian
motion $(X_N(t), t \geq0)$. We denote by $M_N(t)$ the number of
particles at time $t$, and we denote the positions of these particles
by $X_{1,N}(t) \geq X_{2,N}(t) \geq\cdots\geq X_{M_N(t), N}(t)$.
We further define the process $(Z_N(t), t \geq0)$ by setting
%
%e8 ###
%
%e8 #&#
\begin{equation}\label{defL}
L = \tfrac{1}{\sqrt{2}} ( \log N + 3 \log\log N )
\end{equation}
and then letting
%
%e9 ###
%
%e9 #&#
\begin{equation}\label{ZNdef}
Z_N(t) = \sum_{i=1}^{M_N(t)} e^{\mu X_{i,N}(t)} \sin\biggl( \frac{\pi
X_{i,N}(t)}{L} \biggr) \mathbf{1}_{\{X_{i,N}(t) \leq L\}}.
\end{equation}
Note that only particles to the left of $L$ contribute to $Z_N(t)$, and
the level $L$ depends on $N$.
As we will see later, $Z_N(t)$ is a good measure of the ``size'' of the
process at time $t$, in the sense that it predicts the number of
particles shortly after time $t$. Also let
%
%e10 ###
%
%e10 #&#
\begin{equation}\label{YNdef}
Y_N(t) = \sum_{i=1}^{M_N(t)} e^{\mu X_{i,N}(t)}.
\end{equation}

We will see that as the branching Brownian motion evolves, most
particles stay well to the left of $L$, and as long as this is the
case, the number of particles changes little. However, occasionally a
small number of particles get very far to the right. Because the
descendants of these particles are able to avoid the barrier at zero,
the number of particles increases rapidly. Indeed, the increase in the
number of particles is so rapid that when we take the scaling limit as
$N \rightarrow\infty$, we get a process with jumps. The proposition
below shows that this limiting process is a continuous-state branching process.

%
%pr1 #&#
\begin{Prop}\label{ZNCSBP}
For all positive integers $N$, define the process $(V_N(t),\break t \geq0)$ by
%
%e11 ###
%
%e11 #&#
\begin{equation}\label{VNdef}
V_N(t) = \frac{1}{N (\log N)^2} Z_N ( (\log N)^3 t ).
\end{equation}
Suppose as $N \rightarrow\infty$, the distribution of $V_N(0)$
converges to $\nu$, where $\nu$ is a probability distribution on $[0,
\infty)$. Suppose also that $Y_N(0)/N(\log N)^3$ converges to zero in
probability as $N \rightarrow\infty$.
Then there exists a constant $a \in\R$ such that as $N \rightarrow
\infty$, the finite-dimensional distributions of the process $(V_N(t),
t \geq0)$ converge to the finite-dimensional distributions of the
continuous-state branching process with branching mechanism $\Psi(u) =
au + 2 \pi^2 u \log u$ started with distribution $\nu$ at time zero.
\end{Prop}

The condition on $V_N(0)$ ensures that the number of particles in the
system is of order $N$, as shown below with the scaling in Theorem \ref
{MNCSBP}. The condition on $Y_N(0)$ ensures that no single particle at
time $0$ is likely to have descendants that constitute a large fraction
of the population a short time later. If we begin with $N$ particles in
what is a relatively ``stable'' configuration, then the initial
conditions will hold. Furthermore, as shown in Proposition~3 of \cite
{bbs}, if there is initially a single particle near $L$, then these
conditions will be satisfied after a time of order $L^2$.

Note that because the processes $(V_N(t), t \geq0)$ for fixed $N$ can
increase very rapidly in a short time but do not have large jumps, the
sequence of processes $(V_N, N \geq1)$ is not tight, and convergence
in the Skorohod topology does not hold.

The theorem below converts this result about the scaling limit of
$(Z_N(t),\break t \geq0)$ to a result about the number of particles. This
convergence result holds only for $t > 0$. The hypothesis at time $t =
0$ still involves the processes $(V_N(t), t \geq0)$, which may not
imply convergence of the number of particles at time zero. The result
needs to be stated in this way because it is the value of $Z_N(t)$
rather than $M_N(t)$ that predicts the number of particles that will be
alive a short time later.

%
%th2 #&#
\begin{Theo}\label{MNCSBP}
Assume the hypotheses of Proposition~\ref{ZNCSBP} hold. Then as $N
\rightarrow\infty$, the finite-dimensional distributions of the process
\[
\biggl(\frac{1}{2 \pi N} M_N ((\log N)^3 t), t > 0 \biggr)
\]
converge to the finite-dimensional distributions of the
continuous-state branching process with branching mechanism $\Psi(u) =
au + 2 \pi^2 u \log u$ started with distribution $\nu$ at time zero,
where $a$ is the constant from Proposition~\ref{ZNCSBP}.
\end{Theo}

The next result shows that if we pick $n$ particles at random from
branching Brownian motion with absorption at some time and trace back
their ancestral lines, the resulting process, properly scaled,
converges to the Bolthausen--Sznitman coalescent. This is a precise
formulation of the result stated in the \hyperref[sec1]{Introduction}.
Choose $n$ particles uniformly at random from the $M_N((\log N)^3t)$
particles at time $(\log N)^3 t$, and label these particles at random
by the integers $1,\ldots, n$. Fix $t > 0$. For $0 \leq s \leq2 \pi
t$, define $\Pi_N(s)$ to be the partition of $\{1,\ldots, n\}$ such
that $i$ and $j$ are in the same block of $\Pi_N(s)$ if and only if the
particles labeled $i$ and $j$ are descended from the same ancestor at
time $(t - s/2 \pi) (\log N)^3$. Let $(\Pi(s), 0 \leq s \leq2 \pi t)$
be the Bolthausen--Sznitman coalescent run for time $2 \pi t$ and
restricted to $\{1,\ldots, n\}$.

%
%th3 #&#
\begin{Theo}\label{bosz}
Assume the hypotheses of Proposition~\ref{ZNCSBP} hold, and assume that
$\nu(\{0\}) = 0$. Then as $N \rightarrow\infty$, the
finite-dimensional distributions of $(\Pi_N(s), 0 \leq s \leq2 \pi t)$
converge to those of $(\Pi(s), 0 \leq s \leq2 \pi t)$.
\end{Theo}

As discussed earlier, this result is, of course, the analog for this
model of the third conjecture of Brunet et al.~\cite{bdmm1,bdmm2}
stated above. The $(\log N)^3$ time scaling that appears here, as well
as in Proposition~\ref{ZNCSBP} and Theorem~\ref{MNCSBP}, matches the
second conjecture stated above. If two particles are chosen at random,
the time back to their most recent common ancestor is of the order
$(\log N)^3$.

%s1.5 ###
%s1.5 #&#
\subsection{Overview of the proofs}\label{sec15}

Because the proofs of Proposition~\ref{ZNCSBP} and Theorems \ref
{MNCSBP} and~\ref{bosz} are rather long, we outline the basic strategy
here. The key idea is to treat separately the particles that reach
approximately the level $L$. These are the particles that will produce
a large number of descendants within a short time, leading to jumps in
the population size when we look forward in time, and multiple mergers
of ancestral lines going backwards in time.

The first step, carried out in Section~\ref{sec2}, is to collect some results
that we need pertaining to branching Brownian motion in a strip, which
are important both for the proofs and for understanding the heuristics
behind our choices of parameters. Most importantly, we observe that if
a branching Brownian motion is started with a single particle at $x$,
and particles are killed upon reaching $0$ or $L$, then the expected
number of particles in a set $B$ at a sufficiently large time $t$ is
approximately $\int_B p_t(x,y) \, dy$, where
%
%e12 ###
%
%e12 #&#
\begin{equation}\label{newpt}
p_t(x,y) = \frac{2}{L} e^{(1 - \mu^2/2 - \pi^2/2L^2)t} \cdot e^{\mu x}
\sin\biggl( \frac{\pi x}{L} \biggr) \cdot e^{-\mu y} \sin\biggl( \frac
{\pi y}{L} \biggr).
\end{equation}
From this formula, we can make several observations concerning the
behavior of the branching Brownian motion. First, note that the time
parameter $t$ appears in the formula only in the first exponential
factor, so the population size should be roughly constant over time
provided that $1 - \mu^2/2 - \pi^2/2L^2 = 0$. Indeed, we have chosen
the parameters $\mu$ and $L$ above [see (\ref{mu}) and (\ref{defL})] to
satisfy this equation, as this is the drift needed to stabilize the
population size. Second, notice that the formula is proportional to
$e^{\mu x} \sin(\pi x/L)$, which will equal $Z_N(t)$\vadjust{\goodbreak} if we sum over the
positions of all particles at time $t$. Thus, it is $Z_N(t)$ that
predicts the number of particles that will be in a given set at a later
time, which is why $Z_N(t)$ provides a useful measure of the ``size'' of
the process. Third, notice that the formula is proportional to $e^{-\mu
y} \sin(\pi y/L)$. Consequently, regardless of the starting
configuration, once $t$ is large enough for the approximation to be
valid, the particles will have settled into a ``stable'' configuration
in which the ``density'' of particles at position $y$ is proportional to
$e^{-\mu y} \sin(\pi y/L)$. We will see in Lemma~\ref{stripdensity}
that this approximation becomes accurate when $t$ gets to be larger
than $(\log N)^2$.

If we begin at time zero with $N$ particles that are approximately in
the stable configuration, so that their ``density'' is $C L e^{-\mu y}
\sin(\pi y/L)$, where $CL$ is a normalizing constant, then the value of
$Z_N(0)$ should be approximately
\[
N \int_0^L e^{\mu y} \sin\biggl( \frac{\pi y}{L} \biggr) \cdot C L
e^{-\mu y} \sin\biggl( \frac{\pi y}{L} \biggr) \, dy,
\]
which is of the order $NL^2$. On the other hand, if we begin instead
with a single particle at $L$, then one can show typically the
right-most descendant of this particle will reach a level that exceeds
$L$ by only a constant. This is essentially true because critical
branching Brownian motion dies out, and can be seen from Proposition
\ref{EMkProp} below which shows that particles reach $L$ at a much
faster rate than they reach any level that is much greater than $L$.
Consequently, we can estimate the typical contribution of the
descendants of this particle at time $t$ by using (\ref{newpt}) with
$L$ in place of $x$ and $L + \alpha$ in place of $L$, where $\alpha>
0$ is a constant. This means that the value of $Z_N(t)$ should be of
the same order as
\[
\int_0^L e^{\mu y} \sin\biggl( \frac{\pi y}{L} \biggr) \cdot\frac{2}{L
+ \alpha} e^{\mu L} \sin\biggl( \frac{\pi L}{L + \alpha} \biggr) e^{-\mu
y} \sin\biggl( \frac{\pi y}{L+\alpha} \biggr) \, dy,
\]
which is of the order $L^{-1} e^{\mu L}$. We have chosen $L$ so that
particles that reach $L$ produce substantial increases in the
population size. Indeed, note that $L^{-1} e^{\mu L}$ and $N L^2$ are
of the same order precisely when $L$ is within a constant of the value
in (\ref{defL}).

In Section~\ref{sec3}, we therefore define
%
%e13 ###
%
%e13 #&#
\begin{equation}\label{LAdef}
L_A = \tfrac{1}{\sqrt{2}} ( \log N + 3 \log\log N - A ),
\end{equation}
where $A \in\R$, and study the particles that stay to the left of
$L_A$. That is, we consider branching Brownian motion with particles
killed at $0$ and at $L_A$.
Using (\ref{newpt}), it possible to estimate first and second moments
of various quantities. In Section~\ref{sec3}, we apply these results
to calculate the first and second moments of $Z_N(t)$, conditional on
the process a time $\theta(\log N)^3$ earlier, where $\theta$ is a
small constant. The first moment calculation is Lemma~\ref{ZN1exp},
while the variance bound appears in Lemma~\ref{VarZ1}. The variance
bound is sufficient to establish that when $A$ is large, there is a law
of large numbers, with the value of $Z_N(t)$\vadjust{\goodbreak} being close to its
expectation. A similar variance bound for the number of particles is
given in Lemma~\ref{MNVar}. Such results would not be possible without
the truncation at $L_A$, because without truncation the expected number
of particles is dominated by rare events in which one particle moves
far to the right and produces a large number of surviving offspring.
The analysis in Section~\ref{sec3} is motivated by some of the arguments based
on moment bounds in~\cite{kesten}.

In  Section~\ref{sec32} we tackle the question of how many particles
reach the level $L_A$. An estimate of the expected number is given in
Proposition~\ref{EMkProp}. From this result, one can deduce that if we
start with $N$ particles that are in approximately the ``stable''
configuration described above, then the time that it will take before a
particle reaches $L_A$ is of the order $(\log N)^3$, which explains the
$(\log N)^3$ time scaling in our main results. To see heuristically why
this scaling occurs, note that if $\beta> 0$ is a constant, then the
number of particles between $L - \beta$ and $L$ at time $t$ is of the
order
\[
N \int_{L - \beta}^L C L e^{-\mu y} \sin\biggl( \frac{\pi y}{L} \biggr)
\, dy,
\]
which is of the order $1/(\log N)^3$. Such particles have a positive
probability of reaching $L$ between times $t$ and $t + 1$, but the
calculation in Proposition~\ref{EMkProp} shows that particles that are
more than a constant distance from $L$ at time $t$ are unlikely to hit
$L$ by time $t + 1$. Thus, $O(1/(\log N)^3)$ particles hit $L$ per unit time.

Since branching by particles close to $L_A$ may enable several
particles to hit $L_A$ at nearly the same time, we also require the
second moment estimate in Proposition~\ref{EMk2Prop} to establish that
the expected number of particles that reach $L_A$ within a time
interval of length $\theta(\log N)^3$, conditional on at least one
particle reaching $L_A$, is bounded by a constant. Then in Section
\ref{sec33}, we show in Proposition~\ref{GNProp} that a ``good'' event on
which the bounds in  Sections~\ref{sec31} and~\ref{sec32} are valid
occurs with high probability.

In Section~\ref{sec4}, we begin to consider the contribution from
particles after they reach the level $L_A$. The key to this analysis is
Proposition~\ref{nevprop}, which comes from~\cite{nev87}. This result
states that if a particle starts at $L_A$, and $y$ is a large constant,
then the number of descendants of the particle that reach $L_A - y$ is
approximately $y^{-1} e^{\sqrt{2} y} W$, where $W$ is a random
variable. Some analysis that involves a Tauberian theorem leads to
Proposition~\ref{Wprop}, which says that for large $x$, we have $P(W >
x) \sim B/x$. Conceptually, this result is the reason why the genealogy
of the population is described by the Bolthausen--Sznitman coalescent.
The contribution to the population of the particle at $L_A$ will be
approximately proportional to the number of descendants that hit $y$,
if $y$ is sufficiently large. The fact that a jump of size greater than
$x$ results from a particle at $L_A$ with probability proportional to
$1/x$ implies that the L\'evy measure of the limiting continuous-state
branching process will have a density proportional to $x^{-2}$, which
in turn leads to the duality with the Bolthausen--Sznitman coalescent.

In Section~\ref{sec5}, we show how to combine all of the previous
estimates to get sharp results for the behavior of the process
$(Z_N(t), t \geq0)$. The key results are Proposition \ref
{smalljumpexp}, which bounds the expected change in $Z_N$ over a time
interval of length $\theta(\log N)^3$ when there is no large jump, and
Proposition~\ref{rjumplem}, which estimates the probability that $Z_N$
increases by at least $r N (\log N)^2$ over a time interval of length
$\theta(\log N)^3$. These estimates on how the process behaves over a
short time interval can be matched with the infinitesimal generator of
the continuous-state branching process. This work is done in Section
\ref{sec6} and leads to a proof of Proposition~\ref{ZNCSBP}. Once
Proposition~\ref{ZNCSBP} is established, we are able to prove Theorem
\ref{MNCSBP} by arguing that the value of $Z_N(t)$ can be used to
predict accurately the number of particles shortly after time~$t$.

The proof that the genealogy of the process converges to the
Bolthausen--Sznitman coalescent is completed in Section~\ref{sec7}. We
represent the genealogy of the branching Brownian motion using a ``flow
of bridges,'' a tool introduced by Bertoin and Le Gall in~\cite{beleg1}.
Using Proposition~\ref{ZNCSBP} and Theorem~\ref{MNCSBP}, we establish
convergence to the flow of bridges associated with the continuous-state
branching process, which is known to correspond to the
Bolthausen--Sznitman coalescent.

%s1.6 ###
%s1.6 #&#
\subsection{Notational conventions and index of notation}\label{sec16}

For the benefit of the reader, we include in   Table~\ref{tab1} an index of some of the notation.

%t1 #&#
\begin{table}
\tabcolsep=0pt
\caption{Index
of some of the notation that is used throughout the paper}\label{tab1}
\begin{tabular*}{\textwidth}{@{\extracolsep{\fill}}ll@{}}
\hline
$A$ & used to control the level at which particles are killed; see the
definition of $L_A$. \\
$G_{N,k}$ & event that $Z_N(t_j)$ and $Y_N(t_j)$ are sufficiently small
for $j \leq k$. \\
$G_N(\eps)$ & event that $G_{N,k}$ occurs for all $k$. \\
$h(N)$ & slowly increasing function used to upper bound $Y_N$. \\
$L$ & level, given by (\ref{defL}), such that descendants of a particle
that get near this level \\
& will likely constitute a significant fraction of the population in
the future. \\
$L_A$ & level at which particles are killed, defined in (\ref{LAdef}).
\\
$M_N(t)$ & number of particles at time $t$. \\
$R_k$ & number of particles killed at $L_A$ between $t_{k-1}$ and
$t_k$. \\
$s$ & the process $Z_N$ is often studied between times $u (\log N)^3$
and $(u+s)(\log N)^3$. \\
$t_k$ & the process $Z_N$ is frequently studied at the times $t_k$. \\
$u$ & the process $Z_N$ is often studied between times $u (\log N)^3$
and $(u+s)(\log N)^3$. \\
$V_N$ & normalization of the process $Z_N$, defined in (\ref{VNdef}).
\\
$X_N(t)$ & the branching Brownian motion at time $t$. \\
$X_{i,N}(t)$ & position of the $i$th particle from the right at time
$t$. \\
$y$ & large constant; the number of descendants of a particle at $L_A$
that reach \\
& $L_A - y$ plays a central role in the paper. \\
$Y_N(t)$ & weighted sum of particle positions at time $t$, defined in
(\ref{YNdef}), such that a \\
& particle at $x$ contributes $e^{\mu x}$ to the sum. \\
$Z_N(t)$ & measure of the ``size'' of the process at time $t$, defined
in (\ref{ZNdef}), such that a \\
& particle at $x \leq L$ contributes $e^{\mu x} \sin(\pi x/L)$. \\
$Z_{N,1}$ & similar to $Z_N$, but with particles killed at $L_A$,
defined in (\ref{ZN1def}). \\
$Z_{N,1}'$ & similar to $Z_{N,1}$, with $L_A$ used in place of $L$ in
the sine function; see (\ref{ZNprime}). \\
$Z_y$ & number of descendants of a particle at zero that reach $-y$. \\
$\delta$ & small constant used to bound the error in an estimate of a
branching process \\
& limit; see (\ref{WBbound}). \\
$\eps$ & small constant used to bound $Z_N$ above by $\eps^{-1/2} N
(\log N)^2$. \\
$\mu$ & drift of the branching Brownian motion, given by (\ref{mu}).
\\
$\eta$ & small constant used to bound the difference between $Z_y$ and
its limit. \\
$\theta$ & small constant such that $t_k$ and $t_{k+1}$ are $\theta s
(\log N)^3$ apart. \\
$\zeta$ & large constant chosen so that with high probability,
descendants of a particle \\
& at zero will have reached $-y$ by time $\zeta$.\\
\hline
\end{tabular*}
\end{table}

Some constraints on the constants $\eps$, $A$ and $\theta$ are
introduced at the beginning of Section~\ref{sec3}; see equations (\ref
{con0})--(\ref{con3}). Further constraints on these constants, as well
as the choices of the constants $\delta$, $\eta$, $y$ and $\zeta$, are
set out in Section~\ref{sec51}; see equations (\ref{deltasmall})--(\ref{659}).

Throughout the rest of the paper, $C$ will denote a positive finite
constant whose value may change from line to line. The constant $C$ may
depend on $u$ and $s$, but may not depend on $N$ or on the seven
constants $\eps$, $A$, $\theta$, $\delta$, $\eta$, $y$ and $\zeta
$. We
say a sequence of random variables $(R_N)_{N=1}^{\infty}$ is $o(f(N))$
if for any choices of the constants $u$, $s$, $\eps$, $A$, $\theta$,
$\delta$, $\eta$, $y$ and $\zeta$ satisfying the constraints mentioned
above, there is a deterministic sequence $(b_N)_{N=1}^{\infty}$ tending
to zero such that $|R_N| \leq b_N f(N)$ for all $N$. Note in particular
that throughout this paper, the bounds implicit in the notation $o(1)$
or $o(f(N))$ are nonrandom and depend solely on the choices of parameters.

Also, if $g$ is a function of some of the constants $\eps$, $A$,
$\theta
$, $\delta$, $\eta$, $y$, $\zeta$ and $N$, we will occasionally use the
notation $O(g(\eps, A, \theta, \delta, \eta, y, \zeta, N))$ to denote
an expression whose absolute value is bounded by $C g(\eps, A, \theta,
\delta, \eta, y, \zeta, N)$, where $C$ is defined as above.

%s2 ###
%s2 #&#
\section{Branching Brownian motion in a strip}\label{sec2}

Suppose $(B_t)_{t \geq0}$ is Brownian motion started at $x$, with $0 <
x < K$, and assume the process is killed when it hits $0$ or $K$. Then\vadjust{\goodbreak}
(see, e.g., page 188 of~\cite{lawler}) the density of the process at
time $t$, restricted to $(0, K)$, is
%
%e14 ###
%
%e14 #&#
\begin{eqnarray}\label{vtdef}
v_t(x,y) &=& \frac{\pi}{K} u_{\pi^2 t/K^2}(\pi x/K, \pi y/K)\nonumber
\\[-8pt]
\\[-8pt] &=& \frac
{2}{K} \sum_{n=1}^{\infty} e^{-\pi^2 n^2 t/2K^2} \sin\biggl( \frac{n
\pi x}{K} \biggr) \sin\biggl( \frac{n \pi y}{K} \biggr).
\nonumber
\end{eqnarray}

Consider now branching Brownian motion in a strip in which each
particle gives birth at rate one, drifts to the left at rate $\mu> 0$,
and is killed upon reaching $0$ or~$K$. We will need to estimate the
expected number of particles at time $t$ when $t$ is large. Suppose
there is initially a single particle at $x$. The density of particles
at the position $y$ at time $t$ can be calculated using the well-known
many-to-one lemma. The density is a product of $e^t$, which represents
the expected number of particles at time $t$,
a Girsanov factor $e^{\mu(x-y) - \mu^2t/2}$ relating Brownian motion
with drift $-\mu$ to ordinary Brownian motion, and the density of
ordinary Brownian motion killed upon reaching $0$ or $K$. Therefore,
the density of particles at time $t$ is
%
%e15 ###
%
%e15 #&#
\begin{equation}\label{density}
 \qquad q_t(x, y) = e^{(1 - \mu^2/2)t + \mu(x-y)} \cdot\frac{2}{K} \sum
_{n=1}^{\infty} e^{-\pi^2 n^2 t/2K^2} \sin\biggl( \frac{n \pi x}{K}
\biggr) \sin\biggl( \frac{n \pi y}{K} \biggr),
\end{equation}
in the sense that if $B \subset(0, K)$, then the expected number of
particles in $B$ at time $t$ is $\int_B q_t(x,y) \, dy$.

When $t \gg K^2$, the first term in the sum in (\ref{density})
dominates. We make this more precise in Lemma~\ref{stripdensity} below.
We first record the following trigonometric lemma.

%
%le4 #&#
\begin{Lemma}\label{trig}
If $0 \leq y \leq\pi$ and $n \in\N$, then $|\sin ny| \leq n \sin y$.
\end{Lemma}

\begin{pf}
We prove the result by induction. The result is trivial for $n = 1$. If
it is true for $n-1$, then
\begin{eqnarray*}
|\sin ny| &=& \bigl|\sin\bigl((n-1) y\bigr) \cos y + \cos\bigl((n-1)y\bigr) \sin y\bigr| \\
&\leq&\bigl|\sin\bigl((n-1)y\bigr)\bigr| |\cos y| + \bigl|\cos\bigl((n-1)y\bigr)\bigr| |\sin y| \\
&\leq&\bigl|\sin\bigl((n-1)y\bigr)\bigr| + |\sin y| \leq n \sin y,
\end{eqnarray*}
where the last step uses the induction hypothesis.
\end{pf}

By applying Lemma~\ref{trig} to each term in the sum on the right-hand
side of (\ref{density}), we easily get the following estimate. Note
that $p_t(x,y)$ is simply the $n = 1$ term in the expression for
$q_t(x, y)$. The error term $D_t(x,y)$ is small when $t \gg K^2$ and is
bounded above by a constant when $t \geq C_1 K^2$ for some constant $C_1$.

%
%le5 #&#
\begin{Lemma}\label{stripdensity}
Consider branching Brownian motion in a strip in which each particle
gives birth at rate one, drifts to the left at rate $\mu$ and is killed
upon reaching $0$ or $K$. Suppose there is initially a single particle
at $x$. Let
\[
p_t(x,y) = \frac{2}{K} e^{(1 - \mu^2/2 - \pi^2/2K^2)t} \cdot e^{\mu x}
\sin\biggl( \frac{\pi x}{K} \biggr) \cdot e^{-\mu y} \sin\biggl( \frac
{\pi y}{K} \biggr).
\]
Then for all $x, y \in[0, K]$, define $D_t(x,y)$ by
\[
\frac{q_t(x,y)}{p_t(x,y)} = 1 + D_t(x,y).
\]
Then
%
%e16 ###
%
%e16 #&#
\begin{equation}\label{Eterm}
|D_t(x,y)| \leq\frac{\sum_{n=2}^{\infty} n^2 e^{-\pi^2 n^2 t/2
K^2}}{e^{-\pi^2 t/2K^2}}.\vadjust{\goodbreak}
\end{equation}
Therefore, if $B$ is a Borel subset of $(0, K)$, then the expected
number of particles in $B$ at time $t$ may be written as $(\int_B
p_t(x,y) \, dy)(1 + D'_t(x,B))$, where $|D'_t(x,B)|$ is bounded by the
right-hand side of (\ref{Eterm}).
\end{Lemma}

Using these densities, we can estimate the expected values of certain
functions of branching Brownian motion. Lemma~\ref{hhlem}, which is
Lemma 2 of~\cite{hh07}, gives a martingale for branching Brownian
motion in which particles are killed only at zero. Lemma~\ref{YZexplem}
estimates the expected values of three specific functions of branching
Brownian motion in a strip.

%
%le6 #&#
\begin{Lemma}\label{hhlem}
Consider branching Brownian motion in which each particle gives birth
at rate one, drifts to the left at rate $\mu$ and is killed upon
reaching $0$.
Let $M(t)$ be the number of particles at time $t$, and denote the
positions of the particles at time $t$ by $X_1(t),\ldots, X_{M(t)}(t)$.
Let
\[
V(t) = \sum_{i=1}^{M(t)} X_i(t) e^{\mu X_i(t) + (\mu^2/2 - 1) t}.
\]
Then $(V(t), t \geq0)$ is a martingale.
\end{Lemma}

%
%le7 #&#
\begin{Lemma}\label{YZexplem}
Consider branching Brownian motion in a strip in which each particle
gives birth at rate one, drifts to the left at rate $\mu$ and is killed
upon reaching $0$ or $K$.
Let $M(t)$ be the number of particles at time $t$, and denote the
positions of the particles at time $t$ by $X_1(t),\ldots, X_{M(t)}(t)$.
Let
\[
Y(t) = \sum_{i=1}^{M(t)} e^{\mu X_i(t)}, \qquad  Z(t) = \sum
_{i=1}^{M(t)} e^{\mu X_i(t)} \sin\biggl( \frac{\pi X_i(t)}{K} \biggr).
\]
Then
%
%e17 ###
%
%e17 #&#
\begin{equation}\label{Mexp}\qquad
E[M(t)] = \frac{2}{K} e^{(1 - \mu^2/2 - \pi^2/2K^2) t} (1 + D_1) Z(0)
\int_0^K e^{-\mu y} \sin\biggl( \frac{\pi y}{K} \biggr) \, dy
\end{equation}
and
%
%e18 ###
%
%e18 #&#
\begin{equation}\label{Yexp}
E[Y(t)] = \frac{4}{\pi} e^{(1 - \mu^2/2 - \pi^2/2K^2) t} (1 + D_2) Z(0),
\end{equation}
where $|D_1|$ and $|D_2|$ are bounded by the right-hand side of (\ref
{Eterm}). Also,
%
%e19 ###
%
%e19 #&#
\begin{equation}\label{Zexp}
E[Z(t)] = e^{(1 - \mu^2/2 - \pi^2/2K^2) t} Z(0).
\end{equation}
\end{Lemma}

\begin{pf}
To prove (\ref{Mexp}), first suppose there is initially a single
particle at $x$. Lemma~\ref{stripdensity} gives
\begin{eqnarray*}
E[M(t)] &=&\biggl ( \int_0^K p_t(x,y) \, dy \biggr)(1 + D_1) \\
&=& \frac{2}{K} e^{(1 - \mu^2/2 - \pi^2/2K^2) t} (1 + D_1) e^{\mu x}
\sin\biggl( \frac{\pi x}{K} \biggr) \int_0^K e^{-\mu y} \sin\biggl( \frac
{\pi y}{K} \biggr) \, dy,
\end{eqnarray*}
where $|D_1|$ is bounded by the right-hand side of (\ref{Eterm}). The
result now follows by summing over the particles at time zero.

Likewise, to prove (\ref{Yexp}), assume there is initially a single
particle at $x$, and observe that Lemma~\ref{stripdensity} gives
\[
E[Y(t)] = \biggl( \int_0^K e^{\mu y} p_t(x,y) \, dy \biggr) (1 + D_2),
\]
where $|D_2|$ is bounded by the right-hand side of (\ref{Eterm}). Using
\[
\int_0^K \sin\biggl( \frac{\pi y}{K} \biggr) \, dy = \frac{2K}{\pi},
\]
we get
\[
E[Y(t)] = \frac{4}{\pi} e^{(1 - \mu^2/2 - \pi^2/2K^2) t} e^{\mu x}
\sin
\biggl( \frac{\pi x}{K} \biggr) (1 + D_2).
\]
The result again follows by summing over the particles at time zero.

To obtain (\ref{Zexp}), note that if $n$ is a positive integer, then
\[
\int_0^K \sin\biggl( \frac{\pi y}{K} \biggr) \sin\biggl(\frac{n \pi
y}{K} \biggr) \, dy =
\cases{\displaystyle K/2 ,&\quad if $n = 1 $,\cr\displaystyle
0 ,&\quad if $n \geq2$.
}
\]
If at time zero there is just a single particle at $x$, then
\begin{eqnarray*}
&&\hspace*{-5pt}E[Z(t)]\\
 && \hspace*{-5pt}\quad = \int_0^K e^{\mu y} \sin\biggl( \frac{\pi y}{K} \biggr)
q_t(x,y) \, dy \\
&& \hspace*{-5pt}\quad = e^{(1 - \mu^2/2)t + \mu x} \cdot\frac{2}{K} \sum_{n=1}^{\infty}
e^{-\pi^2 n^2t/2K^2} \sin\biggl(\frac{n \pi x}{K} \biggr) \int_0^K \sin
\biggl( \frac{\pi y}{K} \biggr) \sin\biggl( \frac{n \pi y}{K} \biggr) \,
dy \\
&& \hspace*{-5pt}\quad = e^{\mu x} \sin\biggl( \frac{\pi x}{K} \biggr) e^{(1 - \mu^2/2)t}
e^{-\pi^2 t/2K^2} = e^{(1 - \mu^2/2 - \pi^2/2K^2) t} e^{\mu x} \sin
\biggl( \frac{\pi x}{K} \biggr).
\end{eqnarray*}
As before, the result now follows by summing over the particles at time zero.
\end{pf}

For the next result, we will need the Green's function for Brownian motion
in a strip. Let $(B_t, t \geq0)$ be one-dimensional Brownian motion
without drift. Define the Green's function $G(x,y)$ such that if $(B_t, t
\geq0)$ is Brownian motion started from $B_0 = x \in(0, K)$ and if
$\tau= \inf\{t\dvtx  B_t \notin(0,K)\}$, then for all bounded measurable
functions $g$, we have
\[
E \biggl[ \int_0^{\tau} g(B_t) \, dt \biggr] = \int_0^K G(x,y) g(y) \, dy.
\]
The Green's function is given by (see, e.g., (4.4) on page 225 of~\cite{durr96})
%
%e20 ###
%
%e20 #&#
\begin{equation}\label{Green}
G(x,y) =
\cases{\displaystyle2x(K-y)/K ,&\quad if $y \geq x $,\cr
\displaystyle
2y(K-x)/K ,&\quad if $y \leq x$.
}
\end{equation}
To obtain this result from (4.4) in~\cite{durr96}, observe that in the
notation of~\cite{durr96}, we have $\varphi(x) = x$ and $m(x) = 1$ for
ordinary Brownian motion. If $y \leq x$, then $2y(K-x)/K \leq
2x(K-y)/K$. Therefore, for all $x, y \in[0, K]$,
%
%e21 ###
%
%e21 #&#
\begin{equation}\label{Gbound}
G(x,y) \leq2x(K-y)/K.
\end{equation}

To control the fluctuations, we will also need a result about second
moments. The following result, which is a slight extension of Lemma 3.1
of~\cite{kesten}, will be a useful tool.

%
%le8 #&#
\begin{Lemma}\label{kestenlem}
Consider branching Brownian motion with particles killed at both $0$
and $K$. Assume that at time zero there is just a single particle at
$x$, and that the particles at time $t$ are denoted by $X_1(t),\ldots,
X_{M(t)}(t)$. Let $f\dvtx  (0, K) \rightarrow[0, \infty)$ be a measurable
function. Then
\begin{eqnarray*}
E \Biggl[ \Biggl( \sum_{i=1}^{M(t)} f(X_i(t)) \Biggr)^2 \Biggr] &=& \int_0^K
f(y)^2 q_t(x,y) \, dy \\
&&{}+ 2 \int_0^t \int_0^K q_s(x,z) \biggl( \int_0^K
f(y) q_{t-s}(z,y) \, dy \biggr)^2 \, dz \, ds.
\end{eqnarray*}
\end{Lemma}

\begin{pf}
For a Borel set $A \subset(0, K)$, let $N_A(t)$ be the number of
particles in the set $A$ at time $t$. Equation (2.8) of~\cite{saw76} gives
%
%e22 ###
%
%e22 #&#
\begin{equation}\label{ENA}
E[N_A(t)] = \int_A q_t(x,y) \, dy,
\end{equation}
while equations (2.11) and (2.12) of~\cite{saw76} give
%
%e23 ###
%
%e23 #&#
\begin{eqnarray}\label{ENANB}\qquad
E[N_A(t)N_B(t)] &=& E[N_{A \cap B}(t)]+ 2 \int_0^t \int_0^K q_s(x,z)
\biggl( \int_A q_{t-s}(z,w) \, dw \biggr) \nonumber
\\[-8pt]
\\[-8pt]
&&\hphantom{E[N_{A \cap B}(t)]+ 2 \int_0^t \int_0^K}{} \times\biggl( \int_B q_{t-s}(z,y) \,
dy \biggr) \, dz \, ds.
\nonumber
\end{eqnarray}

Suppose $f$ is a simple function, so that
\[
f(x) = \sum_{i=1}^m a_i \mathbf{1}_{A_i},
\]
where the $A_i$ are disjoint Borel subsets of $(0, K)$ and the $a_i$
are positive real numbers. In this case, we have
\[
E \Biggl[ \Biggl( \sum_{i=1}^{M(t)} f(X_i(t)) \Biggr)^2 \Biggr] = \sum
_{i=1}^m \sum_{j=1}^m a_i a_j E[N_{A_i}(t) N_{A_j}(t)].
\]
It is now straightforward to check, using (\ref{ENA}) and (\ref
{ENANB}), that the conclusion of Lemma~\ref{kestenlem} holds in this
case. Since every nonnegative measurable function can be approximated
from below by simple functions, the general result then follows from
the monotone convergence theorem.
\end{pf}

%
%le9 #&#
\begin{Lemma}\label{varZlem}
Assume we are in the setting of Lemma~\ref{YZexplem}. Assume that at
time zero there is just a single particle at $x$. Suppose\vspace*{1pt} that $1 - \mu
^2/2 - \pi^2/2K^2 \leq0$. Also, assume there exist positive constants
$C_1$ and $C_2$ such that $C_1 K^2 \leq t \leq C_2/(1 - \mu^2/2)$. Then
there exists a constant $C$, depending on $\mu$, $C_1$ and $C_2$, but
not on $x$ or $K$, such that
\[
E[Z(t)^2] \leq C e^{\mu x} e^{\mu K} \biggl( \frac{1}{K^2} + \frac
{t}{K^4} \biggr).
\]
\end{Lemma}

\begin{pf}
We apply Lemma~\ref{kestenlem} with $f(y) = e^{\mu y} \sin(\pi y/K)$
to get
%
%e24 ###
%
%e24 #&#
\begin{eqnarray}\label{Z2terms}\qquad
E[Z(t)^2] &=& \int_0^K e^{2 \mu y} \sin\biggl( \frac{\pi y}{K} \biggr)^2
q_t(x,y) \, dy \nonumber
\\[-8pt]
\\[-8pt]
&&{} + 2 \int_0^t \int_0^K q_s(x,z)\biggl ( \int_0^K e^{\mu y}
\sin\biggl( \frac{\pi y}{K} \biggr) q_{t-s}(z,y) \, dy \biggr)^2 \, dz \, ds.
\nonumber
\end{eqnarray}

We begin by bounding the first term in (\ref{Z2terms}). By Lemma \ref
{stripdensity}, for all $x, y \in[0, K]$ we have
%
%e25 ###
%
%e25 #&#
\begin{equation}\label{qtbound}
q_t(x,y) \leq\frac{C}{K} e^{\mu(x-y)} \sin\biggl( \frac{\pi x}{K}
\biggr) \sin\biggl( \frac{\pi y}{K} \biggr),
\end{equation}
where we are using that $1 - \mu^2/2 - \pi^2/2K^2 \leq0$. The
assumption $t \geq C_1 K^2$ ensures that the error term from Lemma \ref
{stripdensity} can be bounded by a constant (throughout the proof, we
allow the value of $C$ to change from line to line). Note that
%
%e26 ###
%
%e26 #&#
\begin{eqnarray}\label{sineq}
\int_0^K e^{\mu y} \sin\biggl( \frac{\pi y}{K} \biggr) \, dy &= &\int_0^K
e^{\mu(K - y)} \sin\biggl( \frac{\pi(K-y)}{K} \biggr) \, dy \nonumber\\
&=& e^{\mu K} \int_0^K e^{-\mu y} \sin\biggl( \frac{\pi y}{K} \biggr) \,
dy \\
&\leq& e^{\mu K} \int_0^K e^{-\mu y} \biggl( \frac{\pi y}{K} \biggr) \, dy
\leq\frac{C e^{\mu K}}{K}.\nonumber
\end{eqnarray}
Here we are using that $\mu> 0$ and that $C$ may depend on $\mu$. Now
using (\ref{qtbound}) and (\ref{sineq}) and the bound $\sin(\pi y/K)^2
\leq1$, we get
%
%e27 ###
%
%e27 #&#
\begin{eqnarray}\label{term1}
   \int_0^K e^{2 \mu y} \sin\biggl( \frac{\pi y}{K} \biggr)^2 q_t(x,y) \,
dy &\leq&\frac{C}{K} \int_0^K e^{2 \mu y} e^{\mu(x-y)} \sin\biggl( \frac
{\pi x}{K} \biggr) \sin\biggl( \frac{\pi y}{K} \biggr) \, dy\nonumber
\\
  &\leq&\frac{C e^{\mu x}}{K} \int_0^K e^{\mu y} \sin\biggl( \frac{\pi
y}{K} \biggr) \, dy \\
&\leq&\frac{C e^{\mu x} e^{\mu K}}{K^2}.
\nonumber
\end{eqnarray}

It remains to bound the second term in (\ref{Z2terms}). Recall that
$v_t(x,y)$, defined in (\ref{vtdef}), denotes the density at time $t$
of Brownian motion started at $x$ and killed when it reaches $0$ or
$K$. Note that
\[
\int_0^{\infty} v_s(x,y) \, ds = G(x,y),
\]
where $G(x,y)$ is Green's function in (\ref{Green}). Since $t \leq
C_2/(1 - \mu^2/2)$, we also have for $s \leq t$,
%
%e28 ###
%
%e28 #&#
\begin{equation}\label{qsvs}
q_s(x,y) = e^{\mu(x-y) + (1 - \mu^2/2)s} v_s(x,y) \leq C e^{\mu
(x-y)} v_s(x,y).
\end{equation}
Since $t \geq C_1K^2$, the bound (\ref{qtbound}) is valid for
$q_{t-s}(x,y)$ when $s \leq t/2$. Using these results and (\ref{Gbound}),
%
%e29 ###
%
%e29 #&#
\begin{eqnarray} \label{term2}
&&\int_0^{t/2} \int_0^K q_s(x,z)\biggl ( \int_0^K e^{\mu y} \sin\biggl(
\frac{\pi y}{K} \biggr) q_{t-s}(z,y) \, dy \biggr)^2 \, dz \, ds
\nonumber\\
&& \qquad\leq\int_0^{t/2} \int_0^K C e^{\mu(x-z)} v_s(x,z) \biggl( \int_0^K
e^{\mu y} \sin\biggl( \frac{\pi y}{K} \biggr) \cdot\frac{C}{K} e^{\mu
(z-y)}\nonumber \\
&&\hphantom{\leq\int_0^{t/2} \int_0^K C e^{\mu(x-z)} v_s(x,z) \biggl( \int_0^K\,} \qquad {}\times\sin\biggl( \frac{\pi z}{K} \biggr) \sin\biggl( \frac{\pi y}{K}
\biggr) \, dy \biggr)^2 \, dz \, ds \nonumber\\
&& \qquad\leq\frac{C e^{\mu x}}{K^2} \int_0^{t/2} \int_0^K e^{\mu
z} v_s(x,z)
\sin\biggl( \frac{\pi z}{K} \biggr)^2 \biggl( \int_0^K \sin\biggl( \frac
{\pi y}{K} \biggr)^2 \, dy \biggr)^2 \, dz \, ds \\
&& \qquad\leq C e^{\mu x} \int_0^K e^{\mu z} \sin\biggl( \frac{\pi z}{K}
\biggr)^2 \biggl( \int_0^{t/2} v_s(x,z) \, ds \biggr) \, dz \nonumber\\
&& \qquad\leq C e^{\mu x} \int_0^K e^{\mu z} \sin\biggl( \frac{\pi z}{K}
\biggr)^2 \frac{2x (K-z)}{K} \, dz \nonumber\\
&& \qquad\leq C e^{\mu x} \int_0^K e^{\mu z} \frac{(K-z)^3}{K^2} \,
dz \leq
\frac{C e^{\mu x} e^{\mu K}}{K^2},
\nonumber
\end{eqnarray}
where for the third inequality, we used that $\sin(\pi y/K)^2 \leq1$,
and for the next-to-last inequality, we used that $\sin(\pi z/K) =
\sin
(\pi(K-z)/K) \leq(K-z)/K$ and $x/K \leq1$.

Next, let $v_t'(x,y)$ be the density at time $t$ of Brownian motion
started at $x$ and killed when it hits $0$. By the Reflection
Principle, for $s \leq t$,
\begin{eqnarray*}
\int_0^K y v'_s(x,y) \, dy &=& \frac{1}{\sqrt{2 \pi s}} \int_0^K \bigl( y
e^{-(x-y)^2/2s} - y e^{-(x+y)^2/2s} \bigr) \, dy \\
&= &\frac{1}{\sqrt{2 \pi s}} \int_{-x}^{K-x} (z + x) e^{-z^2/2s} \,
dz \\
&&{}-
\frac{1}{\sqrt{2 \pi s}} \int_x^{K+x} (z-x) e^{-z^2/2s} \, dz
\\
&\leq&\frac{1}{\sqrt{2 \pi s}} \int_{-x}^x z e^{-z^2/2s} \, dz +
\frac
{2x}{\sqrt{2 \pi s}} \int_{-\infty}^{\infty} e^{-z^2/2s} \, dz = 2x.
\end{eqnarray*}
Therefore, using that $v_s(x,y) = v_s(K-x, K-y) \leq v_s'(K-x, K-y)$,
%
%e30 ###
%
%e30 #&#
\begin{eqnarray}\label{term3}
&&\int_{t/2}^t \int_0^K q_s(x,z) \biggl( \int_0^K e^{\mu y} \sin\biggl(
\frac{\pi y}{K} \biggr) q_{t-s}(z,y) \, dy \biggr)^2 \, dz \, ds
\nonumber\\
&& \qquad\leq\int_{t/2}^t \int_0^K \frac{C}{K} e^{\mu(x-z)} \sin
\biggl( \frac
{\pi x}{K} \biggr) \sin\biggl( \frac{\pi z}{K} \biggr)\nonumber\\
&&\hphantom{\leq\int_{t/2}^t \int_0^K} \qquad {}\times\biggl ( \int_0^K
e^{\mu y} e^{\mu(z-y)} \sin\biggl( \frac{\pi y}{K} \biggr) v_{t-s}(z,y)
\, dy \biggr)^2 \, dz \, ds \nonumber\\
&& \qquad\leq\frac{C e^{\mu x}}{K} \int_{t/2}^t \int_0^K e^{\mu z}
\sin\biggl(
\frac{\pi z}{K} \biggr) \biggl( \int_0^K \sin\biggl( \frac{\pi y}{K}
\biggr) v_{t-s}(z,y) \, dy \biggr)^2 \, dz \, ds \nonumber
\\[-8pt]
\\[-8pt]
&& \qquad\leq\frac{C e^{\mu x}}{K} \int_{t/2}^t \int_0^K e^{\mu z}
\sin\biggl(
\frac{\pi z}{K} \biggr) \biggl( \int_0^K \biggl( \frac{K-y}{K} \biggr)
v_{t-s}(z,y) \, dy \biggr)^2 \, dz \, ds \nonumber\\
&& \qquad\leq\frac{C e^{\mu x}}{K^3} \int_{t/2}^t \int_0^K e^{\mu
z} \sin
\biggl( \frac{\pi z}{K} \biggr)\biggl ( \int_0^K y v_{t-s}'(K-z, y) \, dy
\biggr)^2 \, dz \, ds \nonumber\\
&& \qquad\leq\frac{C e^{\mu x}}{K^3} \int_{t/2}^t \int_0^K e^{\mu
z} \sin
\biggl( \frac{\pi z}{K} \biggr) (K-z)^2 \, dz \, ds \nonumber\\
&& \qquad\leq\frac{C e^{\mu x} t}{K^4} \int_0^K e^{\mu z} (K-z)^3
\, dz \leq
\frac{C e^{\mu x} e^{\mu K}t}{K^4}.\nonumber
\end{eqnarray}
The result follows from (\ref{term1}), (\ref{term2}) and (\ref{term3}).
\end{pf}

%s3 ###
%s3 #&#
\section{Particles hitting the right-boundary}\label{sec3}

Recall that we are considering $(X_N(t), t \geq0)$, which is a
branching Brownian motion with drift $-\mu$ and killing at the origin.
Recall also that Proposition~\ref{ZNCSBP} involves the processes
$(Z_N(t), t \geq0)$, where $Z_N(t)$ is a weighted sum of the positions
of the particles at time $t$. Throughout this entire section, as well
as Sections~\ref{sec5},~\ref{sec6} and~\ref{sec7}, we assume that the hypotheses of Proposition
\ref{ZNCSBP} hold.

%s3.1 ###
%s3.1 #&#
\subsection{The particles that never reach $L_A$}\label{sec31}

To prove Proposition~\ref{ZNCSBP}, we will need to consider these
processes at two times $u$ and $u + s$, where $0 \leq u < u + s$. Fix a
small number $\theta> 0$ such that $\theta^{-1} \in\N$. For $0 \leq k
\leq\theta^{-1}$, define the time $t_k = (u + \theta k s) (\log N)^3$.
We will be interested in the value of the process $Z_N$ at the times
$t_k$. The assumption that $\theta^{-1} \in\N$ is useful for defining
the sequence $\{t_k\}_{0 \le k \le\theta^{-1}}$. However, many of our
results pertain to the state of the process at time $t_k$, conditional
on the state of the process up to time $t_{k-1}$. For these results,
the assumption $\theta^{-1} \in\N$ is not necessary.

Since $Y_N(0)/N(\log N)^3$ converges in probability to zero, there
exists a nonrandom function $h\dvtx  \N\rightarrow(0, \infty)$ such that
$h(N) \rightarrow0$ and $(\log N) h(N) \rightarrow\infty$ as $N
\rightarrow\infty$, and $Y_N(0)/(N (\log N)^3 h(N))$ converges in
probability to zero. [This is a simple consequence of the following
fact: if $X_N \to0 $ in probability, then there exists a nonrandom
sequence $h_N$ such that $h_N \to0$ as $N \to\infty$ and \mbox{$P(X_N >
h_N) \to0$}.] Let $\eps> 0$. For $0 \leq k \leq\theta^{-1}$, let
$G_{N,k}$ be the event that for $j = 0, 1,\ldots, k$, the following two
events occur:
\begin{itemize}
\item We have $Z_N(t_j) \leq\eps^{-1/2} N (\log N)^2$.

\item We have $Y_N(t_j) \leq N (\log N)^{3} h(N)$.
\end{itemize}
Finally, let $G_N(\eps) = G_{N, \theta^{-1}}$. Let $({\cal F}_t, t
\geq
0)$ be the natural filtration of $(X_N(t),  t \geq0)$. This filtration,
of course, depends on $N$, but we suppress this dependence in the
notation. We will need to consider the conditional distribution of
$Z_N(t_k)$ given ${\cal F}_{t_{k-1}}$. Note that the event $G_{N, k-1}$
is in ${\cal F}_{t_{k-1}}$.

In this section, we will consider the particles that would still be
alive if, between times $t_{k-1}$ and $t_k$, we killed particles that
hit $L_A$, where $L_A$ was defined in (\ref{LAdef}). Recall that both
$L_A$ and the drift $\mu$ depend on $N$. We will always assume that $N$
is large enough that $L_A > 0$ and
%
%e31 ###
%
%e31 #&#
\begin{equation}\label{hNbound}
h(N) \leq\frac{e^{\mu L_A}}{N (\log N)^3} = e^{-A} e^{(\mu/\sqrt{2} -
1)(\log N + 3 \log\log N - A)},
\end{equation}
which is possible because, by (\ref{muasymp}), the right-hand side
tends to $e^{-A}$ as $N \rightarrow\infty$. Because $Y_N(t_k) \leq N
(\log N)^3 h(N)$ on $G_{N, k}$, this ensures that on $G_{N,k}$, all
particles at time $t_k$ are to the left of $L_A$, a fact which will be
invoked repeatedly in what follows.

Note that we have defined three constants: $\eps$, $A$ and $\theta$. We
think of $\eps$ as being small. Typically $A$ will be a large positive
constant, but we will also at times consider negative values of $A$.
Finally, $\theta$ will always be a small positive constant. In
particular, we will assume
%
%e35 ###
%e34 ###
%e33 ###
%e32 ###
%
%e32 #&#
%e33 #&#
%e34 #&#
%e35 #&#
\begin{eqnarray}\label{con0}
\theta&\leq&1, \\\label{con1}
|A| \theta&\leq&1, \\\label{con2}
4 \pi^2 A \theta s \eps^{-1/2} &\leq& e^{-A/4}, \\\label{con3}
\theta e^A \eps^{-1/2} &\leq&1.
\end{eqnarray}
These assumptions will be in force through the rest of this section,
except in Proposition~\ref{GNProp} below, where it will be convenient
to allow $\theta$ to be any number with $\theta^{-1} \in\N$. A
stronger set of restrictions on $\theta$ will then be introduced at the
beginning of Section~\ref{sec5}.

For $t \in[t_{k-1}, t_k]$, we say $i \in S(t)$ if for all $v \in
[t_{k-1}, t]$, the particle at time $v$ that is the ancestor of
$X_{i,N}(t)$ is in $(0, L_A)$. Consequently, for $t_{k-1} \leq t \leq
t_k$, the positions of the particles in $S(t)$ follow a branching
Brownian motion with drift $-\mu$, with particles killed when they
reach $0$ or $L_A$. Define
%
%e36 ###
%
%e36 #&#
\begin{equation}
\label{ZN1def}
Z_{N,1}(t_k) = \sum_{i=1}^{M_N(t_k)} e^{\mu X_{i,N}(t_k)} \sin\biggl(
\frac{\pi X_{i,N}(t_k)}{L} \biggr) \mathbf{1}_{\{i \in S(t_k)\}},
\end{equation}
and for $t \in[t_{k-1}, t_k]$, define
%
%e37 ###
%
%e37 #&#
\begin{equation}\label{ZNprime}
Z_{N,1}'(t) = \sum_{i=1}^{M_N(t)} e^{\mu X_{i,N}(t)} \sin\biggl( \frac
{\pi X_{i,N}(t)}{L_A} \biggr) \mathbf{1}_{\{i \in S(t)\}}.
\end{equation}
Although our interest is in $Z_{N,1}(t_k)$, we will need to approximate
this random variable by $Z_{N,1}'(t_k)$, which is defined in the same
way except with $L_A$ in place of~$L$.
The next result shows that the difference between these quantities is small.

%
%le10 #&#
\begin{Lemma}\label{Zprimelem}
On $G_{N, k-1}$, both $|Z_{N,1}'(t_{k-1}) - Z_N(t_{k-1})|$ and $E [
| Z_{N,1}'(t_k) - Z_{N,1}(t_k)| | {\cal F}_{t_{k-1}} ]$ are
$o(N (\log N)^2)$.
\end{Lemma}

\begin{pf}
If $a > 0$, then
\[
\biggl| \frac{d}{dx} \sin\biggl( \frac{a}{x} \biggr) \biggr| = \biggl| \frac
{a}{x^2} \cos\biggl( \frac{a}{x} \biggr) \biggr| \leq\frac{a}{x^2}.
\]
Therefore, if $0 \leq x \leq L_A$, then
\[
\biggl| \sin\biggl( \frac{\pi x}{L} \biggr) - \sin\biggl( \frac{\pi
x}{L_A} \biggr) \biggr| \leq\frac{|L - L_A| \pi x}{\min\{L_A, L\}^2}
\leq\frac{\pi|A| L_A}{\sqrt{2} \min\{L_A, L\}^2}.
\]
On $G_{N, k-1}$, all particles at time $t_{k-1}$ are to the left of
both $L_A$ and $L$ for sufficiently large $N$. The indicators are
therefore not needed in (\ref{ZNdef}) and (\ref{ZNprime}) when $t =
t_{k-1}$, and we get
%
%e38 ###
%
%e38 #&#
\begin{eqnarray}\label{ZNYbound}
|Z_{N,1}'(t_{k-1}) - Z_N(t_{k-1})| &\leq&\frac{\pi|A| L_A}{\sqrt{2}
\min
\{L_A, L\}^2} \sum_{i=1}^{M_N(t_{k-1})} e^{\mu X_{i,N}(t_{k-1})}\nonumber
\\[-8pt]
\\[-8pt] &=&
\frac
{\pi|A| L_A Y_N(t_{k-1})}{\sqrt{2} \min\{L_A, L\}^2},
\nonumber
\end{eqnarray}
which is $o(N (\log N)^2)$ on $G_{N, k-1}$. Applying the same reasoning
at time $t_k$ to the particles in $S(t_k)$, we get
\[
E [ | Z_{N,1}'(t_k) - Z_{N,1}(t_k)| | {\cal F}_{t_{k-1}} ]
\leq\frac{\pi|A| L_A E[Y_N(t_k)|{\cal F}_{t_{k-1}}]}{\sqrt{2} \min
\{
L_A, L\}^2}.
\]
Note that
%
%e39 ###
%
%e39 #&#
\begin{eqnarray}\label{LAexp}
1 - \frac{\mu^2}{2} - \frac{\pi^2}{2 L_A^2} &=& \frac{\pi^2}{(\log
N + 3
\log\log N)^2} - \frac{\pi^2}{(\log N + 3 \log\log N - A)^2}\nonumber
\\[-8pt]
\\[-8pt] &=& -
\frac
{2\pi^2 A}{(\log N)^3} \bigl(1 + o(1) \bigr) .
\nonumber
\end{eqnarray}
Since $t_k - t_{k-1} = (\log N)^3 \theta s$ and (\ref{con1}) holds,
equations (\ref{Yexp}) and (\ref{LAexp}) give $E[Y_N(t_k)|{\cal
F}_{t_{k-1}}] \leq C Z_N(t_{k-1})(1 + o(1))$. It follows that
%
%e40 ###
%
%e40 #&#
\begin{equation}\label{new27}
E [ | Z_{N,1}'(t_k) - Z_{N,1}(t_k)| | {\cal F}_{t_{k-1}} ]
\leq\frac{C |A| L_A Z_N(t_{k-1})(1+o(1))}{\min\{L_A, L\}^2},
\end{equation}
which is $o(N (\log N)^2)$ on $G_{N, k-1}$.
\end{pf}

We now estimate the conditional mean and variance of $Z_{N,1}(t_k)$
given~${\cal F}_{t_{k-1}}$.

%
%le11 #&#
\begin{Lemma}\label{ZN1exp}
On $G_{N, k-1}$, we have
\[
E[Z_{N,1}(t_k)|{\cal F}_{t_{k-1}}] = Z_N(t_{k-1}) \bigl(1 - 2 \pi^2 A
\theta
s + O(A^2 \theta^2)\bigr) + o(N (\log N)^2).
\]
The same bound holds with $E[Z'_{N,1}(t_k)|{\cal F}_{t_{k-1}}]$ on the
left-hand side.
\end{Lemma}

\begin{pf}
By (\ref{Zexp}) and the Markov property of branching Brownian motion
with particles killed at $0$ and $L_A$, we have for sufficiently large
$N$ on~$G_{N, k-1}$,
%
%e41 ###
%
%e41 #&#
\begin{equation}\label{ZN1prime}
E[Z_{N,1}'(t_k)|{\cal F}_{t_{k-1}}] = e^{(1 - \mu^2/2 - \pi
^2/2L_A^2)(t_k - t_{k-1})} Z'_{N,1}(t_{k-1}),
\end{equation}
using the fact that for sufficiently large $N$, on $G_{N, k-1}$ all
particles at time $t_{k-1}$ are to the left of $L_A$. Since $t_k -
t_{k-1} = (\log N)^3 \theta s$, it follows\vadjust{\goodbreak} from (\ref{LAexp}) that
%
%e42 ###
%
%e42 #&#
\begin{eqnarray}\label{expest}
e^{(1 - \mu^2/2 - \pi^2/2L_A^2)(t_k - t_{k-1})} &=& e^{-2 \pi^2 A
\theta
s(1 + o(1))} \nonumber
\\[-8pt]
\\[-8pt]&=& 1 - 2 \pi^2 A \theta s + O(A^2 \theta^2) + o(1),
\nonumber
\end{eqnarray}
where assumption (\ref{con1}) ensures that the error term is $O(A^2
\theta^2)$. The result now follows from equations (\ref{ZN1prime}) and
(\ref{expest}) together with the two bounds in Lemma~\ref{Zprimelem}.
\end{pf}

%
%le12 #&#
\begin{Lemma}\label{VarZ1}
Assume $A \geq0$. On $G_{N, k-1}$, we have
\[
\operatorname{Var}(Z_{N,1}'(t_k)|{\cal F}_{t_{k-1}}) \leq C \theta N
(\log
N)^2 e^{-A} \bigl(Z_N(t_{k-1}) + o(N (\log N)^2)\bigr).
\]
\end{Lemma}

\begin{pf}
For $t \in[t_{k-1}, t_k]$, define $Z_{N,1}'(t)$ as in (\ref{ZNprime}),
and define
\[
Y_N'(t) = \sum_{i=1}^{M_N(t)} e^{\mu X_{i,N}(t)} \mathbf{1}_{\{i \in
S(t)\}}.
\]
Define $t_{k-1} = s_0 < s_1 <\cdots< s_M = t_k$ so that for some
positive constants $C_1$ and $C_2$, we have $C_1 (\log N)^2 \leq s_n -
s_{n-1} \leq C_2 (\log N)^2$ for all $n$. Recall that for any random
variable $X$ and any $\sigma$-fields ${\cal F}$ and ${\cal G}$ with
${\cal F} \subset{\cal G}$, we have
\[
\operatorname{Var}(X|{\cal F}) = E[\operatorname{Var}(X|{\cal
G})|{\cal F}] + \operatorname
{Var}(E[X|{\cal G}]|{\cal F}).
\]
Therefore, for $1 \leq n \leq M$, we have
\[
\operatorname{Var}(Z_{N,1}'(s_n)|{\cal F}_{s_0}) = E[\operatorname
{Var}(Z_{N,1}'(s_n)|{\cal F}_{s_{n-1}})|{\cal F}_{s_0}] + \operatorname
{Var}(E[Z_{N,1}'(s_n)|{\cal F}_{s_{n-1}}]|{\cal F}_{s_0}).
\]
Equation (\ref{Zexp}) implies that $E[Z_{N,1}'(s_n)|{\cal F}_{s_{n-1}}]
= e^{(1 - \mu^2/2 - \pi^2/2L_A^2)(s_n - s_{n-1})}\times  Z_{N,1}'(s_{n-1})$.
Because $A \geq0$ and thus $1 - \mu^2/2 - \pi^2/2L_A^2 \leq0$, it
follows that
\[
\operatorname{Var}(E[Z_{N,1}'(s_n)|{\cal F}_{s_{n-1}}]|{\cal F}_{s_0})
\leq
\operatorname{Var}(Z_{N,1}'(s_{n-1})|{\cal F}_{s_0}).
\]
Therefore,
\[
\operatorname{Var}(Z_{N,1}'(s_n)|{\cal F}_{s_0}) \leq E[\operatorname
{Var}(Z_{N,1}'(s_n)|{\cal F}_{s_{n-1}})|{\cal F}_{s_0}] + \operatorname
{Var}(Z_{N,1}'(s_{n-1})|{\cal F}_{s_0}).
\]
Now $\operatorname{Var}(Z_{N,1}'(s_0)|{\cal F}_{s_0}) = 0$, so by induction,
%
%e43 ###
%
%e43 #&#
\begin{equation}\label{vareq1}
\operatorname{Var}(Z_{N,1}'(s_M)|{\cal F}_{s_0}) \leq\sum_{n=1}^M
E[\operatorname
{Var}(Z_{N,1}'(s_n)|{\cal F}_{s_{n-1}})|{\cal F}_{s_0}].
\end{equation}

Because the particles at time $s_{n-1}$ evolve independently between
times $s_{n-1}$ and $s_n$, the conditional variance $\operatorname
{Var}(Z_{N,1}'(s_n)|{\cal F}_{s_{n-1}})$ is the sum of the conditional
variances of the contributions to $Z_{N,1}'(s_n)$ from the individual
particles at time $s_{n-1}$. We will use the inequality $\operatorname
{Var}(X|{\cal F}) \leq E[X^2|{\cal F}]$ and apply Lemma~\ref{varZlem}
with $K = L_A$ and $t = s_n - s_{n-1}$. The hypotheses are satisfied
because $1 - \mu^2/2 - \pi^2/2L_A^2 \leq0$, and both $s_n - s_{n-1}$
and $1/(1 - \mu^2/2)$ are of the order $(\log N)^2$. Therefore,
\[
\operatorname{Var}(Z_{N,1}'(s_n)|{\cal F}_{s_{n-1}}) \leq C e^{\mu L_A}
Y_N'(s_{n-1})\biggl ( \frac{1}{L_A^2} + \frac{s_n - s_{n-1}}{L_A^4} \biggr).
\]
Now $e^{\mu L_A} \leq N (\log N)^3 e^{-A}$, so
%
%e44 ###
%
%e44 #&#
\begin{equation}\label{vareq2}
\operatorname{Var}(Z_{N,1}'(s_n)|{\cal F}_{s_{n-1}}) \leq C N (\log
N)^3 e^{-A}
\biggl( \frac{1}{L_A^2} + \frac{(\log N)^2}{L_A^4} \biggr).
\end{equation}
From (\ref{Yexp}), we get
%
%e45 ###
%
%e45 #&#
\begin{eqnarray}\label{vareq3}
\max_{2 \leq n \leq M} E[Y_N'(s_{n-1})|{\cal F}_{s_0}] &\leq& C
Z_{N,1}'(s_0)\bigl(1 + o(1)\bigr)\nonumber
\\[-8pt]
\\[-8pt] &=& C Z_{N,1}'(t_{k-1})\bigl(1 + o(1)\bigr).
\nonumber
\end{eqnarray}
Finally, note that $M \leq C \theta(\log N)$. Combining this with
(\ref
{vareq1}), (\ref{vareq2}) and (\ref{vareq3}) gives that on $G_{N, k-1}$,
%
%e46 ###
%
%e46 #&#
\begin{eqnarray} \label{varZNp}
\operatorname{Var}(Z_{N,1}'(t_k)|{\cal F}_{t_{k-1}}) &=& \operatorname
{Var}(Z_{N,1}'(s_M)|{\cal F}_{s_0}) \nonumber\\
&\leq& C N (\log N)^3 e^{-A} \biggl( \frac{1}{L_A^2} + \frac{(\log
N)^2}{L_A^4} \biggr) \nonumber
\\[-8pt]
\\[-8pt]
&&{} \times\bigl( Y_N'(s_0) + C \theta(\log N)
Z_{N,1}'(t_{k-1})\bigl(1 + o(1)\bigr) \bigr) \nonumber\\
&\leq& C \theta N (\log N)^2 e^{-A}\biggl (\frac{Y_N'(t_{k-1})}{\theta
\log N} + Z_{N,1}'(t_{k-1}) \biggr)\bigl (1 + o(1)\bigr).
\nonumber
\end{eqnarray}
The result now follows from Lemma~\ref{Zprimelem} and the fact that
$Y_N'(t_{k-1}) \leq\break Y_N(t_{k-1}) \leq N (\log N)^3 h(N)$ on $G_{N, k-1}$.
\end{pf}

%
%co13 #&#
\begin{Cor}\label{driftcor}
Assume $A \geq0$. On $G_{N, k-1}$, we have
\[
P \bigl( |Z_{N,1}(t_k) - Z_N(t_{k-1})| > 4 e^{-A/4} N (\log N)^2 |
{\cal F}_{t_{k-1}} \bigr) \leq C \theta e^{-A/2} \eps^{-1/2}\bigl(1 + o(1)\bigr).
\]
\end{Cor}

\begin{pf}
By the conditional form of Chebyshev's inequality and Lemma~\ref
{VarZ1}, on $G_{N, k-1}$ we have
%
%e47 ###
%
%e47 #&#
\begin{eqnarray} \label{dcor1}
&&P\bigl(|Z_{N, 1}'(t_k) - E[Z_{N,1}'(t_k)|{\cal F}_{t_{k-1}}]| > e^{-A/4} N
(\log N)^2 |{\cal F}_{t_{k-1}}\bigr)\nonumber\\
&& \qquad \leq\frac{\operatorname{Var}(Z_{N,1}'(t_k)|{\cal
F}_{t_{k-1}})}{e^{-A/2} N^2
(\log N)^4} \\
&& \qquad \leq C \theta e^{-A/2} \eps^{-1/2} \bigl(1 + o(1)\bigr)
\nonumber
\end{eqnarray}
because $Z_N(t_{k-1}) \leq\eps^{-1/2} N (\log N)^2$ on $G_{N, k-1}$.
Using (\ref{LAexp}), some calculus and the assumption that\vadjust{\goodbreak} $A \geq0$,
we get that for $N$ large enough that $A \leq3 \log\log N$,
\[
\bigl| e^{(1 - \mu^2/2 - \pi^2/2L_A^2)(t_k - t_{k-1})} - 1 \bigr| \leq
\biggl| 1 - \frac{\mu^2}{2} - \frac{\pi^2}{2L_A^2} \biggr| \theta s (\log
N)^3 \leq2 \pi^2 A \theta s.
\]
Therefore, by (\ref{ZN1prime}), if $A \leq3 \log\log N$, then
\[
| E[Z_{N,1}'(t_k)|{\cal F}_{t_{k-1}}] - Z'_{N,1}(t_{k-1}) |
\leq2 \pi^2 A \theta s Z_{N,1}'(t_{k-1}).
\]
Because $Z_{N,1}'(t_{k-1}) = Z_N(t_{k-1}) + o(N (\log N)^2) \leq\eps
^{-1/2} N (\log N)^2 +\break o(N (\log N)^2)$ on $G_{N, k-1}$ by Lemma \ref
{Zprimelem} and $2 \pi^2 A \theta s \eps^{-1/2} \leq e^{-A/4}/2$ by
(\ref{con2}), it follows that for sufficiently large $N$,
%
%e48 ###
%
%e48 #&#
\begin{equation}\label{dcor2}
| E[Z_{N,1}'(t_k)|{\cal F}_{t_{k-1}}] - Z_{N,1}'(t_{k-1}) |
\leq e^{-A/4} N (\log N)^2
\end{equation}
on $G_{N, k-1}$. By Lemma~\ref{Zprimelem}, on $G_{N, k-1}$, we have
%
%e49 ###
%
%e49 #&#
\begin{equation}\label{dcor3}
|Z_{N,1}'(t_{k-1}) - Z_N(t_{k-1})| \leq e^{-A/4} N (\log N)^2
\end{equation}
for sufficiently large $N$ and
%
%e50 ###
%
%e50 #&#
\begin{equation}\label{dcor4}
P\bigl(|Z_{N,1}(t_k) - Z_{N,1}'(t_k)| > e^{-A/4} N (\log N)^2|{\cal
F}_{t_{k-1}}\bigr) \rightarrow0
\end{equation}
uniformly as $N \rightarrow\infty$ on $G_{N,k-1}$. The result follows
immediately from (\ref{dcor1}), (\ref{dcor2}), (\ref{dcor3}) and
(\ref{dcor4}).
\end{pf}

%
%pr14 #&#
\begin{Prop}\label{MNVar}
Suppose $A = 0$. Let
\[
M_N'(t_k) = \sum_{i=1}^{M_N(t_k)} \mathbf{1}_{\{i \in S(t_k)\}}
\]
be the number of particles at time $t_k$ whose ancestor at time $t$ is
in $(0, L)$ for all $t \in[t_{k-1}, t_k]$. On $G_{N, k-1}$, there
exists a constant $C$ such that
\[
\operatorname{Var}(M_N'(t_k)|{\cal F}_{t_{k-1}}) \leq C \theta\eps^{-1/2}
N^2\bigl (1 + o(1)\bigr).
\]
\end{Prop}

\begin{pf}
As in the proof of Lemma~\ref{VarZ1}, the conditional variance can be
bounded by the sum of the variances of the contributions to $M'_N(t_k)$
from the individual particles at time $t_{k-1}$. The variance of the
contribution from a particle at $x$ can be bounded by the expected
square of the number of descendants of this particle at time $t_k$.
This expectation is given by Lemma~\ref{kestenlem} with $f(x) = 1$ for
all $x$ and $t_k - t_{k-1}$ in place of $t$. Therefore,
\begin{eqnarray*}
&&\operatorname{Var}(M_N'(t_k)|{\cal F}_{t_{k-1}})\\ && \qquad = \sum_{i=1}^{M_N(t_{k-1})}
\int_0^L q_{t_k - t_{k-1}}(X_{i,N}(t_{k-1}), y) \, dy \\
&& \qquad  \quad {}+ 2 \sum_{i=1}^{M_N(t_{k-1})} \int_{t_{k-1}}^{t_k} \int_0^L q_{t-
t_{k-1}}(X_{i,N}(t_{k-1}), z) \biggl( \int_0^L q_{t_k - t}(z,y) \, dy
\biggr)^2 \, dz \, dt.
\end{eqnarray*}
The first term is $E[M_N'(t_k)|{\cal F}_{t_{k-1}}]$, which by (\ref
{Mexp}) with $K = L$ is at most $C Z_N(t_{k-1})(1 + o(1))/L^2$ because
the integral on the right-hand side of (\ref{Mexp}) is of the order $1/K$.
This expression is $o(N^2)$ on $G_{N, k-1}$.

The argument to bound the second term is similar to the proof of Lemma~\ref{varZlem}
but requires splitting the outer integral into four
pieces. First consider the piece between $t_{k-1}$ and $t_{k-1} + (\log
N)^2$. If $t \leq(\log N)^2$, then (\ref{qsvs}) holds and
%
%e51 ###
%
%e51 #&#
\begin{equation}\label{GreenLz}
\int_0^{\infty} v_t(x,y) \, ds = G(x,y) \leq\frac{2x(L-y)}{L} \leq2(L
- y)
\end{equation}
by (\ref{Gbound}). Since $1 - \mu^2/2 - \pi^2/2L^2 = 0$, Lemma \ref
{stripdensity} gives that on $G_{N, k-1}$,
%
%e52 ###
%
%e52 #&#
\begin{eqnarray}\label{bit1}
&&\sum_{i=1}^{M_N(t_{k-1})} \!\!\int_{t_{k-1}}^{t_{k-1} + (\log N)^2}\!\!
\int
_0^L q_{t - t_{k-1}}(X_{i,N}(t_{k-1}), z) \biggl( \int_0^L q_{t_k -
t}(z,y) \, dy \biggr)^2 \, dz \, dt \nonumber\hspace*{-25pt}\\
&& \quad\leq C \sum_{i=1}^{M_N(t_{k-1})}\!\! \int_0^{(\log N)^2}\!\! \int_0^L
q_t(X_{i,N}(t_{k-1}), z) \nonumber\hspace*{-25pt}\\
&& \quad \hphantom{\leq C \sum_{i=1}^{M_N(t_{k-1})} \!\!\int_0^{(\log N)^2}\!\!
\int_0^L}{}\times\biggl( \int_0^L \frac{2}{L} e^{\mu z} \sin
\biggl( \frac{\pi z}{L} \biggr) e^{-\mu y} \sin\biggl( \frac{\pi y}{L}
\biggr) \, dy \biggr)^2 \, dz \, dt \nonumber\hspace*{-25pt}\\
&& \quad\!\!\leq\frac{C}{L^4} \sum_{i=1}^{M_N(t_{k-1})}\!\! \int_0^{(\log
N)^2}\!\! \int
_0^L q_t(X_{i,N}(t_{k-1}), z) e^{2 \mu z} \sin\biggl( \frac{\pi z}{L}
\biggr)^2 \, dz \, dt\nonumber\hspace*{-25pt}\\
&& \quad\leq\frac{C}{L^4} \sum_{i=1}^{M_N(t_{k-1})} \!\!\int_0^L e^{2
\mu z} \sin
\biggl( \frac{\pi z}{L} \biggr)^2 \biggl( \int_0^{(\log N)^2}
q_t(X_{i,N}(t_{k-1}), z) \, dt \biggr) \, dz  \hspace*{-25pt}\\
&& \quad\leq\frac{C}{L^4} \sum_{i=1}^{M_N(t_{k-1})} \!\!\int_0^L e^{2
\mu z} \sin
\biggl( \frac{\pi z}{L} \biggr)^2 \nonumber\hspace*{-25pt}\\
&& \quad \hphantom{\leq\frac{C}{L^4} \sum_{i=1}^{M_N(t_{k-1})}\!\!
\int_0^L}{}\times e^{\mu(X_{i,N}(t_{k-1}) - z)}\biggl (
\int_0^{\infty} v_t(X_{i,N}(t_{k-1}), z) \, dt \biggr) \, dz \nonumber\hspace*{-25pt}\\
&& \quad\leq\frac{C}{L^4}  \Biggl( \sum_{i=1}^{M_N(t_{k-1})} e^{\mu
X_{i,N}(t_{k-1})}  \Biggr) \int_0^L e^{\mu z} \sin\biggl( \frac{\pi z}{L}
\biggr)^2 (L - z) \, dz \nonumber\hspace*{-25pt}\\
&& \quad\leq\frac{C}{L^4} \cdot Y_N(t_{k-1}) \cdot\frac{e^{\mu
L}}{L^2} \leq
C N^2 h(N).
\nonumber\hspace*{-25pt}
\end{eqnarray}

We next consider the case $t_{k-1} + (\log N)^2 \leq t \leq t_k - (\log
N)^2$, and from Lemma~\ref{stripdensity}, we get that on $G_{N, k-1}$,
%
%e53 ###
%
%e53 #&#
\begin{eqnarray}\label{bit2}
&& \sum_{i=1}^{M_N(t_{k-1})} \int_{t_{k-1} + (\log N)^2}^{t_k - (\log
N)^2}\!\! \int_0^L q_{t - t_{k-1}}(X_{i,N}(t_{k-1}), z) \biggl( \int_0^L
q_{t_k - t}(z,y) \, dy \biggr)^2 \, dz \, dt \nonumber\\
&& \quad\leq C \sum_{i=1}^{M_N(t_{k-1})} \int_{t_{k-1} + (\log
N)^2}^{t_k -
(\log N)^2} \int_0^L \frac{2}{L} e^{\mu X_{i,N}(t_{k-1})} \sin\biggl(
\frac{\pi X_{i,N}(t_{k-1})}{L} \biggr)  \nonumber\\
&&\hspace*{-25pt}\hphantom{\leq C \sum_{i=1}^{M_N(t_{k-1})} \int_{t_{k-1} + (\log
N)^2}^{t_k -
(\log N)^2} \int_0^L}\qquad  {} \times e^{-\mu z} \sin\biggl( \frac{\pi
z}{L} \biggr)\nonumber\\
&& \hspace*{-25pt}\hphantom{\leq C \sum_{i=1}^{M_N(t_{k-1})} \int_{t_{k-1} + (\log
N)^2}^{t_k -
(\log N)^2} \int_0^L}\qquad  {} \times\biggl( \int_0^L \frac{2}{L} e^{\mu z} \sin\biggl(
\frac{\pi z}{L} \biggr) e^{-\mu y}
 \sin\biggl( \frac{\pi y}{L} \biggr) \,
dy \biggr)^2 \, dz \, dt \\
&& \quad\leq\frac{C (t_k - t_{k-1}) Z_N(t_{k-1})}{L^3} \biggl( \int_0^L
e^{\mu
z} \sin\biggl( \frac{\pi z}{L} \biggr)^3 \, dz \biggr) \biggl[ \int_0^L
e^{-\mu y} \sin\biggl( \frac{\pi y}{L} \biggr) \, dy \biggr]^2 \nonumber
\\
&& \quad\leq\frac{C \theta(\log N)^3 Z_N(t_{k-1})}{L^3} \cdot
\frac{e^{\mu
L}}{L^3} \cdot\frac{1}{L^2} \leq C \theta\eps^{-1/2} N^2.\nonumber
\end{eqnarray}

Consider now the case $t_k - (\log N)^2 \leq t \leq t_k - (\log
N)^{7/4}$. Note that if $t \leq C (\log N)^2$, then $e^{(1 - \mu^2/2)t}
\leq C$, so by (\ref{density}) and Lemma~\ref{trig},
\[
q_t(x,y) \leq\frac{C}{L} e^{\mu x} \sin\biggl( \frac{\pi x}{L} \biggr)
e^{-\mu y} \sin\biggl( \frac{\pi y}{L} \biggr) \sum_{n=1}^{\infty} n^2
e^{-\pi^2 n^2 t/L^2}.
\]
Breaking up the sum into blocks of size $M = \lceil L/\sqrt{t} \rceil$
gives
\begin{eqnarray*}
\sum_{n=1}^{\infty} n^2 e^{-\pi^2 n^2 t/L^2} &\leq&\sum_{\ell=
0}^{\infty} M \bigl(M (\ell+ 1)\bigr)^2 e^{-\pi^2 (M \ell)^2 t/L^2} \leq M^3
\sum
_{\ell= 0}^{\infty} (\ell+ 1)^2 e^{-\pi^2 \ell^2}\\
 &\leq&\frac{C
L^3}{t^{3/2}}.
\end{eqnarray*}
Therefore,
%
%e54 ###
%
%e54 #&#
\begin{eqnarray}\label{bit3}
&&\sum_{i=1}^{M_N(t_{k-1})} \int_{t_k - (\log N)^2}^{t_k - (\log
N)^{7/4}} \int_0^L q_{t - t_{k-1}}(X_{i,N}(t_{k-1}), z) \biggl( \int_0^L
q_{t_k - t}(z,y) \, dy \biggr)^2 \, dz \, dt \nonumber\\
&& \quad\leq C \sum_{i=1}^{M_N(t_{k-1})} \int_{(\log
N)^{7/4}}^{(\log N)^2}
\int_0^L \frac{2}{L} e^{\mu X_{i,N}(t_{k-1})} \sin\biggl( \frac{\pi
X_{i,N}(t_{k-1})}{L} \biggr) e^{-\mu z} \sin\biggl( \frac{\pi z}{L}
\biggr) \nonumber\\
&& \hspace*{93pt}\qquad {} \times\biggl( \int_0^L \frac{L^2}{t^{3/2}} e^{\mu z}
\sin
\biggl( \frac{\pi z}{L} \biggr) e^{-\mu y} \sin\biggl( \frac{\pi y}{L}
\biggr) \, dy \biggr)^2 \, dz \, dt \nonumber
\\[-12pt]
\\[-4pt]
&& \quad\leq C L^3 Z_N(t_{k-1}) \biggl( \int_{(\log N)^{7/4}}^{(\log N)^2}
\frac
{1}{t^3} \, dt \biggr)\biggl ( \int_0^L e^{\mu z} \sin\biggl( \frac{\pi
z}{L} \biggr)^3 \, dz \biggr)\nonumber\\
&& \qquad{}\times\biggl [ \int_0^L e^{-\mu y} \sin\biggl( \frac
{\pi y}{L} \biggr) \, dy \biggr]^2 \nonumber\\
&& \quad\leq C L^3 Z_N(t_{k-1}) \cdot\frac{1}{(\log N)^{7/2}} \cdot
\frac
{e^{\mu L}}{L^3} \cdot\frac{1}{L^2}
\leq\frac{C \eps^{-1/2} N^2}{(\log N)^{1/2}}.\nonumber
\end{eqnarray}

Next, consider the case $t_k - (\log N)^{7/4} \leq t \leq t_k - 1$.
Using (\ref{qsvs}), and the obvious fact that the density $v_t(x,y)$ of
Brownian motion killed at 0 and $L$ is dominated by the transition
probabilities of standard Brownian motion, for $t \leq(\log N)^2$, we have
%
%e55 ###
%
%e55 #&#
\begin{equation}\label{qshf}
 \qquad q_t(x,y) \leq C e^{\mu(x-y)} v_t(x,y) \leq C e^{\mu x} e^{-\mu y}
\cdot
\frac{1}{t^{1/2}} e^{-(x-y)^2/2t} \leq\frac{C e^{\mu x} e^{-\mu y}}{t^{1/2}}.
\end{equation}
We split the integral over $z$ into two pieces and obtain
%
%e56 ###
%
%e56 #&#
\begin{eqnarray}\label{bit4}
&&\sum_{i=1}^{M_N(t_{k-1})} \int_{t_k - (\log N)^{7/4}}^{t_k - 1}
\int
_0^{2L/3} q_{t - t_{k-1}}(X_{i,N}(t_{k-1}), z) \biggl( \int_0^L q_{t_k -
t}(z,y) \, dy \biggr)^2 \, dz \, dt \nonumber\hspace*{-25pt}\\
&& \quad\leq\frac{C Z_N(t_{k-1})}{L} \int_1^{(\log N)^{7/4}} \int_0^{2L/3}
e^{-\mu z} \sin\biggl( \frac{\pi z}{L} \biggr)\biggl ( \int_0^L q_t(z,y)
\, dy \biggr)^2 \, dz \, dt \nonumber\hspace*{-25pt}\\
&& \quad\leq\frac{C Z_N(t_{k-1})}{L} \int_1^{(\log N)^{7/4}} \int_0^{2L/3}
e^{-\mu z} \sin\biggl( \frac{\pi z}{L} \biggr)\biggl ( \int_0^L \frac
{e^{\mu z} e^{-\mu y}}{t^{1/2}} \, dy \biggr)^2 \, dz \, dt\hspace*{-25pt} \\
&& \quad\leq\frac{C Z_N(t_{k-1})}{L} \biggl( \int_1^{(\log N)^{7/4}}
\frac
{1}{t}\, dt \biggr)\biggl ( \int_0^{2L/3} e^{\mu z} \sin\biggl( \frac{\pi
z}{L} \biggr) \, dz \biggr)\biggl ( \int_0^L e^{-\mu y} \, dy \biggr)
\nonumber\hspace*{-25pt}\\
&& \quad\leq\frac{C Z_N(t_{k-1})}{L} \cdot\log\log N \cdot e^{2
\mu L/3} = o(N^2)\nonumber\hspace*{-25pt}
\end{eqnarray}
and
%
%e57 ###
%
%e57 #&#
\begin{eqnarray}\label{bit5}
&&\sum_{i=1}^{M_N(t_{k-1})} \int_{t_k - (\log N)^{7/4}}^{t_k - 1}
\int
_{2L/3}^L q_{t - t_{k-1}}(X_{i,N}(t_{k-1}), z) \biggl( \int_0^L q_{t_k -
t}(z,y) \, dy \biggr)^2 \, dz \, dt \nonumber\hspace*{-25pt}\\
&& \quad\leq\frac{C Z_N(t_{k-1})}{L} \int_1^{(\log N)^{7/4}} \int_{2L/3}^L
e^{-\mu z} \sin\biggl( \frac{\pi z}{L} \biggr)\biggl ( \int_0^L q_t(z,y)
\, dy \biggr)^2 \, dz \, dt \nonumber\hspace*{-25pt}\\
&& \quad\leq\frac{C Z_N(t_{k-1})}{L} \nonumber\hspace*{-25pt}
\\[-8pt]
\\[-8pt]
&& \quad  \quad {}\times\int_1^{(\log N)^{7/4}} \int_{2L/3}^L
e^{\mu z}\biggl ( \int_0^{L/3} e^{-\mu y} e^{-(z-y)^2/2t} \, dy + \int
_{L/3}^L e^{-\mu y} \, dy \biggr)^2 \, dz \, dt \nonumber\hspace*{-25pt}\\
&& \quad\leq\frac{C Z_N(t_{k-1})}{L} (\log N)^{7/4} e^{\mu L}\biggl (
\int
_0^{L/3} e^{-\mu y} e^{-(\log N)^{1/4}/36} \, dy + e^{-\mu L/3}
\biggr)^2 \nonumber\hspace*{-25pt}\\
&& \quad\leq C N (\log N)^{15/4} Z_N(t_{k-1}) \bigl( e^{-(\log
N)^{1/4}/36} +
e^{-\mu L/3} \bigr)^2 = o(N^2).\nonumber\hspace*{-25pt}
\end{eqnarray}

Finally, if $0 \leq t \leq1$, then $\int_0^L q_t(z,y) \, dy \leq e^t
\leq C$, so
%
%e58 ###
%
%e58 #&#
\begin{eqnarray}\label{bit6}
&&\sum_{i=1}^{M_N(t_{k-1})} \int_{t_k - 1}^{t_k} \int_0^L q_{t -
t_{k-1}}(X_{i,N}(t_{k-1}), z) \biggl( \int_0^L q_{t_k - t}(z,y) \, dy
\biggr)^2 \, dz \, dt \nonumber\\
&& \qquad\leq\frac{C Z_N(t_{k-1})}{L} \int_0^1 \int_0^L e^{-\mu z}
\sin\biggl(
\frac{\pi z}{L} \biggr)\biggl ( \int_0^L q_t(z,y) \, dy \biggr)^2 \, dz \,
dt \\
&& \qquad\leq\frac{C Z_N(t_{k-1})}{L} \int_0^L e^{-\mu z} \sin\biggl(
\frac
{\pi z}{L} \biggr) \, dz \leq\frac{C Z_N(t_{k-1})}{L^2} = o(N^2).\nonumber
\end{eqnarray}
The result now follows from (\ref{bit1})--(\ref{bit3}) and
(\ref{bit4})--(\ref{bit6}).
\end{pf}

%s3.2 ###
%s3.2 #&#
\subsection{The number of particles that hit $L_A$}\label{sec32}

For $k \in\N$, let $R_k$ denote the number of times $t$ between
$t_{k-1}$ and $t_k$ that a particle reaches $L_A$ at time $t$ and, for
all $u \in[t_{k-1}, t)$,
the ancestor of this particle at time $u$ was in $(0, L_A)$.
Equivalently, $R_k$ is the number of particles that are killed by
hitting $L_A$ between times $t_{k-1}$ and $t_k$. Note that particles
can reach $L_A$ before time $t_{k-1}$ and still contribute to $R_k$.
Below we calculate the conditional mean and second moment of $R_k$
given ${\cal F}_{t_{k-1}}$.

%
%le15 #&#
\begin{Lemma}\label{tdeltalem}
Suppose there is a single particle at $x$ at time zero, where $x \in
(0, L_A)$. Suppose particles undergo branching Brownian motion with
drift $-\mu$ and are killed when they reach $0$ or $L_A$. Let $R$ be
the number of particles that hit $L_A$ between times $t$ and $t +
\kappa
$, where $0 < \kappa< 1$. Then
%
%e59 ###
%
%e59 #&#
\begin{eqnarray}\label{finupper}
E[R] &=& 2 \pi e^A \kappa e^{(1 - \mu^2/2 - \pi^2/2L_A^2) t} \nonumber
\\[-8pt]
\\[-8pt]&&{}\cdot
e^{\mu
x}\sin\biggl(\frac{\pi x}{L_A} \biggr) \frac{(1 + D)(1 + o(1))(1 +
o(\kappa))}{N (\log N)^5},
\nonumber
\end{eqnarray}
where $|D|$ is bounded by the right-hand side of (\ref{Eterm}) with
$L_A$ in place of $K$, and $o(\kappa)$ is a term whose absolute value
is bounded by $g(\kappa)$ for some bounded function $g\dvtx  (0, 1)
\rightarrow(0, \infty)$ with $\lim_{\kappa\rightarrow0} g(\kappa)
= 0$.
\end{Lemma}

\begin{pf}
Let $(B_t, t \geq0)$ be standard Brownian motion started at the
origin. Suppose that (for the branching Brownian motion), there is a
particle at $y$ at time $t$, and let $0 < \kappa< 1$. The expected
number of descendants of the particle at time $t + \kappa$
is $e^{\kappa}$, and the drift of $- \mu$ can only reduce the
probability that a Brownian particle reaches $L_A$.
Therefore, an upper bound for the expected number of descendants that
reach $L_A$ at by time $t + \kappa$ is
%
%e60 ###
%
%e60 #&#
\begin{eqnarray}\label{updelta}
e^{\kappa} P \Bigl( \max_{0 \leq t \leq\kappa} B_t \geq L_A - y \Bigr) &=&
2 e^{\kappa} P(B_{\kappa} \geq L_A - y)\nonumber
\\[-8pt]
\\[-8pt] &=& e^{\kappa} \sqrt{\frac
{2}{\kappa\pi}} \int_{L_A - y}^{\infty} e^{-z^2/2 \kappa} \, dz,
\nonumber
\end{eqnarray}
where the first equality follows from the reflection principle. To get
a lower bound, we may ignore the branching, and bound the probability
that Brownian motion with drift $-\mu$ reaches $L_A$ by time $\kappa$
without hitting the origin by the probability that ordinary Brownian
motion reaches $L_A + \mu\kappa$ by time $\kappa$ without hitting
$\mu
\kappa$. For $y \geq\mu\kappa$, this leads to a lower bound of
%
%e61 ###
%
%e61 #&#
\begin{eqnarray}\label{lowdelta}
&&P \Bigl( \max_{0 \leq t \leq\kappa} B_t \geq L_A - y + \mu\kappa
\Bigr) - P \Bigl( \min_{0 \leq t \leq\kappa} B_t \leq-y +\mu\kappa\Bigr)
\nonumber
\\[-8pt]
\\[-8pt]
&& \qquad= \sqrt{\frac{2}{\kappa\pi}}\biggl ( \int_{L_A - y + \mu
\kappa
}^{\infty} e^{-z^2/2 \kappa} \, dz - \int_{y - \mu\kappa}^{\infty}
e^{-z^2/2\kappa} \, dz \biggr).
\nonumber
\end{eqnarray}

From Lemma~\ref{stripdensity}, the expected number of particles in the
set $B$ at time $t$ is $(\int_B p_t(x,y) \, dy)(1 + D'_t(x,B))$. Now
integrating over $y$ and applying (\ref{updelta}), we get
%
%e62 ###
%
%e62 #&#
\begin{eqnarray}\label{uphit}
E[R] &\leq&(1 + D) \int_0^{L_A} p_t(x, y) \cdot e^{\kappa} \sqrt
{\frac
{2}{\kappa\pi}} \biggl( \int_{L_A - y}^{\infty} e^{-z^2/2 \kappa} \, dz
\biggr) \, dy \nonumber\\
&=& (1 + D) e^{\kappa} \sqrt{ \frac{2}{\kappa\pi}} \cdot\frac{2}{L_A}
e^{(1 - \mu^2/2 - \pi^2/2L_A^2)t} e^{\mu x} \sin\biggl( \frac{\pi
x}{L_A} \biggr) \\
&&{} \times\int_0^{L_A} \int_{L_A - y}^{\infty} e^{-\mu y}
\sin\biggl( \frac{\pi y}{L_A} \biggr) e^{-z^2/2 \kappa}\,dz \, dy,\nonumber
\end{eqnarray}
where $|D|$ is bounded by the right-hand side of (\ref{Eterm}) with
$L_A$ in place of $K$.
Interchanging the roles of $y$ and $L_A - y$, then using Fubini's
theorem followed by the bound $\sin y \leq y$ for $y \geq0$ gives
%
%e63 ###
%
%e63 #&#
\begin{eqnarray}\label{uphit2}
&&\int_0^{L_A} \int_{L_A - y}^{\infty} e^{-\mu y} \sin\biggl( \frac{\pi
y}{L_A} \biggr) e^{-z^2/2 \kappa} \, dz \, dy\nonumber\\
 && \qquad = e^{-\mu L_A} \int
_0^{L_A} \int_y^{\infty} e^{\mu y} \sin\biggl( \frac{\pi y}{L_A} \biggr)
e^{-z^2/2 \kappa} \, dz \, dy \nonumber\\
 && \qquad = e^{-\mu L_A} \int_0^{\infty} \int_0^{\min\{z, L_A\}} e^{\mu y}
\sin
\biggl( \frac{\pi y}{L_A} \biggr) e^{-z^2/2 \kappa} \, dy \, dz
\nonumber
\\[-8pt]
\\[-8pt]
&& \qquad \leq e^{-\mu L_A} \int_0^{\infty} \biggl( \int_0^z \frac{\pi y}{L_A} \,
dy \biggr) e^{\mu z} e^{-z^2/2 \kappa} \, dz \nonumber\\
 && \qquad = \frac{\pi e^{-\mu L_A}}{2 L_A} \int_0^{\infty} z^2 e^{\mu z}
e^{-z^2/2 \kappa} \, dz \nonumber\\
 && \qquad = \frac{\pi e^{-\mu L_A}}{2 L_A} \biggl(\kappa^{3/2} \sqrt{\frac{\pi
}{2}} + \int_0^{\infty} z^2 (e^{\mu z} - 1) e^{-z^2/2 \kappa} \, dz
\biggr).\nonumber
\end{eqnarray}
The substitution $y = z/\sqrt{\kappa}$ gives
%
%e64 ###
%
%e64 #&#
\begin{equation}\label{uphit3}
\int_0^{\infty} z^2 (e^{\mu z} - 1) e^{-z^2/2 \kappa} \, dz = \kappa
^{3/2} \int_{0}^{\infty} (e^{\mu y \sqrt{\kappa}} - 1) y^2 e^{-y^2/2}
\, dy,
\end{equation}
and the last integral goes to zero as $\kappa\rightarrow0$ by the
dominated convergence theorem. Therefore, combining (\ref{uphit}),
(\ref
{uphit2}) and (\ref{uphit3}), we get
\begin{eqnarray*}
E[R] &\leq&(1 + D) e^{\kappa} \sqrt{ \frac{2}{\kappa\pi}} \cdot
\frac
{2}{L_A} e^{(1 - \mu^2/2 - \pi^2/2L_A^2) t} e^{\mu x} \sin\biggl( \frac
{\pi x}{L_A} \biggr) \\
&&{}\times\frac{\pi^{3/2} e^{-\mu L_A} \kappa
^{3/2}}{2^{3/2} L_A} \bigl(1 + o(\kappa)\bigr).
\end{eqnarray*}
Since $L_A = (2^{-1/2} \log N) (1 + o(1))$ and $e^{-\mu L_A} = e^A(1 +
o(1))/(N (\log N)^3)$, it follows that $E[R]$ is bounded above by the
right-hand side of (\ref{finupper}).

We next establish the lower bound. Truncating the outer integral at
$L_A/2$ and using (\ref{lowdelta}), we get, for some $D$ whose absolute
value is bounded by the right-hand side of (\ref{Eterm}) with $L_A$ in
place of $K$,
%
%e65 ###
%
%e65 #&#
\begin{eqnarray}\label{bigneweq}
E[R] &\geq&(1 + D) \int_{L_A/2}^{L_A} p_t(x,y) \sqrt{\frac
{2}{\kappa\pi
}}\nonumber\\
&&\hphantom{(1 + D) \int_{L_A/2}^{L_A}}{}\times \biggl( \int_{L_A - y + \mu\kappa}^{\infty} e^{-z^2/2 \kappa} \, dz
- \int_{y - \mu\kappa}^{\infty} e^{-z^2/2 \kappa} \, dz \biggr) \, dy
\nonumber\\
&=& (1 + D)\sqrt{ \frac{2}{\kappa\pi}} \cdot\frac{2}{L_A} e^{(1 -
\mu
^2/2 - \pi^2/2L_A^2)t} e^{\mu x} \sin\biggl( \frac{\pi x}{L_A} \biggr)
\nonumber\\
&&{} \times\int_{L_A/2}^{L_A} e^{-\mu y} \sin\biggl( \frac
{\pi y}{L_A} \biggr)\biggl ( \int_{L_A - y + \mu\kappa}^{\infty}
e^{-z^2/2 \kappa} \, dz - \int_{y - \mu\kappa}^{\infty} e^{-z^2/2
\kappa} \, dz \biggr) \, dy \nonumber
\\[-8pt]
\\[-8pt]
&=& (1 + D) \sqrt{ \frac{2}{\kappa\pi}} \cdot\frac{2}{L_A} e^{(1
- \mu
^2/2 - \pi^2/2L_A^2)t} e^{\mu x} \sin\biggl( \frac{\pi x}{L_A} \biggr)
\nonumber\\
&&{} \times\biggl( \int_0^{L_A} \int_{L_A - y + \mu\kappa
}^{\infty} e^{-\mu y} \sin\biggl( \frac{\pi y}{L_A} \biggr) e^{-z^2/2
\kappa} \, dz \, dy \nonumber\\
&&\hphantom{{} \times\biggl(}{} - \int_{L_A/2}^{L_A} e^{-\mu y} \sin\biggl( \frac{\pi
y}{L_A} \biggr)\biggl ( \int_{y - \mu\kappa}^{\infty} e^{-z^2/2\kappa}
\, dz \biggr) \, dy \nonumber\\
&&\hphantom{{} \times\biggl(}{} - \int_0^{L_A/2} \int_{L_A - y + \mu\kappa}^{\infty}
e^{-\mu y} \sin\biggl( \frac{\pi y}{L_A} \biggr) e^{-z^2/2 \kappa} \, dz
\, dy \biggr).
\nonumber
\end{eqnarray}
To bound the second term in (\ref{bigneweq}), note that by substituting
$w = z/\sqrt{\kappa}$ and using the fact that $\int_x^{\infty}
e^{-w^2/2} \, dw \leq x^{-1} e^{-x^2/2}$, we get
\[
\int_{\fracd{1}{2}L_A - \mu\kappa}^{\infty} e^{-z^2/2 \kappa} \, dz
\leq\frac{2 \sqrt{\kappa}}{L_A - 2 \mu\kappa} e^{-(L_A - 2 \mu
\kappa
)^2/8 \kappa}.
\]
Therefore,
%
%e66 ###
%
%e66 #&#
\begin{eqnarray}\label{newtm1}
&&\int_{L_A/2}^{L_A} e^{-\mu y} \sin\biggl( \frac{\pi y}{L_A} \biggr)
\biggl( \int_{y - \mu\kappa}^{\infty} e^{-z^2/2 \kappa} \, dz \biggr) \,
dy \nonumber\\
&& \qquad \leq\int_{L_A/2}^{L_A} e^{-\mu y} \sin\biggl( \frac{\pi y}{L_A}
\biggr)\biggl ( \int_{\fracd{1}{2}L_A - \mu\kappa}^{\infty} e^{-z^2/2
\kappa} \, dz \biggr) \, dy\nonumber
\\[-8pt]
\\[-8pt]
&& \qquad  \leq\frac{2 \sqrt{\kappa}}{L_A - 2 \mu\kappa} e^{-(L_A
- 2 \mu\kappa)^2/8 \kappa} \int_{L_A/2}^{L_A} e^{-\mu y} \sin\biggl(
\frac{\pi y}{L_A} \biggr) \, dy \nonumber\\
&& \qquad =\biggl ( \frac{e^{-\mu L_A}}{L_A} \biggr) o(1).
\nonumber
\end{eqnarray}
To bound the third term in (\ref{bigneweq}), note that
%
%e67 ###
%
%e67 #&#
\begin{eqnarray}\label{newtm2}
&&\int_0^{L_A/2} \int_{L_A - y + \mu\kappa}^{\infty} e^{-\mu y}
\sin
\biggl( \frac{\pi y}{L_A} \biggr) e^{-z^2/2 \kappa} \, dz \, dy
\nonumber\\
&& \qquad= e^{-\mu L_A} \int_{L_A/2}^{L_A} \int_{y + \mu\kappa
}^{\infty}
e^{\mu y} \sin\biggl( \frac{\pi y}{L_A} \biggr) e^{-z^2/2 \kappa} \, dz
\, dy \nonumber
\\[-8pt]
\\[-8pt]
&& \qquad= e^{-\mu L_A} \int_{L_A/2}^{\infty} \int_{L_A/2}^{\min\{
z - \mu\kappa
, L_A\}} e^{\mu y} \sin\biggl( \frac{\pi y}{L_A} \biggr) e^{-z^2/2 \kappa
} \, dy \, dz \nonumber\\
&& \qquad\leq\frac{\pi e^{-\mu L_A}}{L_A} \int_{L_A/2}^{\infty}
z^2 e^{\mu z}
e^{-z^2/2 \kappa} \, dz = \biggl( \frac{e^{-\mu L_A}}{L_A} \biggr) o(1).
\nonumber
\end{eqnarray}
For the first term, we argue as in the proof of the upper bound, and
then use that $e^{\mu y} \geq1$ and $\sin y \geq y - y^3/6$ for all $y
\geq0$ to get
\begin{eqnarray*}
&&\int_0^{L_A} \int_{L_A - y + \mu\kappa}^{\infty} e^{-\mu y} \sin
\biggl( \frac{\pi y}{L_A} \biggr) e^{-z^2/2 \kappa} \, dz \, dy \\
&& \quad= e^{-\mu L_A} \int_0^{\infty}
\int_0^{\min\{z - \mu\kappa, L_A\}} e^{\mu y} \sin\biggl( \frac{\pi
y}{L_A} \biggr) e^{-z^2/2 \kappa} \, dy \, dz \\
&& \quad\geq\frac{\pi e^{-\mu L_A}}{L_A} \int_0^{\infty}
\int_0^{\min\{z - \mu\kappa, L_A\}}
y e^{-z^2/2 \kappa} \, dy \, dz\\
&&\qquad{}-\frac{\pi^3 e^{-\mu L_A}}{6 L_A^3}
\int_0^{\infty} \int_0^z y^3 e^{-z^2/2 \kappa} \, dy \, dz \\
&& \quad= \frac{\pi e^{-\mu L_A}}{2 L_A} \int_0^{\infty}
\min\{z - \mu\kappa, L_A\}^2 e^{-z^2/2 \kappa} \, dz
- \frac{\pi^3
e^{-\mu L_A}}{24 L_A^3} \int_0^{\infty}
z^4 e^{-z^2/2 \kappa} \, dz.
\end{eqnarray*}
The second integral is a constant times $\kappa^{5/2}$. The first
integral would be $\kappa^{3/2} \sqrt{\pi/2}$ if we had $z^2$ in the
integrand in place of $\min\{z - \mu\kappa, L_A\}^2$. Also,
\[
\int_0^{\infty} (z^2 - \min\{z - \mu\kappa, L_A\}^2 )
e^{-z^2/2 \kappa} \, dz \leq\int_0^{\infty} \max\bigl\{2 \mu\kappa z, z^2
\mathbf{1}_{\{z \geq L_A\}} \bigr\} e^{-z^2/2\kappa} \, dz.
\]
If we use $2 \mu\kappa z$ in the integrand, the integral is bounded by
$C \kappa^2$. If we use $z^2 \mathbf{1}_{\{z \geq L_A\}}$, the integral
divided by $\kappa^{3/2}$ tends to zero uniformly over $\kappa\in
(0,1)$ as $N \rightarrow\infty$, so the integral is $\kappa^{3/2}
o(1)$. These observations, combined with (\ref{bigneweq}) and the
bounds in (\ref{newtm1}) and (\ref{newtm2}), imply that $E[R]$ is
bounded below by the right-hand side of (\ref{finupper}).
\end{pf}

%
%pr16 #&#
\begin{Prop}\label{EMkProp}
We have
\begin{eqnarray*}
E[R_k|{\cal F}_{t_{k-1}}] &=& 2 \pi e^{A} \cdot\frac{Z_{N,1}'(t_{k-1})
\theta s}{N (\log N)^2} \bigl(1 + O(|A| \theta) + o(1)\bigr) \\
&&{}+ \frac{C e^A
Y_N(t_{k-1})(1 + o(1))}{N (\log N)^3}.
\end{eqnarray*}
On $G_{N, k-1}$, we have
\[
E[R_k|{\cal F}_{t_{k-1}}] = 2 \pi e^{A} \cdot\frac{Z_N(t_{k-1})
\theta
s}{N (\log N)^2} \bigl(1 + O(|A| \theta)\bigr) + o(1).
\]
\end{Prop}

\begin{pf}
We first consider the particles that reach $L_A$ between times
$t_{k-1}$ and $t_{k-1} + (\log N)^2$. Define $R_k(t)$ in the same
manner as $R_k$, but counting only particles that reach $L_A$ between
times $t_{k-1}$ and $t$. Let $R_{k,1} = R_k(t_{k-1} + (\log N)^2)$. We
now consider the martingale from Lemma~\ref{hhlem}. Since $\mu^2/2 - 1
< 0$, this process will still be a supermartingale if particles are
stopped, but not killed, when reaching $L_A$.
More precisely, for $t_{k-1} \leq t \leq t_k$, let $X_{i,N}^{L_A}(t) =
X_{i,N}(t)$ if, for all\vspace*{-1pt} $u \in[t_{k-1}, t)$, the ancestor at time $u$
of the individual $X_{i,N}(t)$ is in $(0, L_A)$, and let
$X_{i,N}^{L_A}(t) = 0$ otherwise, and then for $t_{k-1} \leq t \leq
t_k$, define
\[
V_A(t) = R_k(t) L_A e^{\mu L_A + (\mu^2/2 - 1)(t - t_{k-1})} + \sum
_{i=1}^{M_N(t)} X_{i,N}^{L_A}(t) e^{\mu X_{i,N}^{L_A}(t) + (\mu^2/2 -
1)(t - t_{k-1})}.
\]
Then $(V_A(t), t_{k-1} \leq t \leq t_k)$ is a supermartingale with
respect to $({\cal F}_{t}, t_{k-1} \leq t \leq t_k)$. Therefore,
\begin{eqnarray*}
V_A(t_{k-1}) &\geq& E\bigl[V_A\bigl(t_{k-1} + (\log N)^2\bigr)|{\cal F}_{t_{k-1}}\bigr]
\\
&\geq& E\bigl[R_{k,1} L_A e^{\mu L_A + (\mu^2/2 - 1) (\log N)^2}|{\cal
F}_{t_{k-1}}\bigr] \\
&=& L_A e^{\mu L_A + (\mu^2/2 - 1)(\log N)^2} E[R_{k,1}|{\cal
F}_{t_{k-1}}].
\end{eqnarray*}
Note that since $X_{i,N}^{L_A}(t) \leq L_A$, we have $V_A(t_{k-1}) \leq
L_A Y_N(t_{k-1})$, which means
\[
E[R_{k,1}|{\cal F}_{t_{k-1}}] \leq\frac{Y_N(t_{k-1})}{e^{\mu L_A +
(\mu
^2/2 - 1)(\log N)^2}}.
\]
Since $(1 - \mu^2/2) \leq C/(\log N)^2$ and $e^{\mu L_A} \geq N (\log
N)^3 e^{-A} (1 + o(1))$, we have
%
%e68 ###
%
%e68 #&#
\begin{equation}\label{Mk1}
E[R_{k,1}|{\cal F}_{t_{k-1}}] \leq\frac{C e^A Y_N(t_{k-1}) (1 +
o(1))}{N (\log N)^3}.
\end{equation}

We next consider the particles that reach $L_A$ between times $t_{k-1}
+ (\log N)^2$ and $t_k$. The strategy will be to choose a small number
$\delta$, break the time interval $[t_{k-1} + (\log N)^2, t_k]$ into
time intervals of length $\delta$ and then use Lemma~\ref{tdeltalem} to
estimate the number of particles that reach $L_A$ in each of these
intervals. We first make three remarks concerning the application of
Lemma~\ref{tdeltalem}. First, note that if the interval starts at time
$t$, then $t - t_{k-1}$ plays the role of $t$ in Lemma~\ref{tdeltalem}.
Since $t - t_{k-1} \leq(\log N)^3 \theta s$, equation (\ref{LAexp})
implies that
%
%e69 ###
%
%e69 #&#
\begin{equation}\label{exptheta}
\bigl|e^{(1 - \mu^2/2 - \pi^2/2L_A^2)(t - t_{k-1})} - 1\bigr| \leq C |A| \theta
+ o(1).
\end{equation}
Second, we need to consider all particles at time $t$ rather than just
a single particle at $x$, so in place of $e^{\mu x} \sin(\pi x/L_A)$,
we have the expression $Z_{N,1}'(t_{k-1})$ from (\ref{ZNprime}). Third,
note that by (\ref{Eterm}) with $K = L_A$, the error term $|D|$ is
bounded by $C(1 + o(1))$ for $t \geq t_{k-1} + (\log N)^2$, and $|D|$
is $o(1)$ for $t \geq t_{k-1} + (\log N)^{5/2}$.

Let $R_{k,2}$ be defined in the same way as $R$, but counting only
particles that reach $L_A$ between times $t_{k-1} + (\log N)^2$ and
$t_{k-1} + (\log N)^{5/2}$. We can divide this time interval into at
most $\delta^{-1} (\log N)^{5/2}$ time intervals of length $\delta$, so
by Lemma~\ref{tdeltalem},
%
%e70 ###
%
%e70 #&#
\begin{eqnarray}\label{Mk2}
E[R_{k,2}|{\cal F}_{t_{k-1}}] &\leq& C e^A \delta\cdot\frac
{Z_{N,1}'(t_{k-1})}{N (\log N)^5} \cdot\frac{(\log N)^{5/2}}{\delta}
\cdot\bigl(1 + o(1)\bigr) \nonumber
\\[-8pt]
\\[-8pt]&\leq&\frac{C e^A Z_{N,1}'(t_{k-1})(1 + o(1))}{N (\log
N)^{5/2}}.
\nonumber
\end{eqnarray}

Let $R_{k,3}$ be defined in the same way as $R$, but counting only
particles that reach $L_A$ between times $t_{k-1} + (\log N)^{5/2}$ and
$t_k$. This interval can be divided into $\delta^{-1} (\log N)^3
\theta
s (1 + o(1))$ intervals of length $\delta$, so by Lemma~\ref{tdeltalem}
and (\ref{exptheta}),
%
%e71 ###
%
%e71 #&#
\begin{eqnarray}\label{Mk3}
E[R_{k,3}|{\cal F}_{t_{k-1}}] &=& 2 \pi e^A \delta\bigl(1 + O(|A| \theta)\bigr)
\nonumber\\
&&{}\times\frac{Z_{N,1}'(t_{k-1})(1 + o(1))(1 + o(\delta))}{N (\log N)^5}
\cdot\frac{(\log N)^3 \theta s}{\delta} \\
&=& 2 \pi e^A \frac{Z_{N,1}'(t_{k-1}) \theta s}{N (\log N)^2} \bigl(1 + O(|A|
\theta)\bigr)\bigl(1 + o(\delta)\bigr)\bigl(1 + o(1)\bigr) + o(1).
\nonumber
\end{eqnarray}
The first statement of the proposition follows from (\ref{Mk1}), (\ref
{Mk2}) and (\ref{Mk3}) by choosing $\delta$ as a function of $N$ so
that $\delta\rightarrow0$ as $N \rightarrow\infty$. On $G_{N, k-1}$,
the second statement follows from Lemma~\ref{Zprimelem}.
\end{pf}

The corollary below follows immediately from the above proof, because
the number of particles that hit $L_A$ between $t_k - (\log N)^{5/2}$
and $t_k$ can be bounded in the same manner as $E[R_{k,3}|{\cal
F}_{t_{k-1}}]$, and the number of intervals of length $\delta$ is only
$(\log N)^{5/2}/\delta$.

%
%co17 #&#
\begin{Cor}\label{EMkCor}
Define ${\tilde R}_k$ the same way as $R_k$, except only counting
particles that reach $L_A$ between $t_k - (\log N)^{5/2}$ and $t_k$.
Then $E[{\tilde R}_k|{\cal F}_{t_{k-1}}]$ is $o(1)$ on $G_{N, k-1}$.
\end{Cor}

%
%pr18 #&#
\begin{Prop}\label{EMk2Prop}
Assume $A \geq0$. On $G_{N, k-1}$, we have
\[
E[R_k^2|{\cal F}_{t_{k-1}}] \leq\frac{C \theta e^A Z_N(t_{k-1})}{N
(\log N)^2} + o(1).
\]
\end{Prop}

\begin{pf}
For the purposes of this proof, we may assume that particles are killed
upon reaching $L_A$. Note that $R_k^2 = R_k + 2Y$, where $Y$ is the
number of distinct pairs of particles that get killed upon reaching
$L_A$. We may further write $Y = Y_1 + Y_2$, where $Y_1$ denotes the
number of pairs of particles that get killed upon reaching $L_A$ whose
most recent common ancestor is before time $t_{k-1}$ and $Y_2$ counts
the other pairs of particles.
Proposition~\ref{EMkProp} and (\ref{con1}) give that on $G_{N, k-1}$,
%
%e72 ###
%
%e72 #&#
\begin{equation}\label{newEMk}
E[R_k|{\cal F}_{t_{k-1}}] \leq\frac{C \theta e^A Z_N(t_{k-1})}{N
(\log
N)^2} + o(1).
\end{equation}

If there is a particle at $x$ at time $t_{k-1}$ and a descendant of
this particle reaches $L_A$ by time $t_k$, then the number of pairs in
$Y_1$ involving this descendant will be precisely the number of
particles descended from particles other than the particle at $x$ at
time $t_{k-1}$ that reach $L_A$ by time $t_k$, which is bounded by
$R_k$. Because descendants of different particles evolve independently,
it follows that
%
%e73 ###
%
%e73 #&#
\begin{equation}\label{condEY1}
E[Y_1|{\cal F}_{t_{k-1}}] \leq(E[R_k|{\cal F}_{t_{k-1}}])^2 \leq C
\biggl( \frac{\theta e^A Z_N(t_{k-1})}{N (\log N)^2} \biggr)^2 + o(1).
\end{equation}

It remains to consider $Y_2$. Because pairs of particles contributing
to $Y_2$ have the same ancestor at time $t_{k-1}$, we may consider
separately the contributions of the particles at time $t_{k-1}$. Assume
for now that there is a single particle at $x$ at time $t_{k-1}$, and
we denote the number of associated pairs of particles contributing to
$Y_2$ by $Y_2^x$. Let $h(t,y)$ be the expected number of offspring of a
single particle that is at $y$ at time $t_{k-1} + t$ that will hit
$L_A$ before time $t_k$. A branching event at $(t_{k-1}+t, y)$
produces, on average, $h(t,y)^2$ pairs of particles that hit $L_A$ and
have their most recent common ancestor at time $t$. Since each particle
branches at rate 1,
%
%e74 ###
%
%e74 #&#
\begin{equation}\label{Y2x}
E[Y_2^x] = \int_0^{(\log N)^3 \theta s} \int_0^{L_A} q_t(x, y) h(t,y)^2
\, dy \, dt.
\end{equation}
Since $h(t,y) \leq h(0,y)$, it follows from Proposition~\ref{EMkProp} that
%
%e75 ###
%
%e75 #&#
\begin{eqnarray}
h(t,y) &\leq&\frac{C e^A \theta e^{\mu y}}{N (\log N)^2} \sin\biggl(
\frac{\pi y}{L_A} \biggr)\bigl(1 + o(1)\bigr) + \frac{C e^A e^{\mu y} (1 +
o(1))}{N (\log N)^3} \nonumber
\\[-8pt]
\\[-8pt]
&\leq&\frac{C e^A(1 + o(1))}{N (\log N)^2} \biggl( \theta e^{\mu y} \sin
\biggl( \frac{\pi y}{L_A} \biggr) + \frac{e^{\mu y}}{\log N} \biggr),
\nonumber
\end{eqnarray}
where the $o(1)$ term tends to zero uniformly in $y$ as $N \rightarrow
\infty$.

We first evaluate the portion of the integral in (\ref{Y2x}) when $t
\leq(\log N)^2$. Recall from (\ref{qsvs}) that when $t \leq(\log
N)^2$,
\[
q_t(x, y) \leq C e^{\mu(x - y)} v_t(x,y),
\]
where $v_t(x,y)$ is the density of Brownian motion in the strip $(0,
L_A)$, defined as in (\ref{vtdef}) with $L_A$ in place of $K$.
Therefore, changing the order of integration,
\begin{eqnarray*}
&&\int_0^{(\log N)^2} \int_0^{L_A} q_t(x, y) h(t,y)^2 \, dy \, dt
\\
&& \qquad \leq\frac{C e^{2A} e^{\mu x}(1 + o(1))}{N^2 (\log N)^4}\\
&&\qquad \quad {}\times
\int_0^{L_A} \biggl( \theta e^{\mu y} \sin\biggl( \frac{\pi y}{L_A}
\biggr) + \frac{e^{\mu y}}{\log N} \biggr)^2 e^{-\mu y} \biggl( \int_0^{(\log
N)^2} v_t(x,y) \, dt \biggr) \, dy.
\end{eqnarray*}
By (\ref{GreenLz}), $\int_0^{\infty} v_t(x,y) \, dt \leq2(L_A - y)$.
Using also that $(a + b)^2 \leq C(a^2 + b^2)$, that $\sin(\pi y/L_A) =
\sin(\pi(L_A - y)/L_A) \leq\pi(L_A - y)/L_A \leq C (L_A - y)/(\log
N)$, that $e^{\mu L_A} = N (\log N)^3 e^{-A} (1 + o(1))$, and that
(\ref
{con0}) holds, we get
%
%e76 ###
%
%e76 #&#
\begin{eqnarray}\label{intln2}
&&\int_0^{(\log N)^2} \int_0^{L_A} q_t(x, y) h(t,y)^2 \, dy \, dt
\nonumber\\
&& \qquad\leq\frac{C e^{2A} e^{\mu x}(1 + o(1))}{N^2 (\log N)^4}
\int_0^{L_A}
\biggl( \frac{\theta^2 e^{\mu y} (L_A - y)^3}{(\log N)^2} + \frac{e^{\mu
y} (L_A - y)}{(\log N)^2} \biggr) \, dy \\
&& \qquad\leq\frac{C e^{2A} e^{\mu x}(1 + o(1))}{N^2 (\log N)^6} (
\theta
^2 e^{\mu L_A} + e^{\mu L_A} ) \leq\frac{C e^A e^{\mu x}(1 +
o(1))}{N (\log N)^3}.\nonumber
\end{eqnarray}

When $t \geq(\log N)^2$, Lemma~\ref{stripdensity} implies that
\[
q_t(x,y) \leq C p_t(x,y) \bigl(1 + o(1)\bigr) \leq\frac{C(1 + o(1))}{\log N}
e^{\mu x} \sin\biggl( \frac{\pi x}{L_A} \biggr) e^{-\mu y} \sin\biggl(
\frac{\pi y}{L_A} \biggr)
\]
because $A \geq0$. Therefore,
%
%e77 ###
%
%e77 #&#
\begin{eqnarray}\label{intln3}
&&\int_{(\log N)^2}^{(\log N)^3 \theta s} \int_0^{L_A} q_t(x, y)
h(t,y)^2 \, dy \, dt \nonumber\\
&& \qquad\leq\frac{C e^{2A}(1 + o(1))}{N^2 (\log N)^5} e^{\mu x}
\sin\biggl(
\frac{\pi x}{L_A} \biggr) \nonumber\\
&& \qquad \quad{} \times\int_{(\log N)^2}^{(\log N)^3 \theta s}
\int
_0^{L_A} e^{-\mu y} \sin\biggl( \frac{\pi y}{L_A} \biggr)\biggl ( \theta
e^{\mu y} \sin\biggl( \frac{\pi y}{L_A} \biggr) + \frac{e^{\mu y}}{\log
N} \biggr)^2 \, dy \, dt \nonumber
\\[-8pt]
\\[-8pt]
&& \qquad\leq\frac{C \theta e^{2A}(1 + o(1))}{N^2 (\log N)^2} e^{\mu
x} \sin
\biggl( \frac{\pi x}{L_A} \biggr)\nonumber
\\
&& \qquad  \quad {}\times \int_0^{L_A} \biggl( \frac{\theta^2
e^{\mu y} (L_A - y)^3}{(\log N)^3} + \frac{e^{\mu y} (L_A - y)}{(\log
N)^3} \biggr) \, dy \nonumber\\
&& \qquad\leq\frac{C \theta e^A(1 + o(1))}{N (\log N)^2} \cdot
e^{\mu x} \sin
\biggl( \frac{\pi x}{L_A} \biggr),\nonumber
\end{eqnarray}
using that the last integral can be bounded by $C e^{\mu L_A}/(\log
N)^3 \leq C e^{-A} N$. By combining (\ref{intln2}) and (\ref{intln3})
and summing over the contributions from different particles, we get on
$G_{N, k-1}$,
%
%e78 ###
%
%e78 #&#
\begin{eqnarray}\label{condEY2}
E[Y_2|{\cal F}_{t_{k-1}}] &\leq&\biggl( \frac{C e^A Y_N(t_{k-1})}{N (\log
N)^3} + \frac{C \theta e^A Z_{N,1}'(t_{k-1})}{N (\log N)^2} \biggr)\bigl(1 +
o(1)\bigr) \nonumber
\\[-8pt]
\\[-8pt]&\leq&\frac{C \theta e^A Z_N(t_{k-1})}{N (\log N)^2} + o(1),
\nonumber
\end{eqnarray}
where the last inequality uses Lemma~\ref{Zprimelem}.
The result now follows from (\ref{newEMk}), (\ref{condEY1}) and (\ref
{condEY2}), using the assumption from (\ref{con3}) that $\theta e^A
\eps
^{-1/2} \leq1$.
\end{pf}

%
%re19 #&#
\begin{Rmk}
By Proposition~\ref{EMkProp} and Markov's inequality, we see that when
$\theta$ is small, during most intervals $[t_{k-1}, t_k]$, no particles
reach $L_A$. Using Propositions~\ref{EMkProp} and~\ref{EMk2Prop} and
the second moment method, we get that on~$G_{N, k-1}$,
%
%e79 ###
%
%e79 #&#
\begin{eqnarray}
P(R_k > 0|{\cal F}_{t_{k-1}}) &\geq&\frac{(E[R_k|{\cal
F}_{t_{k-1}}])^2}{E[R_k^2|{\cal F}_{t_{k-1}}]} \nonumber\\
&\geq& C\biggl ( \frac{\theta e^A Z_N(t_{k-1})}{N (\log N)^2} \biggr)^2
\biggl( \frac{N (\log N)^2}{\theta e^A Z_N(t_{k-1})} \biggr)\bigl(1 + o(1)\bigr)
\\
&\geq&
\frac{C \theta e^A Z_N(t_{k-1})(1 + o(1))}{N (\log N)^2}.\nonumber
\end{eqnarray}
Thus, it follows that
\[
E[R_k|R_k > 0, {\cal F}_{t_{k-1}}] = \frac{E[R_k|{\cal
F}_{t_{k-1}}]}{P(R_k > 0|{\cal F}_{t_{k-1}})} \leq C\bigl(1 + o(1)\bigr).
\]
That is, conditional on the event that at least one particle reaches
$L_A$, the expected number of particles that reach $L_A$ is bounded by
a constant that does not depend on
$\theta$ or $A$.
\end{Rmk}

%s3.3 ###
%s3.3 #&#
\subsection{\texorpdfstring{The probability of $G_N(\eps)$}{The probability of GN(epsilon)}}\label{sec33}

We have now acquired enough tools to prove that the probability of
$G_N(\eps)$ is close to 1 when $\eps$ is small and $N$ is large. It is
this result that allows us to work on the event $G_N$ throughout much
of the paper.

Throughout this section, we will assume that particles are killed upon
reaching $L_A$. Define ${\bar Z}_N(t)$ and ${\bar Y}_N(t)$ in the same
way as $Z_N(t)$ and $Y_N(t)$ in (\ref{ZNdef}) and (\ref{YNdef}), but
for this modified process in which particles are killed upon reaching
$L_A$. Also, we use $L_A$ rather than $L$ in the definition of ${\bar Z}_N(t)$.

%
%le20 #&#
\begin{Lemma}\label{YNbar}
For any fixed $A \in\R$ and any fixed $t \geq0$, under the hypotheses
of Proposition~\ref{ZNCSBP}, we have
\[
\lim_{N \rightarrow\infty} P\bigl({\bar Y}_N(t (\log N)^3) > N (\log N)^3
h(N)\bigr) = 0.
\]
\end{Lemma}

\begin{pf}
Let $B_N$ be the event that all particles at time zero are in $(0,
L_A)$. We have $P(B_N) \rightarrow1$ as $N \rightarrow\infty$ because
on the event that there is a particle to the right of $L_A$ at time
zero, we have $Y_N(0) \geq N (\log N)^3 h(N)$ by (\ref{hNbound}),
and $Y_N(0)/(N (\log N)^3 h(N))$ converges in probability to zero by
the definition of $h(N)$. On $B_N$, we have $Y_N(0) = {\bar Y}_N(0)$,
so the result holds when $t = 0$.

Suppose instead $t > 0$. By (\ref{Yexp}), on $B_N$ we have
\[
E[{\bar Y}_N(t (\log N)^3)|{\cal F}_0] = \frac{4}{\pi} e^{(1 - \mu^2/2
- \pi^2/2L_A^2)t (\log N)^3}(1 + D) Z_N(0),
\]
where $|D|$ is bounded by the right-hand side of (\ref{Eterm}) with
$L_A$ in place of $K$. Therefore, using (\ref{LAexp}), on $B_N$ we have
\[
E[{\bar Y}_N(t (\log N)^3)|{\cal F}_0] \leq C e^{-2 \pi^2 A t (1 +
o(1))} Z_N(0)\bigl(1+o(1)\bigr).
\]
Therefore, by Markov's inequality, on $B_N$ we have
\begin{eqnarray*}
P\bigl({\bar Y}_N(t(\log N)^3) > N (\log N)^3 h(N)|{\cal F}_0\bigr) &\leq&\frac
{E[{\bar Y}_N(t (\log N)^3)|{\cal F}_0]}{N (\log N)^3 h(N)} \\
&\leq&\frac{C e^{-2 \pi^2 A t(1 + o(1))} Z_N(0)(1+o(1))}{N (\log N)^3
h(N)}.
\end{eqnarray*}
The right-hand side converges in probability to zero because $Z_N(0)/N
(\log N)^2$ converges in distribution to $\nu$ and $(\log N) h(N)
\rightarrow\infty$ as $N \rightarrow\infty$. The result follows.
\end{pf}

For the rest of Section~\ref{sec33}, we will assume that $A < 0$, so
that $L_A > L$.

%
%le21 #&#
\begin{Lemma}\label{Zbar}
Fix $A < 0$ and $t \geq0$. For all $\kappa> 0$, there exists a
positive constant $C_1$, depending on $\kappa$ but not on $A$ or $t$,
such that under the hypotheses of Proposition~\ref{ZNCSBP},
\[
P  \biggl( \max_{0 \leq r \leq t (\log N)^3} {\bar Z}_N(r) > \frac
{1}{2}\eps^{-1/2} N (\log N)^2  \biggr) \leq\kappa+ C_1 \eps^{1/2}
e^{-2 \pi^2 A t (1 + o(1))}.
\]
\end{Lemma}

\begin{pf}
By (\ref{Zexp}), if $U(r) = e^{-(1 - \mu^2/2 - \pi^2/2L_A^2)s} {\bar
Z}_N(r)$, then $(U(r), r \geq0)$ is a martingale. Since $A < 0$, we
have $1 - \mu^2/2 - \pi^2/2L_A^2 > 0$. Therefore, by Doob's maximal
inequality (see, e.g., the $p = 1$ case of Theorem 1.4 in~\cite{cw})
and (\ref{LAexp}),
\begin{eqnarray*}
&&P \biggl( \max_{0 \leq r \leq t (\log N)^3} {\bar Z}_N(r) > \frac{1}{2}
\eps^{-1/2} N (\log N)^2 \Big| {\cal F}_0 \biggr) \\
&& \qquad\leq P \biggl( \max_{0 \leq s \leq t (\log N)^3} U(r) > \frac
{1}{2} \eps
^{-1/2} N (\log N)^2 e^{-(1 - \mu^2/2 - \pi^2/2L_A^2)t(\log N)^3}
\Big| {\cal F}_0 \biggr) \\
&& \qquad\leq\frac{2{\bar Z}_N(0) \eps^{1/2} e^{(1 - \mu^2/2 - \pi
^2/2L_A^2)t(\log N)^3}}{N (\log N)^2} \leq\frac{2{\bar Z}_N(0) \eps
^{1/2} e^{-2 \pi^2 A t(1 + o(1))}}{N (\log N)^2}.
\end{eqnarray*}
Because the distribution of $Z_N(0)/N (\log N)^2$, and therefore that
of ${\bar Z}_N(0)/\break N (\log N)^2$, converges to $\nu$, there exists a
constant $C_1$ such that $P({\bar Z}_N(0)/\break N (\log N)^2 > C_1) \leq
\kappa$. The result follows.
\end{pf}

%
%le22 #&#
\begin{Lemma}\label{hitLA}
Let $\kappa> 0$ and $t > 0$. Under the hypotheses of Proposition \ref
{ZNCSBP}, there exist positive constants $C_2$ and $\gamma$, depending
on $\kappa$ and $t$, such that for all $A < 0$, the probability that
some particle reaches $L_A$ before time $t (\log N)^3$ is at most $C_2
e^{\gamma A} + \kappa$ for sufficiently large $N$.
\end{Lemma}

\begin{pf}
Let $J = \lceil4 \pi^2 t \rceil$. For $1 \leq j \leq J$, let $A_j =
2^{j-J} A$ and $s_j =\break (j/4 \pi^2)\* (\log N)^3$. Let $s_0 = 0$. Consider
a modified branching Brownian motion, defined up to time $t (\log
N)^3$, in which particles that reach $L_{A_j}$ between times $s_{j-1}$
and $s_j$ are killed. Because $L_{A_j} \leq L_A$ for all $j$, it
suffices to bound the probability that at least one particle gets
killed in this new modified branching Brownian motion.

Let ${\tilde M}_N(r)$ denote the number of particles alive at time $r$,
and denote the positions of these particles by ${\tilde X}_{1,N}(r)
\geq\cdots\geq{\tilde X}_{M_N(r), N}(r)$. For $r \in[s_{j-1}, s_j]$,
define
\[
{\tilde Z}_{N,j}(r) = \sum_{i=1}^{{\tilde M}_N(r)} e^{\mu{\tilde
X}_{i,N}(r)} \sin\biggl( \frac{\pi{\tilde X}_{i,N}(r)}{L_{A_j}} \biggr)
\]
and
\[
{\tilde Y}_N(r) = \sum_{i=1}^{{\tilde M}_N(r)} e^{\mu{\tilde X}_{i,N}(r)}.
\]
For all $r \in[s_{j-1}, s_j]$, we have, using (\ref{Zexp}) and (\ref{LAexp}),
%
%e80 ###
%
%e80 #&#
\begin{eqnarray}\label{ZNj}
E[{\tilde Z}_{N,j}(r)|{\cal F}_{s_{j-1}}] &=& e^{(1 - \mu^2/2 - \pi
^2/2L_{A_j}^2)(r - s_{j-1})} {\tilde Z}_{N,j}(s_{j-1}) \nonumber\\
&\leq& e^{- 2 \pi^2 A_j (s_j - s_{j-1})(1 + o(1))/(\log N)^3} {\tilde
Z}_{N,j}(s_{j-1}) \\
&=& e^{- A_j(1 + o(1))/2} {\tilde Z}_{N,j}(s_{j-1}).\nonumber
\end{eqnarray}

This bound allows us to bound the probability that a particle reaches
${\tilde L}_{A_j}$ between times $s_{j-1}$ and $s_j$ using Proposition
\ref{EMkProp}. We divide the interval from $s_{j-1}$ to $s_j$ into
smaller subintervals of length approximately $\theta(\log N)^3$.
More precisely, define times $s_{j-1} = u_0 < u_1 <\cdots< u_D = s_j$
such that $\theta(\log N)^3 \leq u_k - u_{k-1} \leq2 \theta(\log N)^3$
for all $k$, which is possible as long as we choose $\theta\leq1/4
\pi
^2$. Letting ${\tilde R}_{j,k}$ denote the number of particles that
reach $L_{A_j}$ between times $u_{k-1}$ and $u_k$, we get from
Proposition~\ref{EMkProp},
\[
E[{\tilde R}_{j,k}|{\cal F}_{u_{k-1}}] \leq\frac{C e^{A_j} {\tilde
Z}_{N, j}(u_{k-1}) \theta(1 + o(1))}{N (\log N)^2} + \frac{C e^{A_j}
{\tilde Y}_N(u_{k-1})(1 + o(1))}{N (\log N)^3}.
\]
By Markov's inequality,
\begin{eqnarray*}
P({\tilde R}_{j,k} > 0|{\cal F}_{u_{k-1}}) &\leq&\frac{C e^{A_j}
{\tilde
Z}_{N,j}(u_{k-1}) \theta(1 + o(1))}{N (\log N)^2}\\
&&{} + \min\biggl\{ \frac
{C e^{A_j} {\tilde Y}_N(u_{k-1})(1 + o(1))}{N (\log N)^3}, 1 \biggr\}.
\end{eqnarray*}
By Lemma~\ref{YNbar} applied at the times $u_0,\ldots, u_D$, the second
term is $o_p(1)$, which means it tends to zero in probability as $N
\rightarrow\infty$ for any fixed values of the parameters $A$ and
$\theta$.
Therefore,
\begin{eqnarray*}
P({\tilde R}_{j,k} > 0|{\cal F}_{s_{j-1}}) &=& E[P({\tilde R}_{j,k} >
0|{\cal F}_{u_{k-1}})|{\cal F}_{s_{j-1}}] \\
&\leq&\frac{C e^{A_j} E[{\tilde Z}_{N,j}(u_{k-1})|{\cal F}_{s_{j-1}}]
\theta(1 + o(1))}{N (\log N)^2} + o_p(1).
\end{eqnarray*}
Let ${\tilde R}_j = \sum_{k=1}^D {\tilde R}_{j,k}$. By (\ref{ZNj}) and
the fact that $D \leq C/\theta$,
%
%e81 ###
%
%e81 #&#
\begin{equation}\label{tRj}\qquad
P({\tilde R}_j > 0|{\cal F}_{s_{j-1}}) \leq\frac{C e^{A_j(1 + o(1))/2}
{\tilde Z}_{N,j}(s_{j-1})(1 + o(1))}{N (\log N)^2} + o_p(1).\vadjust{\goodbreak}
\end{equation}

Let ${\tilde G}_j$ be the event that ${\tilde Y}_N(s_k) \leq N (\log
N)^3 h(N)$ for $k = 0,\ldots, j$, and let ${\tilde G} = {\tilde
G}_{J-1}$. We have $P({\tilde G}) = 1 - o(1)$ by Lemma~\ref{YNbar}.
We now\vspace*{1pt} show by induction that for $j = 1,\ldots, J$, on ${\tilde
G}_{j-1}$ we have
%
%e82 ###
%
%e82 #&#
\begin{equation}\label{indNTS}\qquad
E[{\tilde Z}_{N,j}(s_j)|{\cal F}_{s_0}] \leq C e^{-(A_1 +\cdots+
A_j)(1 + o(1))/2} Z_N(0) + o(N (\log N)^2).
\end{equation}
The $j = 1$ case follows from (\ref{ZNj}) and the fact that $|{\tilde
Z}_{N,1}(0) - Z_N(0)| \leq C |A_1| {\tilde Y}_N(0)/(\log N)$ by (\ref
{ZNYbound}), the difference between the expressions coming from the
fact that $L_{A_1}$ is used in the definition of ${\tilde Z}_{N,1}(0)$
and $L$ is used in the definition of $Z_N(0)$. Suppose the result is
true for $j - 1$. On ${\tilde G}_{j-1}$, we have, using (\ref{ZNj}) and
the argument leading to (\ref{ZNYbound}),
\begin{eqnarray*}
E[{\tilde Z}_{N,j}(s_j)|{\cal F}_{s_{j-1}}] &\leq& e^{-A_j(1 + o(1))/2}
{\tilde Z}_{N,j}(s_{j-1}) \\
&\leq& e^{-A_j(1 + o(1))/2} {\tilde Z}_{N,j-1}(s_{j-1})\\
&&{} + e^{-A_j(1 +
o(1))/2} |{\tilde Z}_{N,j}(s_{j-1}) - {\tilde Z}_{N,j-1}(s_{j-1})|
\\
&\leq& e^{-A_j(1 + o(1))/2} {\tilde Z}_{N,j-1}(s_{j-1})\\
&&{} + C e^{-A_j(1 +
o(1))/2} |A_j - A_{j-1}| {\tilde Y}_N(s_{j-1})/(\log N) \\
&\leq& e^{-A_j(1 + o(1))/2} {\tilde Z}_{N,j-1}(s_{j-1}) + o(N (\log
N)^2).
\end{eqnarray*}
Taking conditional expectations with respect to ${\cal F}_{s_0}$ and
applying the induction hypothesis gives (\ref{indNTS}). The result
(\ref
{indNTS}) for all $j = 1,\ldots, J$ on ${\tilde G}_{j-1}$ follows by induction.

We now take conditional expectations with respect to ${\cal F}_{s_0}$
on both sides of (\ref{tRj}). Using that $|{\tilde Z}_{N,j}(s_{j-1}) -
{\tilde Z}_{N,j-1}(s_{j-1})| = o(N(\log N)^2)$ on ${\tilde G}_{j-1}$ as
shown above and that
\[
A_j - (A_1 +\cdots+ A_{j-1}) = 2^{j-J}A - (2^j - 2)2^{-J} A =
A_1=2^{1-J}A = A_1,
\]
we get
\begin{eqnarray*}
P({\tilde R}_j > 0|{\cal F}_{s_0}) &\leq&\frac{C e^{A_j(1 + o(1))/2}
E[{\tilde Z}_{N,j}(s_{j-1})|{\cal F}_{s_0}] (1 + o(1))}{N (\log N)^2} +
o_p(1) \\
&\leq&\frac{C e^{(A_j - (A_1 +\cdots+ A_{j-1}))(1 + o(1))/2} Z_N(0) (1
+ o(1))}{N (\log N)^2} + o_p(1)\\
&\leq&\frac{C e^{A_1(1+o(1))/2} Z_N(0)(1 + o(1))}{N (\log N)^2} +
o_p(1).
\end{eqnarray*}
Therefore, the probability, conditional on ${\cal F}_{s_0}$, that some
particle reaches $A = A_J$ by time $t (\log N)^3$ is at most
\[
\frac{C J e^{A2^{-J}(1 + o(1))} Z_N(0)(1 + o(1))}{N (\log N)^2} + o_p(1).
\]
Since $Z_N(0)/N (\log N)^2$ converges in distribution to $\nu$ as $N
\rightarrow\infty$, there is a constant $c$ such that
$P(Z_N(0)/N(\log
N)^2 > c) < \kappa/2$ for sufficiently large $N$. The result follows.
\end{pf}

%
%pr23 #&#
\begin{Prop}\label{GNProp}
Under the hypotheses of Proposition~\ref{ZNCSBP}, we have
\[
\lim_{\eps\rightarrow0} \sup_{\theta} \Bigl( \limsup_{N
\rightarrow\infty} \bigl(1 - P(G_N(\eps))\bigr) \Bigr) = 0,
\]
where the supremum is taken over all values of $\theta> 0$ such that
$\theta^{-1} \in\N$.
\end{Prop}

\begin{pf}
Let $0 < \eps< 1$, and let $\kappa> 0$. Choose $C_1$ as in Lemma \ref
{Zbar}, and choose $\gamma$ and $C_2$ as in Lemma~\ref{hitLA}, with $u
+ s$ in place of $t$. Choose $A = (\log\eps)/(8 \pi^2 (u+s)) < 0$. We
now assume that $\theta$ is small enough that assumptions (\ref
{con0})--(\ref{con3}) hold for these choices of $\eps$ and $A$, so that
previous results in this section may be applied. This assumption is
permissible because dividing $\theta$ by a positive integer to make it
small enough to satisfy these conditions can only reduce the value of
$P(G_N(\eps))$ by adding additional times at which conditions on $Y_N$
and $Z_N$ must hold.

By Lemma~\ref{hitLA}, the probability that some particle reaches $L_A$
by time $(u + s)\*(\log N)^3$ is at most
%
%e83 ###
%
%e83 #&#
\begin{equation}\label{check1}
\kappa+ C_2 e^{\gamma(\log\eps)/(8 \pi^2 (u+s))}
\end{equation}
for sufficiently large $N$. By Lemma~\ref{YNbar},
%
%e84 ###
%
%e84 #&#
\begin{equation}\label{check2}
\lim_{N \rightarrow\infty} P\bigl({\bar Y}_N(t_j) > N (\log N)^3 h(N)
\mbox
{ for some }j \leq\theta^{-1}\bigr) = 0.
\end{equation}
By Lemma~\ref{Zbar},
\begin{eqnarray*}
&&\limsup_{N \rightarrow\infty} P\biggl ({\bar Z}_N(t_j) > \frac{1}{2}\eps
^{-1/2} N (\log N)^2 \mbox{ for some }j \leq\theta^{-1} \biggr)\\
 && \qquad \leq
\kappa+ C_1 \eps^{1/2} e^{- (\log\eps)(1 + o(1))/4} \\
&& \qquad \leq\kappa+ C_1 \eps^{(1 + o(1))/4}.
\end{eqnarray*}
Using (\ref{ZNYbound}), on the event that ${\bar Y}_N(t_j) \leq N
(\log
N)^3 h(N)$ and no particle reaches $L_A$ by time $(u + s)(\log N)^3$,
we have
\[
Z_N(t_j) \leq{\bar Z}_N(t_j) + \frac{\pi|A| L_A N (\log N)^3
h(N)}{\sqrt{2} L^2} \leq{\bar Z}_N(t_j) + \frac{1}{2} \eps
^{-1/2}N(\log N)^2
\]
for sufficiently large $N$. Thus,
%
%e85 ###
%
%e85 #&#
\begin{eqnarray}\label{check3}
&&\limsup_{N \rightarrow\infty} P \bigl(Z_N(t_j) > \eps^{-1/2} N(\log N)^2
\mbox{ for some }j \leq\theta^{-1} \bigr) \nonumber
\\[-8pt]
\\[-8pt]
&& \qquad\leq2 \kappa+ C_1 \eps^{1/4} + C_2 e^{\gamma(\log\eps
)/(8 \pi^2 (u
+ s))}.
\nonumber
\end{eqnarray}
As ${\bar Y}_N(t_j) = Y_N(t_j)$ and ${\bar Z}_N(t_j) = Z_N(t_j)$ when
no particles reach\vspace*{1pt} $L_A$ by time $(u + s)^3(\log N)^3$, we see that $1
- P(G_N(\eps))$ is bounded by the sum of the probabilities in (\ref
{check1}), (\ref{check2}) and (\ref{check3}). Since none of the bounds
depends on $\theta$, it follows that
\[
\limsup_{\eps\rightarrow0} \sup_{\theta} \Bigl( \limsup_{N
\rightarrow\infty} \bigl(1 - P(G_N(\eps))\bigr) \Bigr) \leq3 \kappa,
\]
and the result follows by letting $\kappa\rightarrow0$.
\end{pf}

%s4 ###
%s4 #&#
\section{Critical branching Brownian motion with killing at $-y$}\label{sec4}

Consider a branching Brownian motion with drift $-\sqrt{2}$ started
with a single particle at $0$. From this process, a modified process
can be constructed in which particles that reach $-y$ are killed. Let
$Z_y$ denote the number of particles that reach $-y$ and are killed.
Note that $Z_y$ has the same distribution as the number of particles
that hit zero in branching Brownian motion with drift $-\sqrt{2}$ and
absorption at zero, started with a single particle at $y$. Because this
process almost surely goes extinct by Theorem 1.1 of~\cite{kesten}, and
it is easy to verify that infinitely many particles will not reach the
origin within any finite time interval, we see that $Z_y$ is almost
surely finite for every $0\le y< \infty$. Furthermore, because each
particle that reaches $-x$ behaves thereafter like another particle
started at zero, the number of particles that reach $-(x + y)$
conditional on $Z_x$ is the same as the distribution of the sum of
$Z_x$ independent random variables with the same distribution as $Z_y$.
Consequently, the process $(Z_y)_{y \geq0}$ is a continuous-time
branching process, as is shown in Section 5 of~\cite{nev87}. As noted
in~\cite{nev87} and in the more recent work of Maillard \cite
{maillard}, this branching process is not in the $L \log L$ class.
However, the following proposition appears on page 238 of~\cite{nev87}.

%
%pr24 #&#
\begin{Prop}\label{nevprop}
There exists a random variable $W$ such that almost surely
\[
\lim_{y \rightarrow\infty} y e^{-\sqrt{2} y} Z_y = W.
\]
Furthermore, for all $u \in\R$, we have
%
%e86 ###
%
%e86 #&#
\begin{equation}\label{LapW}
E\bigl[e^{-e^{\sqrt{2} u} W}\bigr] = \psi(u),
\end{equation}
where $\psi\dvtx  \R\rightarrow(0, 1)$ solves Kolmogorov's equation
\[
\tfrac{1}{2} \psi'' - \sqrt{2} \psi' = \psi(1 - \psi).
\]
\end{Prop}

%
%co25 #&#
\begin{Cor}\label{nevcor}
Let $\eta> 0$. There exists $y$ such that
%
%e87 ###
%
%e87 #&#
\begin{equation}\label{neveta}
P(|y e^{-\sqrt{2} y} Z_y - W| > \eta) < \eta.
\end{equation}
Moreover, there exists $\zeta> 0$ such that if particles are killed
when they reach $-y$, the probability that any particle remains alive
after time $\zeta$ is less than~$\eta$.
\end{Cor}

\begin{pf}
Equation (\ref{neveta}) is immediate from Proposition~\ref{nevprop}.
The second statement follows from the fact that $Z_y$ is almost surely
finite, and therefore so is the time when the last remaining particle
hits $-y$.
\end{pf}

Our goal in this section is to show that $P(W > x) \sim B/x$ as $x
\rightarrow\infty$ for some constant $B$. The strategy will be to
consider the Laplace transform $E[e^{-\lambda W}]$ for small values of
$\lambda$, and then apply a Tauberian theorem. From (\ref{LapW}), we
see that this requires having asymptotic results for $\psi(u)$ as $u
\rightarrow-\infty$. Equivalently, if we define $w(x) = \psi(-x)$, then
%
%e88 ###
%
%e88 #&#
\begin{equation}\label{wdiffeq}
\tfrac{1}{2} w'' + \sqrt{2} w' + w(w - 1) = 0,
\end{equation}
and we are looking for asymptotic results for $w(x)$ as $x \rightarrow
\infty$. It is well known that
%
%e89 ###
%
%e89 #&#
\begin{equation}\label{wasym}
1 - w(x) \sim C x e^{-\sqrt{2} x};
\end{equation}
see, for example, (11) in~\cite{ls1} or (1.13) in~\cite{bram2}.
However, this result turns out to be insufficient for our purposes. The
asymptotic result that we will need is given in the proposition below.

%
%pr26 #&#
\begin{Prop}\label{KPPasym}
Suppose that $w$ is an increasing function satisfying (\ref{wdiffeq})
with $\lim_{x \rightarrow\infty} w(x) = 1$ and $\lim_{x \rightarrow
-\infty} w(x) = 0$. For all $x$, let $u(x) = 1 - w(x)$ and $v(x) =
u(x)/(x e^{-\sqrt{2} x})$. Then for all $c > 0$, we have
\[
\lim_{x \rightarrow\infty} x\bigl(v(x + c) - v(x)\bigr) = 0.
\]
\end{Prop}

\begin{pf}
Let $x > 0$. Let $(R_t, t \geq0)$ be a three-dimensional Bessel
process with $R_0 = x$. According to (2.6) of~\cite{har99}, the process
\[
X_t = v(R_t) \exp\biggl( - \int_0^t u(R_s) \, ds \biggr)
\]
is a positive local martingale, and therefore a supermartingale. Let $T
= \inf\{t\dvtx  R_t = x + c\}$. By the optional sampling theorem,
\[
v(x) = E[X_0] \geq E[X_T] = v(x + c) E \biggl[ \exp\biggl(- \int_0^T
u(R_s) \, ds \biggr) \biggr],
\]
which means
\[
v(x + c) - v(x) \leq v(x + c) \biggl(1 - E \biggl[ \exp\biggl(- \int_0^T
u(R_s) \, ds \biggr) \biggr] \biggr).
\]

Let $0 < \gamma< 1$, and let $A$ be the event that $R_t \leq\gamma x$
for some $t \leq T$. That is, $A$ is the event that the Bessel process
reaches $\gamma x$ before reaching $x + c$. By Corollary~3.4 on page
253 of~\cite{revyor}, we have
\[
P(A) = \frac{(x + c)^{-1} - x^{-1}}{(x + c)^{-1} - (\gamma x)^{-1}} =
\frac{c \gamma}{c + (1 - \gamma)x}.
\]
In view of (\ref{wasym}), there are constants $C_1$ and $C_2$ such that
for sufficiently large $x$, we have $v(x + c) \leq C_1$ and
\[
\max_{\gamma x \leq y \leq x + c} u(y) \leq C_2 x e^{-\sqrt{2} \gamma x}.
\]
It follows that for sufficiently large $x$,
\begin{eqnarray*}
v(x + c) - v(x) &\leq& C_1 E \biggl[ 1 - \exp\biggl( - \int_0^T u(R_s) \,
ds \biggr) \biggr] \\
&\leq& C_1 E \biggl[ \mathbf{1}_A + \biggl( \int_0^T u(R_s) \, ds \biggr)
\mathbf{1}_{A^c} \biggr] \\
&\leq& C_1 P(A) + C_1 C_2 x e^{-\sqrt{2} \gamma x} E[T] .
\end{eqnarray*}

To bound $E[T]$, note that using $E_x$ to denote expectation for the
Bessel process started at $x$, and $\tau_z$ to be the first time that
the Bessel process hits $z$, we have $E_0[\tau_{x+c}] = E_0[\tau_x] +
E_x[\tau_{x + c}]$ by the strong
Markov property.
Therefore, $E[T] = E_x[\tau_{x + c}] \leq E_0[\tau_{x + c}]$.
Furthermore, the three-dimensional Bessel process is the
Euclidean norm of three-dimensional Brownian motion, which is bounded
below by the absolute value of the first coordinate, which is a
one-dimensional Brownian motion. Therefore, $E_0[\tau_{x + c}]$ is at
most the the expected time for a one-dimensional Brownian motion to
reach $-(x + c)$ or $x + c$, which for sufficiently large $x$ is at
most $C_3 x^2$ for some constant $C_3$. It follows that
\begin{eqnarray*}
\limsup_{x \rightarrow\infty} x\bigl(v(x + c) - v(x)\bigr) &\leq&\limsup_{x
\rightarrow\infty} \biggl( x \cdot C_1 \frac{c \gamma}{c + (1 - \gamma)
x} + C_1 C_2 C_3 x^4 e^{-\sqrt{2} \gamma x} \biggr)\\
& =& \frac{c \gamma
C_1}{1 - \gamma}.
\end{eqnarray*}
Because this holds for any $\gamma> 0$, and $C_1$ does not depend on
$\gamma$, the result follows.
\end{pf}

%
%pr27 #&#
\begin{Prop}\label{Wprop}
Let $W$ be the limiting random variable in Proposition~\ref{nevprop}.
Then, there exists a constant $B > 0$ such that as $x \rightarrow
\infty
$,
\[
P(W > x) \sim\frac{B}{x}.
\]
\end{Prop}

\begin{pf}
Let $\phi(\lambda) = E[e^{-\lambda W}]$.
According to the discussion on page 335 of~\cite{regvar}, the condition
that $P(W > x) \sim B/x$ as $x \rightarrow\infty$ is equivalent to the
condition that the function $f(z) = z (1 - \phi(1/z))$ has $B$-index 1,
meaning (see page 128 of~\cite{regvar}) that for all $r \geq1$, we
have
\[
\lim_{z \rightarrow\infty} \bigl(f(rz) - f(z)\bigr) = B \log r.
\]
That is, $P(W > x) \sim B/x$ is equivalent to the condition that for
all $r \geq1$, we have
\[
\lim_{z \rightarrow\infty} rz \bigl(1 - \phi(1/rz)\bigr) - z\bigl(1 - \phi(1/z)\bigr)
= B
\log r,
\]
or equivalently, letting $\lambda= 1/z$,
%
%e90 ###
%
%e90 #&#
\begin{equation}\label{lamx}
\lim_{\lambda\rightarrow0} \frac{r(1 - \phi(\lambda/r)) - (1 -
\phi
(\lambda))}{\lambda} = B \log r.
\end{equation}
Consequently we need to show that (\ref{lamx}) holds for all $r \geq1$.

By (\ref{LapW}), we have $\phi(\lambda) = \psi((\log\lambda
)/\sqrt
{2})$. Let $w(x) = \psi(-x)$, so $w$ satisfies (\ref{wdiffeq}). For all
$x$, let $u(x) = 1 - w(x)$ and $v(x) = u(x)/(x e^{-\sqrt{2} x})$, as in
Proposition~\ref{KPPasym}. Then
\begin{eqnarray*}
1 - \phi(\lambda) &=& 1 - \psi\biggl( \frac{\log\lambda}{\sqrt{2}}
\biggr) = u \biggl( \frac{-\log\lambda}{\sqrt{2}} \biggr) = u \biggl( \frac{\log
(1/\lambda)}{\sqrt{2}} \biggr) \\
&=& v \biggl( \frac{\log(1/\lambda)}{\sqrt
{2}} \biggr)\biggl ( \frac{\log(1/\lambda)}{\sqrt{2}} \biggr) \lambda.
\end{eqnarray*}
Likewise,
\[
1 - \phi(\lambda/r) = v \biggl( \frac{\log(r/\lambda)}{\sqrt{2}} \biggr)
\biggl( \frac{\log(r/\lambda)}{\sqrt{2}} \biggr) \frac{\lambda}{r}.
\]
Letting $x = (\log(1/\lambda))/\sqrt{2}$ and $c = (\log r)/\sqrt{2}$,
it follows that
\begin{eqnarray*}
&&\frac{r(1 - \phi(\lambda/r)) - (1 - \phi(\lambda))}{\lambda}\\
 && \qquad = v
\biggl( \frac{\log(r/\lambda)}{\sqrt{2}} \biggr)\biggl ( \frac{\log(r/\lambda
)}{\sqrt{2}} \biggr) - v \biggl( \frac{\log(1/\lambda)}{\sqrt{2}} \biggr)
\biggl( \frac{\log(1/\lambda)}{\sqrt{2}} \biggr) \\
&& \qquad = v(x + c)(x + c) - v(x) x = x\bigl(v(x + c) - v(x)\bigr) + c v(x + c).
\end{eqnarray*}
As $x \rightarrow\infty$, we have $v(x + c) \rightarrow C$, where $C$
is the constant from (\ref{wasym}), and $x(v(x + c) - v(x))
\rightarrow
0$ by Proposition~\ref{KPPasym}. Therefore,
\begin{eqnarray*}
\lim_{\lambda\rightarrow0} \frac{r(1 - \phi(\lambda/r)) - (1 -
\phi
(\lambda))}{\lambda} &=& \lim_{x \rightarrow\infty} \bigl[x\bigl(v(x + c) - v(x)\bigr)
+ c v(x + c)\bigr] \\
&=& \frac{C \log r}{\sqrt{2}},
\end{eqnarray*}
so (\ref{lamx}) holds with $B = C/\sqrt{2}$.
\end{pf}

We will see later in the proof of Proposition~\ref{rjumplem} that $B =
1/\sqrt{2}$.

%
%co28 #&#
\begin{Cor}\label{WCor}
There is a constant $C$ such that $P(W > x) \leq C/x$ for all $x$, and
$E[W \mathbf{1}_{\{W \leq x\}}] \leq C \log x$ and $E[W^2 \mathbf
{1}_{\{
W \leq x\}}] \leq Cx$ for all $x \geq2$.
\end{Cor}

\begin{pf}
The first statement is immediate from Proposition~\ref{Wprop}. Since
\[
E\bigl[W \mathbf{1}_{\{W \leq x\}}\bigr] \leq\int_0^x P(W \geq y) \, dy \leq1 +
\int_1^x \frac{C}{y} \, dy = 1 + C \log x
\]
and
\[
E\bigl[W^2 \mathbf{1}_{\{W \leq x\}}\bigr] \leq\int_0^x 2y P(W \geq y) \, dy
\leq1 + 2 \int_1^x y \cdot\frac{C}{y} \, dy \leq1 + 2Cx,
\]
the other two statements follow.
\end{pf}

%s5 ###
%s5 #&#
\section{The particles after hitting $L_A$}\label{sec5}

Recall that in Section~\ref{sec31}, we obtained estimates on the
number of particles in branching Brownian motion that never reach the
level $L_A$, while in Section~\ref{sec32} we estimated the number of
particles that reach $L_A$. In this section, we determine how much the
descendants of the particles that reach $L_A$ will contribute to the
process at later times. The basic strategy will be to argue that if a
particle reaches $L_A$, then the number of descendants that it will
have in the population a long time into the future can be approximated
by the number of its descendants that reach $L_A - y$, where $y$ is
some large constant. The number of descendants that reach $L_A - y$ can
be approximated using the random variable $W$ in Proposition~\ref{nevprop}.

%s5.1 ###
%s5.1 #&#
\subsection{Notation and constants}\label{sec51}

Recall from  Section~\ref{sec32} that $R_k$ particles reach $L_A$ between
times $t_{k-1}$ and $t_k$. By Propositions~\ref{EMkProp} and \ref
{EMk2Prop}, on $G_{N, k-1}$ we have
%
%e91 ###
%
%e91 #&#
\begin{equation} \label{EMkV}
E[R_k|{\cal F}_{t_{k-1}}] \leq C \theta e^A \eps^{-1/2} + o(1)
\end{equation}
and
%
%e92 ###
%
%e92 #&#
\begin{equation} \label{EMk2V}
E[R_k^2|{\cal F}_{t_{k-1}}] \leq C \theta e^A \eps^{-1/2} + o(1).
\end{equation}
These moment estimates will be used repeatedly in what follows.
Denote by $u_1 < u_2 <\cdots< u_{R_k}$ the times at which these
particles reach $L_A$. Recalling (\ref{ZNdef}) and (\ref{ZN1def}), define
\[
Z_{N,2}(t_k) = \sum_{i=1}^{M_N(t_k)} e^{\mu X_{i,N}(t_k)} \sin\biggl(
\frac{\pi X_{i,N}(t_k)}{L} \biggr) \mathbf{1}_{\{i \notin S(t_k)\}}
\mathbf{1}_{\{X_{i,N}(t_k) \leq L\}}.
\]
Note that
\[
Z_N(t_k) = Z_{N,1}(t_k) + Z_{N,2}(t_k),
\]
and the particles contributing to $Z_{N,2}(t_k)$ are precisely the
particles at time $t_k$ that are descended from the particles that
reach $L_A$ at one of the times $u_1,\ldots, u_{R_k}$.

Our aim in this section will be to estimate, on $G_{N, k-1}$, the expectation
%
%e93 ###
%
%e93 #&#
\begin{equation}\label{truncexp}
E\bigl[\bigl(Z_N(t_k) - Z_N(t_{k-1})\bigr) \mathbf{1}_{\{Z_N(t_k) - Z_N(t_{k-1}) \leq
\eps N (\log N)^2\}}|{\cal F}_{t_{k-1}}\bigr],
\end{equation}
as well as probabilities of the form
%
%e94 ###
%
%e94 #&#
\begin{equation}\label{bigjeq}
P\bigl(Z_N(t_k) - Z_N(t_{k-1}) > r N (\log N)^2|{\cal F}_{t_{k-1}}\bigr)
\end{equation}
for $r \geq\eps$.
We apply the truncation at $\eps N (\log N)^2$ to focus separately on
particles reaching $L_A$ that make a small addition to the value of the
process, whose contributions are counted in (\ref{truncexp}), and
particles reaching $L_A$ that lead to large jumps in the value of the
process, an event whose probability is estimated by (\ref{bigjeq}).

Estimating these quantities precisely will involve manipulating seven constants.
Recall that we have been already working with the three constants $\eps
$, $A$ and $\theta$. Throughout this section, $\eps$ will be a fixed
number with $0 < \eps< 1$. We will also introduce a new constant
$\delta> 0$ and in fact will fix
%
%e95 ###
%
%e95 #&#
\begin{equation}
\delta\le\eps^7
\label{deltasmall}.
\end{equation}
By Proposition~\ref{Wprop}, one can choose $x$ large enough that if $z
\geq x$, then
%
%e96 ###
%
%e96 #&#
\begin{equation}\label{WBbound}
\frac{(1 - \delta) B}{z} \leq P(W > z) \leq\frac{(1 + \delta) B}{z},
\end{equation}
where $B$ comes from Proposition~\ref{Wprop}. We will then choose $A
\geq1$ large enough that
%
%e98 ###
%e97 ###
%
%e97 #&#
%e98 #&#
\begin{eqnarray} \label{651}
2 \sqrt{2} \pi e^{-A} x &\leq&\eps, \\ \label{652}
4 e^{-A/9} &\leq&\delta/6.
\end{eqnarray}
Once $\eps$, $\delta$ and $A$ are chosen, we will choose $\theta> 0$
small enough to satisfy the following equations:
%
%e105 ###
%e104 ###
%e103 ###
%e102 ###
%e101 ###
%e100 ###
%e99 ###
%
%e99 #&#
%e100 #&#
%e101 #&#
%e102 #&#
%e103 #&#
%e104 #&#
%e105 #&#
\begin{eqnarray}\label{8613}
A \theta&\leq&1; \\\label{8623}
4 \pi^2 A \theta s \eps^{-1/2} &\leq& e^{-A/4}; \\\label{654}
4 \theta^{1/4} &\leq&\delta/6; \\\label{655}
\theta^{1/4} e^A &\leq&\delta; \\\label{657}
\theta A^2 &\leq&\delta^{1/2}; \\\label{65star}
\theta A^2 e^A \eps^{-1/2}& \leq&1; \\\label{658}
C_0 A \theta^{1/2} &\leq&1,
\end{eqnarray}
where $C_0$ is a constant to be defined later in (\ref{EcondG}).
Note that (\ref{8613}) and (\ref{8623}) were already assumed in (\ref
{con1}) and (\ref{con2}), while (\ref{65star}) implies (\ref{con3})
because $A \geq1$.
In this section, we will also work with the additional constants $\eta
$, $y$ and $\zeta$ from Corollary~\ref{nevcor}. We will choose $\eta=
\theta$. We will then choose $y$ to be large enough to satisfy both
(\ref{neveta})
and the equation
%
%e106 ###
%
%e106 #&#
\begin{equation}\label{659}
1 \leq\theta y.
\end{equation}
We finally choose $\zeta$ to satisfy the conditions of Corollary \ref
{nevcor} for these values of $\eta$ and $y$.

Consider the particle that reaches $L_A$ at time $u_j \in(t_{k-1},
t_k]$. Denote by $V_{j,k}$ the number of descendants of this particle
that, at some time $t > u_j$, reach $L_A - y + (t - u_j)(\sqrt{2} -
\mu
)$ and have the property that, for all $u \in[u_j, t)$, the ancestor
of this particle was in the interval $(L_A - y + (u - u_j)(\sqrt{2} -
\mu), \infty)$.
This is equivalent to the number of descendant particles that would
reach $L_A - y + (t - u_j)(\sqrt{2} - \mu)$ at time $t$ for some $t$ if
particles were killed upon reaching this level. Denote the first times
at which these $V_{j,k}$ particles reach level $L_A - y + (t -
u_j)(\sqrt{2} - \mu)$ by $r_{1,j,k} < r_{2,j,k} <\cdots< r_{V_{j,k},
j, k}$. Note that $V_{j,k}$ has the same distribution as the random
variable $Z_y$ of Proposition~\ref{nevprop},
and the adjustment of $(t - u_j)(\sqrt{2} - \mu)$ is necessary because
particles drift to the left at rate $\mu$, rather than at rate $\sqrt
{2}$ as in the setting of Proposition~\ref{nevprop}. Now let
\[
W_{j,k}' = y e^{-\sqrt{2} y} V_{j,k}.
\]
By Corollary~\ref{nevcor}, there exists a random variable $W_{j,k}$
with the same distribution as the random variable $W$ in Corollary \ref
{nevcor} such that $P(|W_{j,k}' - W_{j,k}| > \eta) < \eta$.
Furthermore, it is clear that for fixed $k$, conditional on ${\cal
F}_{t_{k-1}}$ and conditional on $R_k = r$, the random variables
$W_{1,k}',\ldots, W_{r,k}'$ are independent and have\vspace*{-2pt} the same
distribution as $y e^{-\sqrt{2} y} Z_y$.
Likewise, the random variables $W_{j,k}$ can be chosen such that
conditional on ${\cal F}_{t_{k-1}}$ and conditional on $R_k = r$,
$W_{1,k},\ldots, W_{r,k}$ are independent and have the same
distribution as the random variable $W$ in Corollary~\ref{nevcor}.

%s5.2 ###
%s5.2 #&#
\subsection{The contribution of one particle at $L_A$}\label{sec52}

In this subsection, we show that the contribution to $Z_{N,2}(t_k)$
from the $j$th particle to hit $L_A$ can be approximated by $\pi\sqrt
{2} e^{-A} N (\log N)^2 W_{j,k}$. As a result, typically $Z_N(t_k) -
Z_N(t_{k-1}) > \eps N (\log N)^2$ precisely when $W_{j,k} > \eps/(\pi
\sqrt{2} e^{-A})$ for some $j \leq R_k$. Establishing the validity of
this approximation requires bounding the probabilities of several
unlikely events.

%
%le29 #&#
\begin{Lemma}\label{B1Lem}
Let $B_1$ be the event that there exist $j_1, j_2 \leq R_k$ with $j_1
\neq j_2$ such that\vspace*{1pt} $W_{j_1, k} \geq e^{2A/3}$ and $W_{j_2, k} \geq
e^{2A/3}$. Then on $G_{N, k-1}$, we have $P(B_1|{\cal F}_{t_{k-1}})
\leq C \theta e^{-A/3} \eps^{-1/2} + o(1).$
\end{Lemma}

\begin{pf}
Conditional on ${\cal F}_{t_{k-1}}$ and $R_k$, the expected number of
pairs $(j_1, j_2)$ with $j_1 \neq j_2$ such that $W_{j_1, k} \geq
e^{2A/3}$ and $W_{j_2, k} \geq e^{2A/3}$ is $ {R_k\choose2} P(W \geq
e^{2A/3})^2$, where $W$ is the random variable defined in Corollary
\ref
{nevcor}. By Proposition~\ref{Wprop}, $P(W \geq e^{2A/3}) \leq C
e^{-2A/3}$, so
\[
P(B_1|{\cal F}_{t_{k-1}}) \leq C E[R_k^2|{\cal F}_{t_{k-1}}] e^{-4A/3}
\leq C \theta e^{-A/3} \eps^{-1/2} + o(1),
\]
where the last inequality uses (\ref{EMk2V}).
\end{pf}

%
%le30 #&#
\begin{Lemma}\label{B2Lem}
Fix $r \geq\eps$, and let $B_2$ be the event that
\[
\frac{r - 4 e^{-A/4} - e^{-A/9} - 4\theta^{1/4}}{\pi\sqrt{2} e^{-A}}
\leq W_{j,k} \leq\frac{r + 4 e^{-A/4} + 4\theta^{1/4}}{\pi\sqrt{2} e^{-A}}
\]
for some $j \leq R_k$. On $G_{N, k-1}$, we have $P(B_2|{\cal
F}_{t_{k-1}}) \leq C \theta\delta\eps^{-5/2} + o(1)$,
where the constant $C$ does not depend on $r$.
\end{Lemma}

\begin{pf}
Let $\gamma= 4e^{-A/4} + e^{-A/9} + 4\theta^{1/4}$. Note that $\gamma
\leq\delta/2 \leq\eps/2$ because $4e^{-A/4} \leq4e^{-A/9} \leq
\delta
/6$ and $4 \theta^{1/4} \leq\delta/6$ by (\ref{652}) and (\ref{654}).
Assume $x$ is chosen so that (\ref{WBbound}) holds for $z \geq x$. By
(\ref{651}), we have $(r - \gamma)/(\pi\sqrt{2} e^{-A}) \geq(\eps-
\gamma)/(\pi\sqrt{2} e^{-A}) \geq x$. Therefore,
\begin{eqnarray*}
&&P\biggl ( \frac{r - \gamma}{\pi\sqrt{2} e^{-A}} \leq W \leq\frac{r +
\gamma}{\pi\sqrt{2} e^{-A}} \biggr) \\
&& \qquad \leq\frac{B (1 + \delta) \pi
\sqrt
{2} e^{-A}}{r - \gamma} - \frac{B (1 - \delta) \pi\sqrt{2}
e^{-A}}{r +
\gamma} \\
&& \qquad \leq C e^{-A} \biggl( \frac{1 + \delta}{r - \gamma} - \frac{1 - \delta
}{r + \gamma} \biggr) = C e^{-A} \biggl( \frac{2 \gamma+ 2 r \delta}{r^2
- \gamma^2} \biggr) \leq\frac{C e^{-A} \delta}{\eps^2}.
\end{eqnarray*}
It follows from this bound and Markov's inequality that $P(B_2|{\cal
F}_{t_{k-1}}) \leq C e^{-A} \delta\*\eps^{-2} E[R_k|{\cal
F}_{t_{k-1}}]$. The result now follows from (\ref{EMkV}).
\end{pf}

%
%le31 #&#
\begin{Lemma}\label{zigzag}
Let $B_3$ be the event that for some $j$, the particle that reaches
$L_A$ at time $u_j$ has a descendant that at some time $t \in(u_j, t_k]$
reaches $L_A - y + (t - u_j)(\sqrt{2} - \mu)$, and that this descendant
itself has a descendant that reaches $L_A$ before time $t_k$. Then on
$G_{N, k-1}$, we have $P(B_3|{\cal F}_{t_{k-1}}) \leq C e^A \theta
^{3/2} \eps^{-1/2} + o(1)$.
\end{Lemma}

\begin{pf}
The particle that reaches $L_A$ at time $u_j$ has $V_{j,k}$ descendants
that reach $L_A - y + (t - u_j)(\sqrt{2} - \mu)$ at some time $t >
u_j$. Let $A_{j,k}$ be the event that one of these particles reaches
$L_A - y + (t - u_j)(\sqrt{2} - \mu)$ at some time $t > u_j + \zeta$.
By Corollary~\ref{nevcor} and Proposition~\ref{Wprop}, since $\theta
=\eta<1$,
%
%e107 ###
%
%e107 #&#
\begin{eqnarray}\label{AWprob}
&&P(A_{j,k} \cup\{W_{j,k}' > \theta^{-1/2}\}  \mbox{ for some }j
\leq
R_k|{\cal F}_{t_{k-1}})\nonumber\\
 && \qquad \leq E[R_k|{\cal F}_{t_{k-1}}] \bigl(2 \eta+ P(W >
\theta^{-1/2} - \eta)\bigr) \\
&& \qquad \leq C E[R_k|{\cal F}_{t_{k-1}}] \bigl(\eta+ \sqrt{\theta}\bigr).
\nonumber
\end{eqnarray}
At most $y^{-1} e^{\sqrt{2} y} W_{j,k}'$ descendants of the particle
that reaches $L_A$ at time $u_j$ will hit $L_A - y + (t - u_j)(\sqrt{2}
- \mu)$ at some time $t \leq t_k$. This is an upper bound rather than
an equality because some particles may reach this level after time
$t_k$. We now consider $N$ large enough that $y \geq\zeta(\sqrt{2} -
\mu)$. On the event $A_{j,k}^c \cap\{W'_{j,k} \leq\theta^{-1/2}\}$,
the probability that a descendant of one of these particles reaches
$L_A$ by time $t_k$ can be bounded above by $y^{-1} e^{\sqrt{2} y}
\theta^{-1/2}$ times the probability that a single particle at $L_A - y
+ \zeta(\sqrt{2} - \mu)$ has a descendant that reaches $L_A$ by time
$(\log N)^3 \theta s$. Using Markov's inequality to bound this latter
probability by the expectation of the number of
such descendants, it follows from Proposition~\ref{EMkProp} that the
probability is bounded above by
\begin{eqnarray*}
&&\frac{C e^A}{N (\log N)^2} \biggl( \theta e^{\mu(L_A - y + \zeta(\sqrt
{2} - \mu))} \sin\biggl( \frac{\pi(L_A - y + \zeta(\sqrt{2} - \mu
))}{L_A} \biggr)\\
&&\hspace*{131pt}\hphantom{\frac{C e^A}{N (\log N)^2} \biggl(}{} + \frac{e^{\mu(L_A - y + \zeta(\sqrt{2} - \mu
))}}{\log N} \biggr)\bigl(1 + o(1)\bigr).
\end{eqnarray*}
Note that we are applying Proposition~\ref{EMkProp} in the case when $k
= 1$, and there is just a single particle initially at the location
$L_A - y + \zeta(\sqrt{2} - \mu)$.
Since $e^{\mu L_A} = N (\log N)^3 e^{-A} (1 + o(1))$, $\sin(\pi(L_A -
y + \zeta(\sqrt{2} - \mu))/L_A) \leq(Cy/\log N)(1 + o(1))$, and
$e^{\mu\zeta(\sqrt{2} - \mu)}$ is $1 + o(1)$ this expression can be
bounded above by
\begin{eqnarray*}
&&\frac{C e^A}{N (\log N)^2} \bigl ( \theta y e^{-\mu y} N (\log N)^2
e^{-A} + e^{-\mu y} N (\log N)^2 e^{-A}  \bigr)\bigl(1 + o(1)\bigr)\\
&& \qquad  \leq C e^{-\mu
y} (\theta y + 1)\bigl(1 + o(1)\bigr).
\end{eqnarray*}
Combining these observations gives
\begin{eqnarray*}
&&P(B_3|{\cal F}_{t_{k-1}})\\
 && \qquad \leq C E[R_k|{\cal F}_{t_{k-1}}] \bigl( \eta
+ \sqrt{\theta} + y^{-1} e^{\sqrt{2} y} \theta^{-1/2} \cdot e^{-\mu y}
(\theta y + 1) \bigr)\bigl(1 + o(1)\bigr) \\
&& \qquad \leq C E[R_k|{\cal F}_{t_{k-1}}] \bigl( \eta+ \sqrt{\theta} + \theta
^{-1/2} y^{-1} \bigr)\bigl(1 + o(1)\bigr).
\end{eqnarray*}
The result now follows from (\ref{EMkV}) and the assumptions that
$\eta
= \theta$ and $1 \leq\theta y$.
\end{pf}

Recall that the particles at time $t_k$ contributing to $Z_{N,2}(t_k)$
are precisely the particles at time $t_k$
that are descended from the particles that reach $L_A$ at one of the
times $u_1,\ldots, u_{R_k}$. To separate the contributions from each of
these particles, write $i \in S_j$ if the particle at $X_{i,N}(t_k)$ at
time $t_k$ is descended from the particle that was at $L_A$ at time
$u_j$. Then for $1 \leq j \leq R_k$, define
%
%e108 ###
%
%e108 #&#
\begin{equation}\label{ZN2j}\qquad
Z_{N,2,j}(t_k) = \sum_{i=1}^{M_N(t_k)} e^{\mu X_{i,N}(t_k)} \sin\biggl(
\frac{\pi X_{i,N}(t_k)}{L} \biggr) \mathbf{1}_{\{i \in S_j\}} \mathbf
{1}_{\{X_{i,N}(t) \leq L\}}.
\end{equation}
Note that $Z_{N,2}(t_k) = \sum_{j=1}^{R_k} Z_{N,2,j}(t_k)$. The next
lemma shows that $Z_{N,2,j}(t_k)$ is approximately determined by the
random variable $W_{j,k}$.

%
%le32 #&#
\begin{Lemma}\label{B4Lem}
Let $B_4$ be the event that for some $j \leq R_k$, we have
\[
\bigl| Z_{N,2,j}(t_k) - \pi\sqrt{2} e^{-A} N (\log N)^2 W_{j,k} \bigr| >
4N (\log N)^2 \theta^{1/4}.
\]
On $G_{N, k-1}$, we have $P(B_4|{\cal F}_{t_{k-1}}) \leq C e^A \theta
^{5/4} \eps^{-1/2} + o(1)$.
\end{Lemma}

\begin{pf}
Define a new random variable $Z_{N,2,j}'(t_k)$ by modifying
$Z_{N,2,j}(t_k)$ in the following three ways:
\begin{itemize}
\item We set $Z_{N,2,j}'(t_k)$ to zero if $u_j > t_k - (\log N)^{5/2}$.

\item We set $Z_{N,2,j}'(t_k)$ to zero if $r_{V_{j,k}, j, k} > u_j +
\zeta.$\vspace*{1pt}

\item We modify $S_j$ to exclude particles that, after time $u_j$,
reach $L_A - y + (t - u_j)(\sqrt{2} - \mu)$ at some time $t \in(u_j,
t_k]$ but then reach $L_A$ again before time~$t_k$.
[Note that this modification is equivalent to killing particles that\vspace*{1pt}
reach $L_A$ after they reach $L_A - y + (t - u_j)(\sqrt{2} - \mu)$ at
some time $t > u_j$.]
\end{itemize}
Then define $Z_{N,2,j}''(t_k)$ by making these three modifications and
replacing $L$ by $L_A$ in the definition (\ref{ZN2j}).

By Corollary~\ref{EMkCor} and Markov's inequality, $P(u_{R_k} > t_k -
(\log N)^{5/2}|{\cal F}_{t_{k-1}}) = o(1)$ on $G_{N, k-1}$. This
implies that the first of the four modifications above is unlikely to occur.
By Corollary~\ref{nevcor} and (\ref{EMkV}),
\[
P(r_{V_{j,k}, j, k} > u_j + \zeta\mbox{ for some }j|{\cal
F}_{t_{k-1}}) \leq\eta E[R_k|{\cal F}_{t_{k-1}}] \leq C \eta\theta
e^A \eps^{-1/2} + o(1),
\]
which bounds the probability of the second type of modification. Lemma
\ref{zigzag} bounds the probability of the third type of modification.
These results and the fact that $\eta= \theta$ imply that on $G_{N, k-1}$,
%
%e109 ###
%
%e109 #&#
\begin{eqnarray}\label{ZNjZNj}
&&P\bigl(Z_{N,2,j}'(t_k) \neq Z_{N,2,j}(t_k) \mbox{ for some }j \leq
R_k|{\cal
F}_{t_{k-1}}\bigr)\nonumber
\\[-8pt]
\\[-8pt] && \qquad \leq C e^A \theta^{3/2} \eps^{-1/2} + o(1).
\nonumber
\end{eqnarray}

Let $\Gamma_j$ be the event that $u_j \leq t_k - (\log N)^{5/2}$, that
$r_{V_{j,k}, j, k} \leq u_j + \zeta$ and that $W_{j,k}' \leq\theta^{-1/4}$.
The probability that either of the first two of these events fails to
occur has already been bounded, so using the argument given in (\ref
{AWprob}), on $G_{N, t_{k-1}}$ we have
%
%e110 ###
%
%e110 #&#
\begin{equation}\label{Gjprime}
P \Biggl( \bigcup_{j=1}^{R_k} \Gamma_j^c \Big| {\cal F}_{t_{k-1}}
\Biggr) \leq C e^A \theta^{5/4} \eps^{-1/2} + o(1).
\end{equation}
Let ${\cal H}_{k-1} = \sigma({\cal F}_{t_{k-1}}, V_{1,k},\ldots,
V_{R_k, k}, u_1,\ldots, u_{R_k}, (r_{i,j,k})_{1 \leq i \leq V_{j,k}, 1
\leq j \leq R_k})$.\vspace*{1pt} Note that $\Gamma_j \in{\cal H}_{k-1}$ for all
$j$, and on $\Gamma_j$ for sufficiently large $N$, the $V_{j,k}$
particles that reach $L_A - y + (t - u_j)(\sqrt{2} - \mu)$ for some $t
> u_j$ are all reaching a level between $L_A - y$ and $L_A - y + \zeta
(\sqrt{2} - \mu)$ at some time between $t_{k-1}$ and $t_k$. These
particles and their descendants then evolve independently until time
$t_k$, and we kill particles that return to $L_A$ if we are evaluating
$Z_{N,2,j}'(t_k)$ or $Z_{N,2,j}''(t_k)$.\vspace*{1pt}

By the argument leading to (\ref{new27}), with the times $r_{i,j,k}$
playing the role of $t_{k-1}$, on $\Gamma_j$ we have
\begin{eqnarray*}
&&E [ |Z_{N,2,j}'(t_k) - Z_{N,2,j}''(t_k)| | {\cal H}_{k-1} ]\\
&& \qquad \leq V_{j,k} \frac{C A e^{\mu(L_A - y)}}{\log N} \bigl(1 + o(1)\bigr)
\\
&& \qquad \leq C y^{-1} e^{\sqrt{2} y} W_{j,k}' \cdot\frac{A e^{-A} N (\log
N)^3 e^{-\mu y}}{\log N} \bigl(1 + o(1)\bigr) \\
&& \qquad \leq C y^{-1} N (\log N)^2 \theta^{-1/4} \bigl(1 + o(1)\bigr).
\end{eqnarray*}
Therefore, by Markov's inequality and assumption (\ref{659}), that $1
\leq\theta y$, on $\Gamma_j$, we have
%
%e111 ###
%
%e111 #&#
\begin{eqnarray}\label{ZN2jp}
 \qquad P \bigl( | Z_{N,2,j}'(t_k) - Z_{N,2,j}''(t_k)| > N (\log N)^2 \theta
^{1/4} | {\cal H}_{k-1} \bigr) &\leq& C y^{-1} \theta^{-1/2} +
o(1)\nonumber
\\[-4pt]
\\[-12pt]
 \qquad &\leq& C \theta^{1/2} + o(1).
\nonumber
\end{eqnarray}
Let
\[
y_{i,j,k} = e^{\mu(L_A - y +(r_{i,j,k} - u_j)(\sqrt{2} - \mu))}
\]
and
\[
z_{i,j,k} = y_{i,j,k} \sin\biggl(\frac{\pi(L_A - y +(r_{i,j,k} -
u_j)(\sqrt{2} - \mu))}{L_A} \biggr).
\]
The $i$th of the $V_{j,k}$ particles that reach $L_A - y + (t -
u_j)(\sqrt{2} - \mu)$ for some $t > u_j$ reaches this level at time
$r_{i,j,k}$. Therefore, by (\ref{Zexp}) and (\ref{LAexp}), on $\Gamma
_j$ the expected contribution to $Z_{N,2,j}''(t_k)$ from descendants of
this particle is given by
\begin{eqnarray*}
&&e^{(1- \mu^2/2 - \pi^2/2L_A^2)(t_k - r_{i,j,k})}z_{i,j,k}\\
&& \qquad = \bigl(1+ O(A\theta) + o(1)\bigr) e^{\mu(L_A- y)}\frac{\pi y}{L_A} \biggl(1-
O \biggl(\frac{\zeta(\sqrt{2}- \mu)}y \biggr) \biggr)\\
&& \qquad  = e^{\mu(L_A- y)}\frac{\pi y}{L_A}\bigl(1+ O(A\theta) + o(1)\bigr).
\end{eqnarray*}
Thus, on $\Gamma_j$,
\begin{eqnarray*}
E[Z_{N,2,j}''(t_k)|{\cal H}_{k-1}] &=& V_{j,k} \biggl( e^{\mu(L_A - y)}
\frac{\pi y}{L_A}\bigl(1 + O(A \theta) + o(1)\bigr) \biggr)\\
&=& V_{j,k}\biggl ( N (\log N)^3 e^{-A} e^{-\mu y} \frac{\pi y}{L_A}
\biggr) \bigl(1 + O(A \theta) + o(1)\bigr) \\
&=& W_{j,k}' \bigl(\pi\sqrt{2} e^{-A} N (\log N)^2\bigr)\bigl(1 + O(A \theta) + o(1)\bigr).
\end{eqnarray*}
This means there is a constant $C_0$ such that for sufficiently large
$N$, on $\Gamma_j$,
\begin{eqnarray*}
\bigl| E[Z_{N,2,j}''(t_k)|{\cal H}_{k-1}] - \pi\sqrt{2} e^{-A} N (\log
N)^2 W_{j,k}' \bigr| &\leq& C_0 N (\log N)^2 W_{j,k}' A \theta\\
&\leq& C_0 N
(\log N)^2 A \theta^{3/4}.
\end{eqnarray*}
Therefore, using (\ref{658}),
%
%e112 ###
%
%e112 #&#
\begin{equation}\label{EcondG}
 \quad \bigl| E[Z_{N,2,j}''(t_k)|{\cal H}_{k-1}] - \pi\sqrt{2} e^{-A} N (\log
N)^2 W_{j,k}' \bigr| \leq N (\log N)^2 \theta^{1/4}
\end{equation}
for sufficiently large $N$. On $\Gamma_j$ we can similarly estimate the
variance of the contribution of each of these particles. We apply (\ref
{varZNp}), with the times $r_{i,j,k}$ playing the role of $t_{k-1}$.
Since the descendants of these particles after times $r_{1,j,k}, \ldots
, r_{V_{j,k}, j,k}$ evolve independently, we get
\[
\operatorname{Var}(Z_{N,2,j}''(t_k)|{\cal H}_{k-1}) \leq\sum_{i=1}^{V_{j,k}}
C \theta N (\log N)^2 e^{-A} \biggl(z_{i,j,k} + \frac{y_{i,j,k}}{\theta
\log N} \biggr)\bigl (1 + o(1)\bigr).
\]
Arguing as above and using (\ref{659}), we find that on $\Gamma_j$,
%
%e113 #&#
\begin{eqnarray} \label{VcondG}
 &&\operatorname{Var}(Z_{N,2,j}''(t_k)|{\cal H}_{k-1})\nonumber\\
  && \qquad \leq
C V_{j,k} \theta N (\log N)^2\nonumber\\
&& \qquad  \quad {}\times e^{-A} \biggl( e^{\mu(L_A - y)} \sin
\biggl( \frac{\pi(L_A - y)}{L_A} \biggr) + \frac{e^{\mu(L_A - y)}}{\theta
\log N} \biggr)\bigl(1 + o(1)\bigr)  \nonumber
\\
  && \qquad \leq C V_{j,k} \theta N (\log N)^2  \\
  && \qquad  \quad {}\times e^{-2A} \bigl( y e^{-\mu y} N (\log
N)^2 + \theta^{-1} e^{-\mu y} N (\log N)^2 \bigr)\bigl(1 + o(1)\bigr) \nonumber
\\
  && \qquad \leq C W_{j,k}' N^2 (\log N)^4 \theta e^{-2A} ( 1 + \theta^{-1}
y^{-1} )\bigl(1 + o(1)\bigr) \nonumber\\
  && \qquad \leq C N^2 (\log N)^4 \theta^{3/4}\bigl(1 + o(1)\bigr).
\nonumber
\end{eqnarray}
By (\ref{EcondG}), (\ref{VcondG}) and the conditional form of
Chebyshev's inequality, on $\Gamma_j$ we have
\begin{eqnarray*}
&&P \bigl( \bigl|Z_{N,2,j}''(t_k) - \pi\sqrt{2} e^{-A} N (\log N)^2 W_{j,k}'
\bigr| > 2 N (\log N)^2 \theta^{1/4} | {\cal H}_{k-1} \bigr) \\
&& \qquad\leq\frac{C N^2 (\log N)^4 \theta^{3/4}(1 + o(1))}{(N
(\log N)^2
\theta^{1/4})^2} \leq C \theta^{1/4} + o(1).
\end{eqnarray*}
Note that $\pi\sqrt{2} e^{-A} \eta\leq\theta^{1/4}$ because $A
\geq
0$, $\eta= \theta$, $\delta\leq1$ by (\ref{deltasmall}), and thus
$\theta^{3/4} \leq1/24^3$ by (\ref{654}). Therefore, since
$P(|W_{j,k}' - W_{j,k}| > \eta) < \eta$, on $\Gamma_j$ we have
\begin{eqnarray*}
&&P \bigl( \bigl|Z_{N,2,j}''(t_k) - \pi\sqrt{2} e^{-A} N (\log N)^2 W_{j,k} \bigr|
> 3 N (\log N)^2 \theta^{1/4} | {\cal H}_{k-1} \bigr) \\
&& \qquad\leq C \theta^{1/4} + \eta+ o(1) \leq C \theta^{1/4} + o(1).
\end{eqnarray*}
Now (\ref{ZN2jp}) leads to
\begin{eqnarray*}
&&P \bigl(\bigl|Z_{N,2,j}'(t_k) - \pi\sqrt{2} e^{-A} N (\log N)^2 W_{j,k}\bigr| > 4
N (\log N)^2 \theta^{1/4} | {\cal H}_{k-1} \bigr) \\
&& \qquad \leq\mathbf
{1}_{\Gamma_j^c} + C \theta^{1/4} + o(1).
\end{eqnarray*}
Taking the union over over $j \leq R_k$ and then taking conditional
expectations of both sides with respect to ${\cal F}_{t_{k-1}}$, we get
%
%e113 ###
%
%e114 #&#
\begin{eqnarray}
&&P \bigl( \bigl|Z_{N,2,j}'(t_k) - \pi\sqrt{2} e^{-A} N (\log N)^2 W_{j,k}\bigr| >
4 N (\log N)^2 \theta^{1/4}\nonumber\\
&& \hspace*{138pt}\qquad  \hphantom{P \bigl(}\mbox{for some }j \leq R_k | {\cal
F}_{t_{k-1}} \bigr) \\
&& \qquad\leq P \Biggl( \bigcup_{j=1}^{R_k} \Gamma_j^c \Big| {\cal F}_{t_{k-1}}
\Biggr) + \bigl(C \theta^{1/4} + o(1)\bigr) E[R_k|{\cal F}_{t_{k-1}}].
\nonumber
\end{eqnarray}
The result now follows from (\ref{ZNjZNj}), (\ref{Gjprime}) and (\ref{EMkV}).
\end{pf}

%
%le33 #&#
\begin{Lemma}\label{B5Lem}
Let $B_5$ be the event that
\[
\sum_{j=1}^{R_k} Z_{N,2,j}(t_k) \mathbf{1}_{\{W_{j,k} \leq e^{2A/3}\}}
> e^{-A/9} N (\log N)^2
\]
or
\[
\sum_{j=1}^{R_k} W_{j,k} \mathbf{1}_{\{W_{j,k} \leq e^{2A/3}\}} >
\frac
{e^{8A/9}}{\pi\sqrt{2}}.
\]
Then $P(B_5|{\cal F}_{t_{k-1}}) \leq C (\theta^{5/4} e^A + \theta
e^{-A/9}) \eps^{-1/2} + o(1)$ on $G_{N, k-1}$.
\end{Lemma}

\begin{pf}
We have
%
%e114 ###
%
%e115 #&#
\begin{eqnarray} \label{B5eq1}\qquad
P(B_5|{\cal F}_{t_{k-1}}) &\leq& P(B_4|{\cal F}_{t_{k-1}}) \nonumber
\\[-8pt]
\\[-8pt]&&{}+ P\Biggl ( \sum
_{j=1}^{R_k} (W_{j,k} + \beta) \mathbf{1}_{\{W_{j,k} \leq e^{2A/3}\}} >
\frac{e^{8A/9}}{\pi\sqrt{2}} \Big| {\cal F}_{t_{k-1}} \Biggr),
\nonumber
\end{eqnarray}
where $\beta= 4 e^A \theta^{1/4}/(\pi\sqrt{2})$, which by (\ref{655})
is bounded by a constant.
Let ${\cal G}_{k-1} = \sigma({\cal F}_{t_{k-1}}, R_k)$. Using Corollary
\ref{WCor} with $x = e^{2A/3}$ and (\ref{EMkV}),
%
%e115 ###
%
%e116 #&#
\begin{eqnarray} \label{condvar1}
&&E \Biggl[ \operatorname{Var}\Biggl ( \sum_{j=1}^{R_k} (W_{j,k} + \beta) \mathbf
{1}_{\{W_{j,k} \leq e^{2A/3}\}} | {\cal G}_{k-1} \Biggr) \Big|
{\cal F}_{t_{k-1}} \Biggr]\nonumber\\
&& \qquad = E\bigl[R_k \operatorname{Var}\bigl((W + \beta) \mathbf{1}_{\{W \leq
e^{2A/3}\}}\bigr)|{\cal
F}_{t_{k-1}}\bigr] \nonumber\\
&& \qquad \leq E\bigl[(W + \beta)^2 \mathbf{1}_{\{W \leq e^{2A/3}\}}\bigr] E[R_k|{\cal
F}_{t_{k-1}}] \\
&& \qquad \leq C e^{2A/3} \cdot C \theta e^A \eps^{-1/2} + o(1)\nonumber \\
&& \qquad \leq C \theta e^{5A/3} \eps^{-1/2} + o(1).\nonumber
\end{eqnarray}
Likewise, using Corollary~\ref{WCor} and (\ref{EMk2V}),
%
%e116 ###
%
%e117 #&#
\begin{eqnarray}\label{condvar2}
&&\operatorname{Var} \Biggl( E \Biggl[ \sum_{j=1}^{R_k} (W_{j,k} + \beta)
\mathbf
{1}_{\{W_{j,k} \leq e^{2A/3}\}} \Big| {\cal G}_{k-1} \Biggr] \Big|
{\cal F}_{t_{k-1}} \Biggr) \nonumber\\
&& \qquad\leq\operatorname{Var} \bigl( R_k E\bigl[(W + \beta)\mathbf{1}_{\{W
\leq e^{2A/3}\}
}\bigr] | {\cal F}_{t_{k-1}} \bigr) \\
&& \qquad\leq\operatorname{Var}(CAR_k|{\cal F}_{t_{k-1}})
\leq C A^2 E[R_k^2|{\cal F}_{t_{k-1}}] \leq C \theta A^2 e^A \eps
^{-1/2} + o(1).\nonumber
\end{eqnarray}
Recall that for all random variables $X$ and $\sigma$-fields ${\cal F}$
and ${\cal G}$ with ${\cal F} \subset{\cal G}$,
\[
\operatorname{Var}(X|{\cal F}) = E[\operatorname{Var}(X|{\cal
G})|{\cal F}] + \operatorname
{Var}(E[X|{\cal G}]|{\cal F}).
\]
Therefore, summing (\ref{condvar1}) and (\ref{condvar2}) gives
\[
\operatorname{Var} \Biggl( \sum_{j=1}^{R_k} (W_{j,k} + \beta) \mathbf
{1}_{\{
W_{j,k} \leq e^{2A/3}\}} \Big| {\cal F}_{t_{k-1}} \Biggr) \leq C \theta
e^{5A/3} \eps^{-1/2} + o(1),
\]
as $A^2 e^{-2A/3}$ is bounded by a constant.
Also, using again Corollary~\ref{WCor} and since $A^2 \theta e^A \eps
^{-1/2} \le1$ by (\ref{65star}),
\begin{eqnarray*}
E \Biggl[ \sum_{j=1}^{R_k} (W_{j,k} + \beta) \mathbf{1}_{\{W_{j,k} \leq
e^{2A/3}\}} \Big| {\cal F}_{t_{k-1}} \Biggr] &\leq& C E(R_k|{\cal
F}_{t_{k-1}}) ( \beta+ 2A/3) \\
&\leq& C \theta\eps^{-1/2} Ae^A + o(1) \leq C + o(1).
\end{eqnarray*}
Thus by the conditional form of Chebyshev's inequality, we get
\begin{eqnarray*}
P \Biggl( \sum_{j=1}^{R_k} (W_{j,k} + \beta) \mathbf{1}_{\{W_{j,k} \leq
e^{2A/3}\}} > \frac{e^{8A/9}}{\pi\sqrt{2}} \Big| {\cal F}_{t_{k-1}}
\Biggr) &\leq&\frac{C \theta e^{5A/3} \eps^{-1/2}}{(e^{8A/9}/\pi\sqrt
{2})^2} + o(1)\\
 &\leq& C \theta e^{-A/9} \eps^{-1/2} + o(1),
\end{eqnarray*}
which, combined with (\ref{B5eq1}) and Lemma~\ref{B4Lem}, gives the result.
\end{pf}

%
%le34 #&#
\begin{Lemma}\label{B6Lem}
Fix $r \geq\eps$.
Consider the event $E$ that $Z_N(t_k) - Z_N(t_{k-1}) > r N (\log N)^2$,
and consider the event $F$ that $W_{j,k} > r/(\pi\sqrt{2} e^{-A})$ for
some $j \leq R_k$. Let $B_6$ be the event that one of these two events
occurs but not the other (i.e., the symmetric difference of these
two events).
Then $P(B_6|{\cal F}_{t_{k-1}}) \leq C \theta\delta\eps^{-5/2} +
o(1)$ on $G_{N,k-1}$, where the constant $C$ does not depend on $r$.
\end{Lemma}

\begin{pf}
Let $B_0$ be the event that $|Z_{N,1}(t_k) - Z_N(t_{k-1})| > 4 e^{-A/4}
N (\log N)^2$. By Corollary~\ref{driftcor} and Lemmas~\ref{B1Lem}--\ref{B5Lem} as well as the
assumptions (\ref{652}) and (\ref{655}), we have on $G_{N, k-1}$,
\[
P \Biggl( \bigcup_{i=0}^5 B_i \Big| {\cal F}_{t_{k-1}} \Biggr) \leq C
\delta\theta\eps^{-5/2} + o(1).
\]
Therefore, it suffices to show that
\[
B_6 \subset B=\bigcup_{i=0}^5 B_i.
\]

Thus, suppose first $\omega\in E^c \cap F$, and let us show that
$\omega\in B$. We have $W_{j,k} > r/(\pi\sqrt{2} e^{-A})$ for some $j
\leq R_k$. It follows that if $ \omega\in B_2^c$, we have $W_{j,k} >
(r + 4 e^{-A/4} + 4 \theta^{1/4})/(\pi\sqrt{2} e^{-A})$. If
furthermore $ \omega\in B_2^c \cap B_4^c$, we have $Z_{N,2,j}(t_k) > N
(\log N)^2 (r + 4e^{-A/4})$. Now if also $\omega\in B_0^c$, we have
$Z_{N,1}(t_k) \geq Z_N(t_{k-1}) - 4e^{-A/4} N (\log N)^2$, so on $B_2^c
\cap B_4^c \cap B_0^c$, we have $Z_N(t_k) \geq Z_{N,1}(t_k) +
Z_{N,2,j}(t_k) > Z_N(t_{k-1}) + r N (\log N)^2$, and so $E$
occurs. Since we have assumed that $\omega\notin E$, it must be that
$\omega\in B_0 \cup B_2 \cup B_4 \subset B$.

Alternatively, suppose $\omega\in E \cap F^c$,
hence $W_{j,k} \leq r/(\pi\sqrt{2} e^{-A})$ for all $j \leq R_k$. It
follows that on $B_2^c$, we have $W_{j,k} \leq(r - 4 e^{-A/4} -
e^{-A/9} - 4 \theta^{1/4})/\break(\pi\sqrt{2} e^{-A})$ for all $j \leq R_k$.
Then on $B_2^c \cap B_4^c$, we have $Z_{N,2,j}(t_k) \leq N (\log N)^2
(r - 4 e^{-A/4} - e^{-A/9})$ for all $j \leq R_k$. On $B_1^c$, there
exists at most one $j \leq R_k$ such that $W_{j,k} \geq e^{2A/3}$.
Therefore, on $B_2^c \cap B_4^c \cap B_1^c \cap B_5^c$, we have
\[
Z_{N,2}(t_k) = \sum_{j=1}^{R_k} Z_{N,2,j}(t_k) \leq N (\log N)^2 (r - 4
e^{-A/4}).
\]
Finally, on $B_0^c$, we have $Z_{N,1}(t_k) \leq Z_N(t_{k-1}) + 4
e^{-A/4} N (\log N)^2$, so on $\bigcap_{i=0}^5 B_i^c$, we have $Z_N(t_k)
\leq Z_N(t_{k-1}) + r N (\log N)^2$ which means that $E$ does not
occur. Since we assumed $\omega\in E$, it must be that
$\omega\in\bigcup_{i=0}^5 B_i = B$, which finishes the proof of the lemma.
\end{pf}

%s5.3 ###
%s5.3 #&#
\subsection{The small jumps}\label{sec53}

In this subsection, we estimate the expectation in (\ref{truncexp}),
which covers the case in which the process $Z_N$ does not jump by more
than $\eps N (\log N)^2$ between times $t_{k-1}$ and $t_k$. We have
\[
Z_N(t_k) - Z_N(t_{k-1}) = \bigl(Z_{N,1}(t_k) - Z_N(t_{k-1})\bigr) + \sum
_{j=1}^{R_k} Z_{N,2,j}(t_k).
\]
Lemma~\ref{B6Lem} with $r = \eps$ shows that with high probability, we
have $Z_N(t_k) - Z_N(t_{k-1}) > \eps N (\log N)^2$ if and only if one
of the random variables $W_{1,k},\ldots,\break W_{R_k, k}$ is greater than\vadjust{\goodbreak}
$\eps/(\pi\sqrt{2} e^{-A})$. Therefore, in view of Lemma~\ref{B4Lem},
we can approximate the quantity in (\ref{truncexp}) by
%
%e117 ###
%
%e118 #&#
\begin{eqnarray}\label{Skdef}
S_k &=& \bigl(Z_{N,1}(t_k) - Z_N(t_{k-1})\bigr) \nonumber
\\[-8pt]
\\[-8pt]&&{}+ \pi\sqrt{2} e^{-A} N (\log N)^2
\sum_{j=1}^{R_k} W_{j,k} \mathbf{1}_{\{W_{j,k} \leq\eps/(\pi\sqrt{2}
e^{-A})\}},
\nonumber
\end{eqnarray}
which omits the contributions from terms with $W_{j,k} > \eps/(\pi
\sqrt
{2} e^{-A})$.
We now calculate the expected value of $S_k$ and will later justify in
Lemma~\ref{Rerr} that this is sufficiently close to the quantity in
(\ref{truncexp}).

%
%le35 #&#
\begin{Lemma}\label{Rexp}
On $G_{N, k-1}$, we have
\begin{eqnarray*}
E[S_k|{\cal F}_{t_{k-1}}] &=& Z_N(t_{k-1}) \theta s\bigl (2 \sqrt{2} \pi^2
E\bigl[W \mathbf{1}_{\{W \leq\eps/(\pi\sqrt{2} e^{-A})\}}\bigr] - 2 \pi^2 A\bigr)
\\
&&{}+ O(A^2 \theta^2 \eps^{-1/2} N (\log N)^2) + o(N (\log N)^2).
\end{eqnarray*}
\end{Lemma}

\begin{pf}
By Lemma~\ref{ZN1exp}, we have on $G_{N, k-1}$
%
%e118 ###
%
%e119 #&#
\begin{eqnarray}\label{driftexp}
&&E[Z_{N,1}(t_k) - Z_N(t_{k-1})|{\cal F}_{t_{k-1}}] \nonumber
\\[-8pt]
\\[-8pt]&& \qquad = -Z_N(t_{k-1})\bigl(2 \pi
^2 A \theta s + O(A^2 \theta^2)\bigr) + o(N (\log N)^2).
\nonumber
\end{eqnarray}
Also, since the random variables $W_{j,k}$ are independent of one
another, and of ${\cal F}_{t_{k-1}}$ and $R_k$, we have
\[
E \Biggl[ \sum_{j=1}^{R_k} W_{j,k} \mathbf{1}_{\{W_{j,k} \leq\eps/(\pi
\sqrt{2} e^{-A})\}} \Big| {\cal F}_{t_{k-1}} \Biggr] = E\bigl[W \mathbf
{1}_{\{W \leq\eps/(\pi\sqrt{2} e^{-A})\}}\bigr] E[R_k|{\cal F}_{t_{k-1}}].
\]
Combining this result with Proposition~\ref{EMkProp}, we get on $G_{N, k-1}$,
%
%e119 ###
%
%e120 #&#
\begin{eqnarray}\label{driftexp2}
&&E \Biggl[ \pi\sqrt{2} e^{-A} N (\log N)^2 \sum_{i=1}^{R_k} W_{j,k}
\mathbf{1}_{\{W_{j,k} \leq\eps/(\pi\sqrt{2} e^{-A})\}} \Big| {\cal
F}_{t_{k-1}} \Biggr] \nonumber
\\
&& \qquad= E\bigl[W \mathbf{1}_{\{W \leq\eps/(\pi\sqrt{2}
e^{-A})\}}\bigr]\\
&& \qquad  \quad {}\times
\bigl(2 \sqrt
{2} \pi^2 \theta s Z_N(t_{k-1})\bigl(1 + O(A \theta)\bigr) + o (N (\log N)^2 ) \bigr).
\nonumber
\end{eqnarray}
Note that from Corollary~\ref{WCor},
%
%e120 ###
%
%e121 #&#
\begin{eqnarray}\label{WCA}
   E\bigl[W \mathbf{1}_{\{W \leq\eps/(\pi\sqrt{2} e^{-A})\}}\bigr] &\leq&1 + C
\log
\biggl( \frac{\eps}{\pi\sqrt{2} e^{-A}} \biggr)\nonumber
\\[-8pt]
\\[-8pt] &\leq&1 + C (\log\eps+
A) \leq CA,
\nonumber
\end{eqnarray}
since $\log\eps< 0$.
The result now follows by combining (\ref{driftexp}) and (\ref
{driftexp2}), and using (\ref{WCA}) to help bound some of the error terms.
\end{pf}

It remains to bound the expected error that is made when approximating
the increment $(Z_N(t_k) - Z_N(t_{k-1})) \mathbf{1}_{\{Z_N(t_k) -
Z_N(t_{k-1}) \leq\eps N (\log N)^2\}}$ by $S_k$.

%
%le36 #&#
\begin{Lemma}\label{Zvar}
We have $E[(Z_{N,1}'(t_k) - Z_N(t_{k-1}))^2|{\cal F}_{t_{k-1}}] \leq C
\theta N^2 (\log N)^4\* (e^{-A} \eps^{-1/2} + o(1))$ on $G_{N, k-1}$.
\end{Lemma}

\begin{pf}
By Lemmas~\ref{ZN1exp} and~\ref{VarZ1}, on $G_{N, k-1}$,
\begin{eqnarray*}
&&E\bigl[\bigl(Z_{N,1}'(t_k) - Z_N(t_{k-1})\bigr)^2|{\cal F}_{t_{k-1}}\bigr]\\
 && \qquad =
\operatorname
{Var}(Z_{N,1}'(t_k)|{\cal F}_{t_{k-1}}) + \bigl(E[Z_{N,1}'(t_k) -
Z_N(t_{k-1})|{\cal F}_{t_{k-1}}]\bigr)^2 \\
&& \qquad \leq C \theta N (\log N)^2 e^{-A} \bigl( Z_N(t_{k-1}) + o(N (\log N)^2)
\bigr) \\
&& \qquad  \quad {}+ \bigl(C A \theta Z_N(t_{k-1}) + o(N (\log N)^2) \bigr)^2
\\
&& \qquad \leq C \theta N^2 (\log N)^4 e^{-A} \eps^{-1/2} + C A^2 \theta^2 N^2
(\log N)^4 \eps^{-1} + o(N^2 (\log N)^4) \\
&& \qquad \leq C \theta N^2 (\log N)^4 \bigl(e^{-A} \eps^{-1/2} + A^2 \theta\eps
^{-1} + o(1)\bigr),
\end{eqnarray*}
and the result follows from (\ref{65star}).
\end{pf}

%
%le37 #&#
\begin{Lemma}\label{Wvar}
On $G_{N, k-1}$, we have
\[
E\Biggl [ \Biggl( \sum_{j=1}^{R_k} W_{j,k} \mathbf{1}_{\{W_{j,k} \leq\eps
/(\pi\sqrt{2} e^{-A})\}} \Biggr)^2 \Big| {\cal F}_{t_{k-1}} \Biggr]
\leq C \theta e^{2A} \eps^{1/2} + o(1).
\]
\end{Lemma}

\begin{pf}
Note that $e^{-A/9} \leq C \delta\leq C \eps$ by (\ref{deltasmall})
and (\ref{652}). Because $A^2 e^{-8A/9}$ is bounded above by a
constant, it follows that $A^2 \leq C \epsilon e^A$. Therefore, by
(\ref
{EMkV}), (\ref{EMk2V}), (\ref{WCA}) and Corollary~\ref{WCor}, on
$G_{N, k-1}$,
\begin{eqnarray*}
&&E \Biggl[\Biggl ( \sum_{j=1}^{R_k} W_{j,k} \mathbf{1}_{\{W_{j,k} \leq
\eps/(\pi\sqrt{2} e^{-A})\}} \Biggr)^2 \Big| {\cal F}_{t_{k-1}} \Biggr]
\\
&& \qquad= E[R_k|{\cal F}_{t_{k-1}}] E \bigl[ W^2 \mathbf{1}_{\{W \leq
\eps/(\pi
\sqrt{2} e^{-A})\}} \bigr] \\
&& \qquad  \quad {}+ E[R_k(R_k - 1)|{\cal F}_{t_{k-1}}] \bigl(E
\bigl[ W \mathbf{1}_{\{W \leq\eps/(\pi\sqrt{2} e^{-A})\}} \bigr] \bigr)^2
\\
&& \qquad\leq\bigl(C \theta e^A \eps^{-1/2} + o(1)\bigr) (\eps e^A + A^2) \leq
C \theta
e^{2A} \eps^{1/2} + o(1)
\end{eqnarray*}
as claimed.
\end{pf}

%
%le38 #&#
\begin{Lemma}\label{Rerr}
On $G_{N, k-1}$, we have
%
%e121 ###
%
%e122 #&#
\begin{eqnarray}\qquad
&&E \bigl[\bigl|S_k -\bigl (Z_N(t_k) - Z_N(t_{k-1})\bigr) \mathbf{1}_{\{Z_N(t_k) -
Z_N(t_{k-1}) \leq\eps N (\log N)^2\}}\bigr| | {\cal F}_{t_{k-1}} \bigr]
\nonumber
\\[-8pt]
\\[-8pt]\qquad
&& \qquad \leq C \theta N (\log N)^2 \delta^{1/2} \eps^{-3} + o(N
(\log N)^2).
\nonumber
\end{eqnarray}
\end{Lemma}

\begin{pf}
Throughout this proof, we work on the event $G_{N, k-1}$.
Choose $r = \eps$, and recall from the proof of Lemma~\ref{B6Lem} that
the event $B = \bigcup_{i=0}^5 B_i$ can also be written\vadjust{\goodbreak} as $B= \bigcup_{i=0}^6 B_i$ since $B_6 \subset\bigcup_{i=0}^5 B_i$.
We will bound the following three terms:
%
%e122 ###
%
%e123 #&#
\begin{eqnarray} \label{absterm1}
& \qquad\ \displaystyle E \bigl[ \mathbf{1}_{B^c}\bigl|S_k - \bigl(Z_N(t_k) - Z_N(t_{k-1})\bigr) \mathbf
{1}_{\{Z_N(t_k) - Z_N(t_{k-1}) \leq\eps N (\log N)^2\}} \bigr| |
{\cal F}_{t_{k-1}} \bigr];&
\\ \label{absterm2}
& \qquad\ \displaystyle E \bigl[ \mathbf{1}_B \bigl|\bigl(Z_N(t_k) - Z_N(t_{k-1})\bigr) \mathbf{1}_{\{
Z_N(t_k) - Z_N(t_{k-1}) \leq\eps N (\log N)^2\}}\bigr| | {\cal
F}_{t_{k-1}} \bigr];&
\\ \label{absterm3}
& \qquad\ \displaystyle E [\mathbf{1}_B |S_k| | {\cal F}_{t_{k-1}} ].&
\end{eqnarray}

We first bound (\ref{absterm1}). On $B_6^c$, we have $Z_N(t_k) -
Z_N(t_{k-1}) > \eps N (\log N)^2$ if and only if $W_{j,k} > \eps/(\pi
\sqrt{2} e^{-A})$ for some $j \leq R_k$. In this case, on the event
that $W_{j_0,k} > \eps/(\pi\sqrt{2} e^{-A})$ for some $j_0 \leq R_k$,
the difference between $S_k$ and $(Z_N(t_k) - Z_N(t_{k-1})) \mathbf
{1}_{\{
Z_N(t_k) - Z_N(t_{k-1}) \leq\eps N (\log N)^2\}}$ will\vspace*{1pt} simply be
$S_k$, as the latter expression will be zero. However, on $B_0^c$, we have
\[
|Z_{N,1}(t_k) - Z_N(t_{k-1})| \leq4 e^{-A/4} N (\log N)^2.
\]
By (\ref{652}) and the fact that $\delta\in(0, \eps)$ we have $\eps
/(\pi\sqrt{2} e^{-A}) \geq e^{2A/3}$. Thus $W_{j_0,k} \ge e^{2A/3}$
and on $B_1^c,$ for all $j\neq j_0, W_{j,k} \le e^{2A/3}$.
Thus, the definition of $B_5$ from Lemma~\ref{B5Lem} implies that on
$B_1^c \cap B_5^c$, we have
\[
\sum_{i=1}^{R_k} W_{j,k} \mathbf{1}_{\{W_{j,k} \leq\eps/(\pi\sqrt{2}
e^{-A})\}} \leq\frac{e^{8A/9}}{\pi\sqrt{2}}.
\]
Therefore, $|S_k| \leq(4 e^{-A/4} + e^{-A/9} ) N (\log N)^2$ on $B^c
\cap\{W_{j_0,k} > \eps/\break(\pi\sqrt{2} e^{-A})  \mbox{for some }j_0
\leq
R_k\}$.
If, however, $W_{j,k} \leq\eps/(\pi\sqrt{2} e^{-A})$ for all\break $j
\leq
R_k$, then on $B_4^c$, the difference between $S_k$ and $(Z_N(t_k) -\break
Z_N(t_{k-1})) \mathbf{1}_{\{Z_N(t_k) - Z_N(t_{k-1}) \leq\eps N (\log
N)^2\}}$ is bounded by $4 R_k N (\log N)^2 \theta^{1/4}$.\break Therefore,
\begin{eqnarray*}
&&\hspace*{-6.5pt}E \bigl[\bigl|S_k - \bigl(Z_N(t_k) - Z_N(t_{k-1})\bigr) \mathbf{1}_{\{Z_N(t_k) -
Z_N(t_{k-1}) \leq\eps N (\log N)^2\}}\bigr| \mathbf{1}_{B^c} | {\cal
F}_{t_{k-1}} \bigr] \\
&& \qquad\leq\bigl((4 e^{-A/4} + e^{-A/9}) \\
&&  \qquad \hphantom{\leq\bigl( }{}\times P\bigl(W_{j_0,k} > \eps/\bigl(\pi\sqrt{2}
e^{-A}\bigr) \mbox{ for some }j_0|{\cal F}_{t_{k-1}}\bigr) \\
&& \hspace*{114pt} \qquad \hphantom{\leq\bigl( }{}+ 4 \theta^{1/4}
E[R_k|{\cal F}_{t_{k-1}}] \bigr) N (\log N)^2.
\end{eqnarray*}
Now (\ref{EMkV}) gives $E[R_k|{\cal F}_{t_{k-1}}] \leq C \theta e^A
\eps
^{-1/2} + o(1)$, and Proposition~\ref{Wprop} implies
\begin{eqnarray*}
P\bigl(W_{j_0,k} \geq\eps/\bigl(\pi\sqrt{2} e^{-A}\bigr) \mbox{ for some
}j_0|{\cal
F}_{t_{k-1}}\bigr) &\leq& C e^{-A} \eps^{-1} E[R_k|{\cal F}_{t_{k-1}}]\\
 &\leq& C
\theta\eps^{-3/2} + o(1).
\end{eqnarray*}
Therefore,
%
%e125 ###
%
%e126 #&#
\begin{eqnarray}\label{abs1bound}\qquad
&&E\bigl [\bigl|S_k -\bigl (Z_N(t_k) - Z_N(t_{k-1})\bigr) \mathbf{1}_{\{Z_N(t_k) -
Z_N(t_{k-1}) \leq\eps N (\log N)^2\}}\bigr| \mathbf{1}_{B^c} | {\cal
F}_{t_{k-1}} \bigr]\nonumber
\\[-8pt]
\\[-8pt]
&& \qquad\leq\bigl(C \theta\eps^{-3/2} e^{-A/9} + C \theta^{5/4} e^A
\eps^{-1/2} +
o(1)\bigr) N (\log N)^2,
\nonumber
\end{eqnarray}
which gives a bound on (\ref{absterm1}).\vadjust{\goodbreak}

We next bound (\ref{absterm2}). By Lemma~\ref{B6Lem} and its proof,
we have
%
%e126 ###
%
%e127 #&#
\begin{equation}\label{PBbound}
P(B|{\cal F}_{t_{k-1}}) \leq C \theta\delta\eps^{-5/2} + o(1).
\end{equation}
The random variable in (\ref{absterm2}) is bounded in absolute value by
$\max\{Z_N(t_{k-1}), \break \eps N (\log N)^2\}$. Therefore, on $G_{N, k-1}$,
%
%e127 ###
%
%e128 #&#
\begin{eqnarray}\label{abs2bound}
&&  E \bigl[\bigl|\bigl(Z_N(t_k) - Z_N(t_{k-1})\bigr) \mathbf{1}_{\{Z_N(t_k) - Z_N(t_{k-1})
\leq\eps N (\log N)^2\}}\bigr| \mathbf{1}_B | {\cal F}_{t_{k-1}} \bigr]
\nonumber\hspace*{-35pt}
\\[-8pt]
\\[-8pt]
&&  \qquad\leq P(B|{\cal F}_{t_{k-1}}) N (\log N)^2 \eps^{-1/2} \leq C
\theta
\delta\eps^{-3} N (\log N)^2 + o(N(\log N)^2).
\nonumber\hspace*{-35pt}
\end{eqnarray}

It remains to bound (\ref{absterm3}). By the conditional
Cauchy--Schwarz inequality, Lemma~\ref{Zvar}, and (\ref{PBbound}),
%
%e128 ###
%
%e129 #&#
\begin{eqnarray}\label{CSeq}
&&E [|Z_{N,1}'(t_k) - Z_N(t_{k-1})| \mathbf{1}_B | {\cal F}_{t_{k-1}}
] \nonumber\\
&& \qquad \leq\sqrt{E\bigl[\bigl(Z_{N,1}'(t_k) - Z_N(t_{k-1})\bigr)^2|{\cal
F}_{t_{k-1}}\bigr] P(B|{\cal F}_{t_{k-1}})}\nonumber
\\[-8pt]
\\[-8pt]
&& \qquad \leq\sqrt{C \theta N^2 (\log N)^4 e^{-A} \eps^{-1/2} \cdot\theta
\delta\eps^{-5/2} \bigl(1 + o(1)\bigr)} \nonumber\\
&& \qquad \leq C \theta e^{-A/2} \delta^{1/2} \eps^{-3/2} N (\log N)^2 \bigl(1 +
o(1)\bigr).
\nonumber
\end{eqnarray}
Likewise, by the conditional Cauchy--Schwarz inequality and Lemma~\ref{Wvar},
%
%e129 ###
%
%e130 #&#
\begin{eqnarray}\label{CSW}
&&E \Biggl[ \Biggl| \pi\sqrt{2} e^{-A} N (\log N)^2 \sum_{j=1}^{R_k}
W_{j,k} \mathbf{1}_{\{W_{j,k} \leq\eps/(\pi\sqrt{2} e^{-A})\}} \Biggr|
\mathbf{1}_B \Big| {\cal F}_{t_{k-1}} \Biggr] \nonumber\\
&& \qquad\leq C e^{-A} N (\log N)^2 \sqrt{\theta e^{2A} \eps^{1/2}
\cdot\theta
\delta\eps^{-5/2} \bigl(1 + o(1)\bigr)} \\
&& \qquad\leq C \theta N (\log N)^2 \delta^{1/2} \eps^{-1} \bigl(1 +
o(1)\bigr).\nonumber
\end{eqnarray}
Now Lemma~\ref{Zprimelem}, (\ref{CSeq}) and (\ref{CSW}) imply
%
%e130 ###
%
%e131 #&#
\begin{eqnarray}\label{abs3bound}
  E [ |S_k| \mathbf{1}_B | {\cal F}_{t_{k-1}} ] &\leq& C \theta N
(\log N)^2 \delta^{1/2}  \nonumber
\\[-8pt]
\\[-8pt]&&{}\times(e^{-A/2} \eps^{-3/2}+ \eps^{-1}) + o (N
(\log N)^2 ).
\nonumber
\end{eqnarray}
The result follows from (\ref{abs1bound}), (\ref{abs2bound}) and
(\ref
{abs3bound}) in view of the inequality (\ref{652}) and $\delta\le
\eps
$, as well as (\ref{655}).
\end{pf}

%
%pr39 #&#
\begin{Prop}\label{smalljumpexp}
There exists a real number $c$ such that
\begin{eqnarray*}
&&E \bigl[\bigl(Z_N(t_k) - Z_N(t_{k-1})\bigr) \mathbf{1}_{\{Z_N(t_k) - Z_N(t_{k-1})
\leq\eps N (\log N)^2\}} | {\cal F}_{t_{k-1}} \bigr] \\
&& \qquad= Z_N(t_{k-1}) \theta s \bigl(c + 2 \pi^2 \log\eps+ g(\epsilon
, A)\bigr) \\
&& \qquad {}+ O
 (\theta N (\log N)^2 \delta^{1/2} \eps^{-3} ) + o(N(\log N)^2),
\end{eqnarray*}
where $g\dvtx  (0, \infty) \times(0, \infty) \rightarrow\R$ is a function
such that $\lim_{y \rightarrow\infty} g(x,y) = 0$ for all $x > 0$.
\end{Prop}

\begin{pf}
By combining Lemmas~\ref{Rexp} and~\ref{Rerr} and using (\ref{657}),
we get
%
%e131 ###
%
%e132 #&#
\begin{eqnarray}
&&E \bigl[\bigl(Z_N(t_k) - Z_N(t_{k-1})\bigr) \mathbf{1}_{\{Z_N(t_k) - Z_N(t_{k-1})
\leq\eps N (\log N)^2\}} | {\cal F}_{t_{k-1}} \bigr] \nonumber
\\
&& \qquad= Z_N(t_{k-1}) \theta s \bigl(2 \sqrt{2} \pi^2 E\bigl[W \mathbf
{1}_{\{W \leq\eps
/(\pi\sqrt{2} e^{-A})\}}\bigr] - 2 \pi^2 A\bigr)\\
&& \qquad  \quad {} + O (\theta N (\log N)^2
\delta
^{1/2} \eps^{-3} ) + o(N (\log N)^2).
\nonumber
\end{eqnarray}
Denote the conditional expectation on the left-hand side of this
equation by $f(N, \eps, \theta)$. Note that this expectation depends on
$N$, $\eps$ and $\theta$, but can not depend on $\delta$ or $A$, as
these constants were introduced just for the proof. Assume for the
moment that $k = 1$, and the initial conditions are chosen so that
$Z_N(0) = N (\log N)^2$.
Then there exists a positive constant $C$ such that
\begin{eqnarray*}
&&\limsup_{\theta\rightarrow0} \limsup_{N \rightarrow\infty} \frac
{f(N, \eps, \theta)}{N (\log N)^2 \theta s} \\
&& \qquad \leq\bigl(2 \sqrt{2} \pi^2 E\bigl[W
\mathbf{1}_{\{W \leq\eps/(\pi\sqrt{2} e^{-A})\}}\bigr] - 2 \pi^2 A\bigr) +
C \delta
^{1/2} \eps^{-3}
\end{eqnarray*}
and likewise
\begin{eqnarray*}
&&\liminf_{\theta\rightarrow0} \liminf_{N \rightarrow\infty} \frac
{f(N, \eps, \theta)}{N (\log N)^2 \theta s} \\
&& \qquad \geq\bigl(2 \sqrt{2} \pi^2 E\bigl[W
\mathbf{1}_{\{W \leq\eps/(\pi\sqrt{2} e^{-A})\}}\bigr] - 2 \pi^2 A\bigr) -
C \delta
^{1/2} \eps^{-3}.
\end{eqnarray*}
We now simultaneously take $\delta\rightarrow0$ and $A \rightarrow
\infty$. This can be done without violating the constraints on the
constants because once $\delta$ is chosen, we can pick $A$ large enough
to satisfy (\ref{651}) and (\ref{652}), and then only consider
$\theta$
small enough that (\ref{8613})--(\ref{658}) are satisfied.
The second term $C\delta^{1/2} \eps^{-3}$ then tends to zero. Since the
left-hand side does not depend on $A$,
the first term must also tend to a limit as $A \rightarrow\infty$.
That is, we know that
%
%e132 ###
%
%e133 #&#
\begin{equation}\label{exlim}
\lim_{A \rightarrow\infty} \bigl(2 \sqrt{2} \pi^2 E\bigl[W \mathbf{1}_{\{W
\leq\eps
/(\pi\sqrt{2} e^{-A})\}}\bigr] - 2 \pi^2 A\bigr) \qquad\mbox{exists}.
\end{equation}
Now let $r = \eps/(\pi\sqrt{2} e^{-A})$, so
\[
A = \log\biggl( \frac{\pi\sqrt{2} r}{\eps} \biggr) = \log\bigl(\pi\sqrt{2}
r\bigr) - \log\eps.
\]
Therefore, the limit in (\ref{exlim}) is equal to
%
%e133 ###
%
%e134 #&#
\begin{equation}\label{findB}\qquad
\lim_{r \rightarrow\infty} 2 \pi^2 \bigl(\sqrt{2} E\bigl[W \mathbf{1}_{\{W
\leq
r\}}\bigr] - \log\bigl(\pi\sqrt{2} r\bigr) + \log\eps\bigr) = c + 2 \pi^2 \log\eps
\end{equation}
for some real number $c$ that does not depend on $\eps$. The
proposition follows.
\end{pf}

%
%re40 #&#
\begin{Rmk}\rm
Equation \eqref{exlim} is a statement which concerns only critical
branching Brownian motion with absorption and does not depend on $N$.
It would be desirable\vadjust{\goodbreak} to find a direct proof of this fact, but we were
not able to obtain one. This would follow if one could show that
\[
\int_1^\infty\biggl| P(W>x) - \frac1{\sqrt{2}x} \biggr|\,dx < \infty.
\]
An explicit expression for the value of the limit in \eqref{exlim}
would also make it possible to identify the constant $a$ appearing in
the statement of Proposition~\ref{ZNCSBP}.
\end{Rmk}

%s5.4 ###
%s5.4 #&#
\subsection{The large jumps}\label{sec54}

We now estimate the probability in (\ref{bigjeq}) that the process
$Z_N$ makes a large jump between times $t_{k-1}$ and $t_k$.

%
%pr41 #&#
\begin{Prop}\label{rjumplem}
For all $r \geq\eps$, on $G_{N, k-1}$ we have
\begin{eqnarray*}
&&P\bigl(Z_N(t_k) - Z_N(t_{k-1}) > r N (\log N)^2|{\cal F}_{t_{k-1}}\bigr)\\
&& \qquad  = \frac
{2 \pi^2 \theta s}{r} \cdot\frac{Z_N(t_{k-1})}{N (\log N)^2} + O
(\theta\delta\eps^{-5/2}) + o(1).
\end{eqnarray*}
\end{Prop}

\begin{pf}
By Lemma~\ref{B6Lem}, we have
\begin{eqnarray*}
&&P\bigl(Z_N(t_k) - Z_N(t_{k-1}) > r N (\log N)^2|{\cal F}_{t_{k-1}}\bigr)
\\
&& \qquad= P\bigl(W_{j,k} > r/\bigl(\pi\sqrt{2} e^{-A}\bigr) \mbox{ for some
}j|{\cal
F}_{t_{k-1}}\bigr) + O(\theta\delta\eps^{-5/2}) + o(1).
\end{eqnarray*}
Recall that $\epsilon/(\pi\sqrt{2} e^{-A}) \geq e^{2A/3}$ by (\ref
{652}) and the fact that $\delta\in(0, \epsilon)$. By Lemma~\ref
{B1Lem}, for sufficiently large $A$ the probability that $W_{j_1,k} >
r/(\pi\sqrt{2} e^{-A})$ and $W_{j_2,k} > r/(\pi\sqrt{2} e^{-A})$ for
some $j_1 \neq j_2$ is at most $C \theta e^{-A/3} \eps^{-1/2} + o(1)$.
Therefore,
\begin{eqnarray*}
&&P\bigl(W_{j,k} > r/\bigl(\pi\sqrt{2} e^{-A}\bigr) \mbox{ for some }j|{\cal
F}_{t_{k-1}}\bigr) \\
&& \qquad= E[R_k|{\cal F}_{t_{k-1}}] P\bigl(W > r/\bigl(\pi\sqrt{2} e^{-A}\bigr)\bigr) +
O(\theta
e^{-A/3} \eps^{-1/2}) + o(1),
\end{eqnarray*}
and the error term is smaller than $O(\theta\delta\eps^{-5/2})$ by
(\ref{652}).

By Proposition~\ref{Wprop}, if we use $\sim$ to mean that the ratio of
the two sides tends to one
as $x \rightarrow\infty$, then
\begin{eqnarray*}
E\bigl[W \mathbf{1}_{\{W \leq x\}}\bigr] &=& \int_0^x P(y \leq W \leq x) \, dy\\
 &=&
\int
_0^x P(W \geq y) \, dy - x P(W > x) \sim B \log x.
\end{eqnarray*}
Therefore, (\ref{findB}) implies that $B = 1/\sqrt{2}$. Therefore, by
(\ref{WBbound}),
\[
\frac{(1 - \delta) \pi}{r e^A} \leq P \biggl( W > \frac{r}{\pi\sqrt{2}
e^{-A}} \biggr) \leq\frac{(1 + \delta) \pi}{r e^A}.
\]
Combining this result with Proposition~\ref{EMkProp}, we get on $G_{N, k-1}$,
%
%e134 ###
%
%e135 #&#
\begin{eqnarray}\label{Rkend}
&&E[R_k|{\cal F}_{t_{k-1}}] P\bigl(W > r/\bigl(\pi\sqrt{2} e^{-A}\bigr)\bigr) \nonumber
\\[-8pt]
\\[-8pt]&& \qquad = \frac{2
\pi
^2 \theta s}{r} \cdot\frac{Z_N(t_{k-1})}{N (\log N)^2} \cdot\bigl(1 + O(A
\theta)\bigr)\bigl(1 + O(\delta)\bigr) + o(1),
\nonumber
\end{eqnarray}
which is enough to imply the result. Since $1/r \leq\eps^{-1}$ and
$Z_N(t_{k-1})/\break N(\log N)^2 \leq\eps^{-1/2}$ on $G_{N, k-1}$, the
dominant error term coming from (\ref{Rkend}) is $O(\theta\delta\eps
^{-3/2})$.
\end{pf}

%s6 ###
%s6 #&#
\section{Convergence to the CSBP}\label{sec6}

In this section, we prove Proposition~\ref{ZNCSBP} and Theorem \ref
{MNCSBP}. Both of these results require proving that a sequence of
processes converges to the continuous-state branching process $(Z(t), t
\geq0)$ with branching mechanism
\[
\Psi(u) = a u + 2 \pi^2 u \log u = - c u + 2 \pi^2 \int_0^{\infty}
\bigl(e^{-ux} - 1 + ux \mathbf{1}_{\{x \leq1\}}\bigr) x^{-2} \, dx,
\]
where $c$ is the constant defined in (\ref{findB}).
We will first establish Proposition~\ref{ZNCSBP}, and then use this
result to deduce Theorem~\ref{MNCSBP}.

%s6.1 ###
%s6.1 #&#
\subsection{The generator of the CSBP}\label{sec61}

Let $C_0([0, \infty))$ be the set of continuous functions $f\dvtx  [0,
\infty
) \rightarrow\R$ that vanish at infinity, endowed with the sup norm so
that for $f \in C_0([0, \infty))$, we have
\[
\|f\| = \sup_{x \geq0} |f(x)|.
\]
For $f \in C_0([0, \infty))$ and $x \in[0, \infty)$, let $T_tf(x) =
E[f(Z(t))|Z(0) = x]$.
It is well-known (see, e.g.,~\cite{cablam}) that $(T_t, t \geq0)$ is a
Feller semigroup. The following result describes the associated
infinitesimal generator. This result is essentially well-known. The
form of the generator appeared in~\cite{silv}, and later in~\cite{dk99}
where a particle representation of continuous-state branching processes
was constructed. The fact that the set ${\cal E}$ defined below is a
core for the generator was established for closely related families of
processes in~\cite{li06,ma09}. However, we give a short proof of the
result below for completeness.

%
%pr42 #&#
\begin{Prop}\label{gencsbp}
Let $A$ be the infinitesimal generator for $(Z(t), t \geq0)$. Let
${\cal E} \subset C_0([0, \infty))$ be the set of functions of the form
%
%e135 ###
%
%e136 #&#
\begin{equation}\label{finE}
f(x) = a_1 e^{-\lambda_1 x} +\cdots+ a_m e^{-\lambda_m x},
\end{equation}
where $a_1,\ldots, a_m \in\R$ and $\lambda_1,\ldots, \lambda_m > 0$.
Then ${\cal E}$ is a core for $A$, and for $f \in{\cal E}$,
%
%e136 ###
%
%e137 #&#
\begin{equation}\label{csbpgen}
   Af(x) = x \biggl(c f'(x) + 2 \pi^2 \int_0^{\infty} \bigl( f(x + y) - f(x)
- y \mathbf{1}_{\{y \leq1\}} f'(x) \bigr) y^{-2} \, dy \biggr).\hspace*{-35pt}
\end{equation}
\end{Prop}

\begin{pf}
If $f(x) = e^{-\lambda x}$, then by (\ref{lapcsbp}) and (\ref
{csbpdiffeq}), we have
\begin{eqnarray*}
Af(x) &=& \lim_{t \rightarrow0} \frac{T_tf(x) - f(x)}{t} = \lim_{t
\rightarrow0} \frac{e^{-x u_t(\lambda)} - e^{-\lambda x}}{t}\\
 &=& \frac
{\partial}{\partial t} e^{-x u_t(\lambda)} \bigg|_{t = 0} = x
e^{-\lambda x} \Psi(\lambda),
\end{eqnarray*}
which equals the right-hand side of (\ref{csbpgen}). The result (\ref
{csbpgen}) then follows for all $f \in{\cal E}$ by linearity. By the
Stone--Weierstrass theorem, ${\cal E}$ is dense in $C_0([0, \infty))$.
By (\ref{lapcsbp}), we have $T_tf \in{\cal E}$ whenever $f \in{\cal
E}$. It now follows from Proposition 3.3 in Chapter 1 of~\cite{ek86}
that ${\cal E}$ is a core for $A$.
\end{pf}

%s6.2 ###
%s6.2 #&#
\subsection{\texorpdfstring{Proof of Proposition \protect\ref{ZNCSBP}}{Proof of Proposition 1}}\label{sec62}

The next result is Theorem 8.2 in Chapter~4 of~\cite{ek86} in the
present context.

%
%pr43 #&#
\begin{Prop}
Suppose the distribution of $V_N(0)$ converges to the distribution of
$Z(0)$ as $N \rightarrow\infty$. Then the finite-dimensional
distributions of $(V_N(t),  t \geq0)$ converge to those of $(Z(t), t
\geq0)$ as $N \rightarrow\infty$ if and only if for all $j \geq0$,
all $0 \leq s_1 < s_2 <\cdots< s_j \leq u < u + s$, all bounded
continuous functions $h_1,\ldots, h_j\dvtx  [0, \infty) \rightarrow\R$, and
all $f \in{\cal E}$, we have
%
%e137 ###
%
%e138 #&#
\begin{eqnarray}\label{mainVN}
&&\lim_{N \rightarrow\infty} E \biggl[ \biggl( f\bigl(V_N(u + s)\bigr) - f(V_N(u))\nonumber
\\[-8pt]
\\[-8pt]&&\hspace*{18.5pt}\hphantom{\lim_{N \rightarrow\infty} E \biggl[ \biggl(}{} -
\int_u^{u + s} Af(V_N(t)) \, dt \biggr) \prod_{i=1}^j h_i(V_N(s_i))
\biggr] = 0.
\nonumber
\end{eqnarray}
\end{Prop}

In view of this result, we will aim to establish (\ref{mainVN}), which
will imply Proposition~\ref{ZNCSBP}. We will assume that $0 \leq s_1 <
s_2 <\cdots< s_j \leq u < u + s$. We also define the times
\[
u = \tau_0 < \tau_1 <\cdots< \tau_{\theta^{-1}} = u+s,
\]
where $\tau_k = t_k/(\log N)^3$ for all $k$.
This means that $V_N(\tau_k) = Z_N(t_k)/(N(\log N)^2)$ for all $k$. We
also assume that the function $f \in{\cal E}$ and the bounded
continuous functions $h_1,\ldots, h_j$ are fixed throughout this subsection.

Since $f$ is of the form given in (\ref{finE}), the norms $\|f\|$, $\|
f'\|$ and $\|f''\|$ are finite and thus can be treated as constants. If
$g(x) = x f(x)$ and $d(x) = x f'(x)$, then $\|g\|$, $\|g'\|$ and $\|d\|
$ are likewise finite.
Also, if we define
%
%e138 ###
%
%e139 #&#
\begin{equation}\label{xfxy}
h(x) = \sup_{y \geq x} x |f''(y)|, \qquad  k(x) = \sup_{y \geq x}
x |f(y)|,
\end{equation}
then it is easy to check that $\|h\| < \infty$ and $\|k\| < \infty$.
Finally, if $y \geq0$, then by Taylor's theorem there is a $z \in[x,
x+y]$ such that $f(x+y) = f(x) + yf'(x) + \frac{1}{2} y^2 f''(z)$. Therefore,
\[
\biggl| x \int_0^1 \bigl( f(x+y) - f(x) - y f'(x) \bigr) y^{-2} \, dy
\biggr| \leq\frac{1}{2} |h(x)|
\]
and
\[
\biggl| x \int_1^{\infty} \bigl( f(x+y) - f(x) \bigr) y^{-2} \, dy \biggr|
\leq|k(x)| + |g(x)|.
\]
It follows that
%
%e139 ###
%
%e140 #&#
\begin{equation}\label{Afbound}
\|Af\| \leq|c| \|d\| + 2 \pi^2 \bigl( \|g\| + \|k\| + \tfrac{1}{2} \|
h\| \bigr) < \infty.
\end{equation}

%
%le44 #&#
\begin{Lemma}\label{lem41}
We have
%
%e140 ###
%
%e141 #&#
\begin{eqnarray}
&&  E \bigl[ \bigl( f(V_N(\tau_k)) - f(V_N(\tau_{k-1})) \bigr) \mathbf{1}_{\{
V_N(\tau_k) - V_N(\tau_{k-1}) \leq\eps\}} | {\cal F}_{t_{k-1}}
\bigr] \mathbf{1}_{G_{N, k-1}} \nonumber\hspace*{-35pt}
\\[-8pt]
\\[-8pt]
&&   \quad= f'(V_N(\tau_{k-1})) V_N(\tau_{k-1}) \theta s (c + 2 \pi
^2 \log\eps)
\mathbf{1}_{G_{N, k-1}} + O(\theta\eps^{1/2}) + o(1).
\nonumber\hspace*{-35pt}
\end{eqnarray}
\end{Lemma}

\begin{pf}
Define
\[
{\bar S}_k = \frac{Z_{N,1}'(t_k) - Z_N(t_{k-1})}{N (\log N)^2} + \pi
\sqrt{2} e^{-A} \sum_{j=1}^{R_k} W_{j,k} \mathbf{1}_{\{W_{j,k} \leq
\eps
/(\pi\sqrt{2} e^{-A})\}}.
\]
Note that ${\bar S}_k$ would be equal to $S_k/(N (\log N)^2)$, where
$S_k$ is defined in (\ref{Skdef}), if $Z_{N,1}'(t_k)$ were replaced in
the definition by $Z_{N,1}(t_k)$. Therefore, by Lemma~\ref{Zprimelem},
%
%e141 ###
%
%e142 #&#
\begin{equation}\label{SkSkbar}
E \biggl[ \biggl|{\bar S}_k - \frac{S_k}{N (\log N)^2} \biggr|\Big | {\cal
F}_{t_{k-1}} \biggr] \mathbf{1}_{G_{N, k-1}} = o(1).
\end{equation}
Thus, by Lemma~\ref{Rerr},
%
%e142 ###
%
%e143 #&#
\begin{eqnarray}\label{newrefeq}
&&E \bigl[ \bigl| {\bar S}_k - \bigl(V_N(\tau_k) - V_N(\tau_{k-1})\bigr)\bigr| \mathbf{1}_{\{
V_N(\tau_k) - V_N(\tau_{k-1}) \leq\eps\}} | {\cal F}_{t_{k-1}}
\bigr] \mathbf{1}_{G_{N, k-1}} \nonumber
\\[-8pt]
\\[-8pt]&& \qquad \leq C \theta\delta^{1/2} \eps^{-3} + o(1).
\nonumber
\end{eqnarray}
It follows from (\ref{newrefeq}) that
%
%e143 ###
%
%e144 #&#
\begin{eqnarray}\label{411}
&&  E \bigl[ \bigl|f(V_N(\tau_k)) - f\bigl(V_N(\tau_{k-1}) + {\bar S}_k\bigr)\bigr| \mathbf
{1}_{\{
V_N(\tau_k) - V_N(\tau_{k-1}) \leq\eps\}} | {\cal F}_{t_{k-1}}
\bigr] \mathbf{1}_{G_{N, k-1}} \nonumber\hspace*{-35pt}
\\[-8pt]
\\[-8pt]
&& \quad  \leq C \| f' \| \theta\delta^{1/2} \eps^{-3} + o(1) \leq C
\theta
\delta^{1/2} \eps^{-3} + o(1).
\nonumber\hspace*{-35pt}
\end{eqnarray}
By Taylor's theorem, there exists $\xi$ between $V_N(\tau_{k-1})$ and
$V_N(\tau_{k-1}) + {\bar S}_k$
such that
%
%e144 ###
%
%e145 #&#
\begin{eqnarray}\label{412}
&&E \bigl[\bigl ( f\bigl(V_N(\tau_{k-1}) + {\bar S}_k\bigr) - f(V_N(\tau_{k-1}))
\bigr) \mathbf{1}_{\{V_N(\tau_k) - V_N(\tau_{k-1}) \leq\eps\}} | {\cal
F}_{t_{k-1}} \bigr] \mathbf{1}_{G_{N, k-1}} \nonumber\hspace*{-35pt}\\
&& \qquad= E \bigl[ \bigl( f'(V_N(\tau_{k-1})) {\bar S}_k + f''(\xi) {\bar
S}_k^2/2 \bigr) \mathbf{1}_{\{V_N(\tau_k) - V_N(\tau_{k-1}) \leq\eps\}}
| {\cal F}_{t_{k-1}} \bigr] \mathbf{1}_{G_{N, k-1}}\nonumber\hspace*{-35pt}
\\[-8pt]
\\[-8pt]
&& \qquad= f'(V_N(\tau_{k-1})) E\bigl[{\bar S}_k \mathbf{1}_{\{V_N(\tau
_k) - V_N(\tau
_{k-1}) \leq\eps\}}|{\cal F}_{t_{k-1}}\bigr] \mathbf{1}_{G_{N, k-1}}\nonumber\hspace*{-35pt}\\
&& \qquad  \quad {} +
O(E[{\bar S}_k^2|{\cal F}_{t_{k-1}}] \mathbf{1}_{G_{N,
k-1}}).\nonumber\hspace*{-35pt}
\end{eqnarray}
Lemma~\ref{Rerr} and (\ref{SkSkbar}) give
%
%e145 ###
%
%e146 #&#
\begin{equation}\label{newww}\qquad
E\bigl[|{\bar S}_k| \mathbf{1}_{\{V_N(\tau_k) - V_N(\tau_{k-1}) > \eps\}
}|{\cal
F}_{t_{k-1}}\bigr] \mathbf{1}_{G_{N, k-1}} \leq C \theta\delta^{1/2} \eps^{-3}
+ o(1).
\end{equation}
Note that $\delta^{1/2} \eps^{-3} \leq\eps^{1/2}$ by (\ref
{deltasmall}), and $A$ can be chosen large enough so that $g(\epsilon,
A) \leq\eps$, where $g$ is the function from Proposition \ref
{smalljumpexp}. Therefore, (\ref{newww}) combined with Lemma \ref
{Rerr}, equation (\ref{SkSkbar}), and Proposition~\ref{smalljumpexp} implies
%
%e146 ###
%
%e147 #&#
\begin{eqnarray}\label{413}
&&  f'(V_N(\tau_{k-1})) E\bigl[{\bar S}_k \mathbf{1}_{\{V_N(\tau_k) -
V_N(\tau
_{k-1}) \leq\eps\}}|{\cal F}_{t_{k-1}}\bigr] \mathbf{1}_{G_{N, k-1}}
\nonumber\hspace*{-35pt}
\\
&&   \qquad= f'(V_N(\tau_{k-1})) E[{\bar S}_k|{\cal F}_{t_{k-1}}]
\mathbf{1}_{G_{N,
k-1}} + O(\theta\eps^{1/2}) + o(1) \hspace*{-35pt}\\
&&   \qquad= f'(V_N(\tau_{k-1})) V_N(\tau_{k-1}) \theta s (c + 2 \pi
^2 \log\eps)
\mathbf{1}_{\{G_{N, k-1}\}} + O(\theta\eps^{1/2}) + o(1).\nonumber\hspace*{-35pt}
\end{eqnarray}
Since $e^{-A} \eps^{-1/2} \leq\eps^{1/2}$ by (\ref{deltasmall}) and
(\ref{652}), it follows from Lemmas~\ref{Zvar} and~\ref{Wvar} that
%
%e147 ###
%
%e148 #&#
\begin{eqnarray}\label{414}
E[{\bar S}_k^2|{\cal F}_{t_{k-1}}] \mathbf{1}_{G_{N, k-1}} &\leq&
C(\theta
e^{-A} \eps^{-1/2} + e^{-2A} \cdot\theta e^{2A} \eps^{1/2}) +
o(1)\nonumber
\\[-8pt]
\\[-8pt]
&\leq& C \theta\eps^{1/2} + o(1).
\nonumber
\end{eqnarray}
The result follows from (\ref{411}), (\ref{412}), (\ref{413}) and
(\ref{414}).
\end{pf}

%
%le45 #&#
\begin{Lemma}\label{lem42}
We have
\begin{eqnarray*}
&&E \bigl[\bigl ( f(V_N(\tau_k)) - f(V_N(\tau_{k-1})) \bigr) \mathbf{1}_{\{
V_N(\tau_k) - V_N(\tau_{k-1}) > \eps\}} | {\cal F}_{t_{k-1}} \bigr]
\mathbf{1}_{G_{N, k-1}} \\
&& \qquad= 2 \pi^2 \theta s V_N(\tau_{k-1}) \mathbf{1}_{G_{N, k-1}}
\int_{\eps
}^{\infty}\bigl ( f\bigl(V_N(\tau_{k-1}) + y\bigr) - f(V_N(\tau_{k-1})) \bigr)
y^{-2} \, dy \\
&& \qquad  \quad {}+ O(\theta\eps^{1/2}) + o(1).
\end{eqnarray*}
\end{Lemma}

\begin{pf}
By Proposition~\ref{rjumplem} with $r = \eps$,
%
%e148 ###
%
%e149 #&#
\begin{eqnarray}\label{421}
&&E\bigl[f(V_N(\tau_{k-1})) \mathbf{1}_{\{V_N(\tau_k) - V_N(\tau_{k-1})
> \eps\}
} | {\cal F}_{t_{k-1}}\bigr] \mathbf{1}_{G_{N, k-1}} \nonumber\hspace*{-35pt}\\
&& \quad= f(V_N(\tau_{k-1})) P\bigl(V_N(\tau_k) - V_N(\tau_{k-1}) >
\eps| {\cal
F}_{t_{k-1}}\bigr) \mathbf{1}_{G_{N, k-1}} \nonumber\hspace*{-35pt}
\\[-8pt]
\\[-8pt]
&& \quad= f(V_N(\tau_{k-1})) V_N(\tau_{k-1}) \cdot\frac{2 \pi^2
\theta s}{\eps
} \mathbf{1}_{G_{N, k-1}} + O(\theta\delta\eps^{-5/2}) + o(1)
\nonumber\hspace*{-35pt}
\\
&& \quad= 2 \pi^2 \theta s V_N(\tau_{k-1}) \mathbf{1}_{G_{N, k-1}}
\int_{\eps
}^{\infty} f(V_N(\tau_{k-1})) y^{-2} \, dy + O(\theta\delta\eps
^{-5/2}) + o(1).
\nonumber\hspace*{-35pt}
\end{eqnarray}
To simplify notation, assume that $\eps^{-1}$ is an integer. Then
%
%e149 ###
%
%e150 #&#
\begin{eqnarray}\label{sumM}
&&E\bigl[f(V_N(\tau_k)) \mathbf{1}_{\{V_N(\tau_k) - V_N(\tau_{k-1}) >
\eps\}} |
{\cal F}_{t_{k-1}}\bigr] \mathbf{1}_{G_{N, k-1}} \nonumber\\
&& \qquad= \sum_{m=\eps^{-1}}^{\infty} E\bigl[f(V_N(\tau_k)) \mathbf
{1}_{\{m \eps^2 <
V_N(\tau_k) - V_N(\tau_{k-1}) \leq(m+1)\eps^2\}}|{\cal F}_{t_{k-1}}\bigr]
\mathbf{1}_{G_{N, k-1}} \nonumber\\
&& \qquad= \sum_{m=\eps^{-1}}^{\eps^{-3} - 1} f\bigl(\eps^2 m +
V_N(\tau_{k-1})\bigr)\nonumber \\
&& \hphantom{= \sum_{m=\eps^{-1}}^{\eps^{-3} - 1} }\qquad {}\times P\bigl(m
\eps^2 < V_N(\tau_k) - V_N(\tau_{k-1}) \leq(m+1)\eps^2|{\cal
F}_{t_{k-1}}\bigr)\mathbf{1}_{G_{N, k-1}} \\
&& \qquad \quad{} + \sum_{m=\eps^{-1}}^{\eps^{-3} - 1} E \bigl[\bigl (
f(V_N(\tau_k)) - f\bigl(\eps^2 m + V_N(\tau_{k-1})\bigr) \bigr) \nonumber\\
&& \hphantom{\quad{} + \sum_{m=\eps^{-1}}^{\eps^{-3} - 1} E \bigl[}\qquad {}\times\mathbf{1}_{\{m
\eps
^2 < V_N(\tau_k) - V_N(\tau_{k-1}) \leq(m+1)\eps^2\}} |{\cal
F}_{t_{k-1}} \bigr] \mathbf{1}_{G_{N, k-1}} \nonumber\\
&& \qquad \quad{} + E\bigl[f(V_N(\tau_k)) \mathbf{1}_{\{V_N(\tau_k) -
V_N(\tau
_{k-1}) > \eps^{-1}\}}| {\cal F}_{t_{k-1}}\bigr] \mathbf{1}_{G_{N,
k-1}}.\nonumber
\end{eqnarray}
Denote the three terms on the right-hand side of (\ref{sumM}) by $T_1$,
$T_2$ and $T_3$. Proposition~\ref{rjumplem} gives
%
%e150 ###
%
%e151 #&#
\begin{eqnarray}\label{2tmbd}
|T_2| &\leq&\eps^2 \|f'\| P\bigl(V_N(\tau_k) - V_N(\tau_{k-1}) > \eps
|{\cal
F}_{t_{k-1}}\bigr)\mathbf{1}_{G_{N, k-1}} \nonumber
\\
&\leq& C \theta\eps V_N(\tau_{k-1}) \mathbf{1}_{G_{N, k-1}} +
O(\theta
\delta\eps^{-1/2}) + o(1) \\
&\leq &C \theta(\eps^{1/2} + \delta\eps
^{-1/2}) + o(1)
\nonumber
\end{eqnarray}
and
%
%e151 ###
%
%e152 #&#
\begin{eqnarray}\label{3tmbd}
|T_3| &\leq&\|f\| P\bigl(V_N(\tau_k) - V_N(\tau_{k-1}) > \eps^{-1}
|{\cal F}_{t_{k-1}}\bigr) \mathbf{1}_{G_{N, k-1}}\nonumber
\\
&\leq& C \theta\eps V_N(\tau_{k-1}) \mathbf{1}_{G_{N, k-1}} +
O(\theta
\delta\eps^{-5/2}) + o(1)\\
& \leq &C \theta(\eps^{1/2}
+ \delta\eps^{-5/2}) + o(1).
\nonumber
\end{eqnarray}
By Proposition~\ref{rjumplem} and the fact that
\[
\frac{1}{\eps^2} \biggl( \frac{1}{m} - \frac{1}{m+1} \biggr) = \frac
{1}{\eps^2 m (m+1)},
\]
we have
\begin{eqnarray*}
&&P\bigl(m \eps^2 < V_N(\tau_k) - V_N(\tau_{k-1})  \leq (m+1)\eps^2|{\cal
F}_{t_{k-1}}\bigr) \mathbf{1}_{G_{N, k-1}} \\
&& \qquad = \frac{2 \pi^2 \theta s V_N(\tau_{k-1})}{\eps^2 m (m+1)} \mathbf
{1}_{G_{N, k-1}} + O(\theta\delta\eps^{-5/2}) + o(1).
\end{eqnarray*}
Adding up at most $\eps^{-3}$ error terms of order $\theta\delta\eps
^{-5/2}$ to get a single error term of order $\theta\delta\eps
^{-11/2}$, we get
%
%e152 ###
%
%e153 #&#
\begin{eqnarray}\label{1tmbdpre}
T_1 &=& 2 \pi^2 \theta s V_N(\tau_{k-1}) \mathbf{1}_{G_{N, k-1}}
\sum_{m =
\eps^{-1}}^{\eps^{-3} - 1} \frac{f(\eps^2 m + V_N(\tau
_{k-1}))}{\eps^2
m(m+1)} \nonumber\\
&&{}+ O(\theta\delta\eps^{-11/2}) + o(1) \nonumber
\\[-8pt]
\\[-8pt]
&=& 2 \pi^2 \theta s V_N(\tau_{k-1}) \mathbf{1}_{G_{N, k-1}}  \sum
_{m = \eps
^{-1}}^{\eps^{-3} - 1} f\bigl(\eps^2 m + V_N(\tau_{k-1})\bigr) \int_{\eps^2
m}^{\eps^2(m+1)} y^{-2} \, dy\nonumber\\
&&{} + O(\theta\delta\eps^{-11/2}) + o(1).
\nonumber
\end{eqnarray}
Because an error of at most $\|f'\| \eps^2$ is made when replacing
$f(\eps^2 m + V_N(\tau_{k-1}))$ by $f(V_N(\tau_{k-1}) + y)$ with
$\eps
^2 m \leq y \leq\eps^2 (m+1)$, we get
\begin{eqnarray*}
&&\Biggl| \sum_{m=\eps^{-1}}^{\eps^{-3} - 1} f\bigl(\eps^2m + V_N(\tau_{k-1})\bigr)
\int_{\eps^2 m}^{\eps^2 (m+1)} y^{-2} \, dy -\int_{\eps}^{\infty}
f\bigl(V_N(\tau_{k-1}) + y\bigr) y^{-2} \, dy \Biggr| \\
&& \qquad\leq C \biggl( \int_{\eps}^{\eps^{-1}} \eps^2 y^{-2} \, dy +
\int_{\eps
^{-1}}^{\infty} y^{-2} \, dy \biggr) \leq C \eps.
\end{eqnarray*}
Combining this with (\ref{1tmbdpre}) gives
%
%e153 ###
%
%e154 #&#
\begin{eqnarray}\label{1tmbd}
T_1 &=& 2 \pi^2 \theta s V_N(\tau_{k-1}) \mathbf{1}_{G_{N, k-1}} \int
_{\eps
}^{\infty} f\bigl(V_N(\tau_{k-1}) + y\bigr) y^{-2} \, dy\nonumber
\\[-8pt]
\\[-8pt] &&{}+ O(\theta\eps
^{1/2}) +
O(\theta\delta\eps^{-11/2}) + o(1).
\nonumber
\end{eqnarray}
Since we have chosen $\delta\leq\eps^7 \leq\eps^6$ by (\ref
{deltasmall}), the lemma now follows by summing (\ref{2tmbd}), (\ref
{3tmbd}) and (\ref{1tmbd}) and subtracting (\ref{421}) from the result.
\end{pf}

Note that
\[
\int_{\eps}^{\infty} y \mathbf{1}_{\{y \leq1\}} f'(x) y^{-2} \, dy =
f'(x) \int_{\eps}^1 y^{-1} \, dy = - f'(x) \log\eps.
\]
Therefore, for every $\eps>0$ one can write
\[
Af(x) = A_1 f(x) + A_2 f(x),
\]
where
%
%e154 ###
%
%e155 #&#
\begin{eqnarray}\label{A1def}
A_1 f(x) &=& x \biggl( (c + 2 \pi^2 \log\eps) f'(x) \nonumber
\\[-8pt]
\\[-8pt]&&\hphantom{x \biggl(}{}+ 2 \pi^2 \int_0^{\eps
} \bigl( f(x + y) - f(x) - y f'(x) \bigr) y^{-2} \, dy \biggr)
\nonumber
\end{eqnarray}
and
\[
A_2 f(x) = x\biggl ( 2 \pi^2 \int_{\eps}^{\infty} \bigl( f(x + y) - f(x)
\bigr) y^{-2} \, dy \biggr).
\]

%
%le46 #&#
\begin{Lemma}\label{uniformlem}
On $G_{N, k-1}$, we have
\[
E \biggl[ \int_{\tau_{k-1}}^{\tau_k} \mathbf{1}_{\{|V_N(t) - V_N(\tau
_{k-1})| > \eps^2\}} \, dt \Big| {\cal F}_{t_{k-1}} \biggr] \leq C
\theta\eps^2 + o(1).
\]
\end{Lemma}

\begin{pf}
Since $\theta\leq\theta^{1/4}$, and since $\delta< \eps^{5/2}$ by
(\ref{deltasmall}), it follows from (\ref{655}) that $\theta e^A \eps
^{-1/2} \leq\eps^2$. Therefore, by Proposition~\ref{EMkProp} and
Markov's inequality, on $G_{N, k-1}$,
\[
P (R_k > 0|{\cal F}_{t_{k-1}} ) \leq C \theta e^A \eps^{-1/2} + o(1)
\leq
C \eps^2 + o(1).
\]
$(Z_{N,1}(t_k) - Z_N(t_{k-1}))/(N(\log N)^2) = V_N(\tau_k) - V_N(\tau
_{k-1})$ on $G_{N, k-1} \cap\{R_k = 0\}$ and $4e^{-A/4} \leq\eps^2$
by (\ref{deltasmall}) and (\ref{652}), it follows from Corollary \ref
{driftcor} that
%
%e155 ###
%
%e156 #&#
\begin{eqnarray}\label{new45}\qquad
P\bigl(|V_N(\tau_k) - V_N(\tau_{k-1})| > \eps^2|{\cal F}_{t_{k-1}}\bigr) &\leq& C
\eps^2 + C \theta e^{-A/2} \eps^{-1/2} + o(1)\nonumber
\\[-8pt]
\\[-8pt] &\leq& C \eps^2 + o(1).
\nonumber
\end{eqnarray}

We claim that (\ref{new45}) also holds with $\tau_k$ replaced by any
$t$ such that $\tau_{k-1} < t < \tau_k$. Applying Corollary \ref
{driftcor} requires specifying five parameters: $u$, $s$, $\eps$, $A$
and $\theta$. To establish the claim, we apply Corollary~\ref{driftcor}
with\vspace*{1pt} new parameters ${\tilde u} = t_{k-1}/(\log N)^3$, ${\tilde s} =
s$, ${\tilde\eps} = \eps$, ${\tilde A} = A$ and ${\tilde\theta} = (t
- \tau_{k-1})/s$. Note that ${\tilde\theta} \leq\theta$, so
conditions (\ref{con0})--(\ref{con3}) continue to hold with the new
parameters. Also, using the new parameters, we get ${\tilde t}_0 =
{\tilde u}(\log N)^3 = t_{k-1}$ and ${\tilde t}_1 = ({\tilde u} +
{\tilde\theta} s)(\log N)^3 = t (\log N)^3$. It thus follows from
Corollary~\ref{driftcor} that
\[
P\bigl(|V_N(t) - V_N(\tau_{k-1})| > \eps^2|{\cal F}_{t_{k-1}}\bigr) \leq C \eps^2
+ C {\tilde\theta} e^{-{\tilde A}/2} {\tilde\eps}^{-1/2} + o(1)
\leq
C \eps^2 + o(1).
\]
Here the constant $C$ does not depend upon the choice of $t$. The
absolute value of the $o(1)$ can be bounded above by $B_N(t)$, where
$B_N(t) \leq1$ for all $N$ and $t$, and $\lim_{N \rightarrow\infty}
B_N(t) = 0$ for every fixed $t$. Thus, by Fubini's theorem and the
dominated convergence theorem,
\[
E \biggl[ \int_{\tau_{k-1}}^{\tau_k} \mathbf{1}_{\{|V_N(t) - V_N(\tau
_{k-1})| > \eps^2\}} \, dt \Big| {\cal F}_{t_{k-1}} \biggr] \leq\int
_{\tau_{k-1}}^{\tau_k} C \eps^2 + B_N(t) \, dt \leq C \theta\eps^2
+ o(1)
\]
as claimed.
\end{pf}

%
%le47 #&#
\begin{Lemma}\label{lem43}
We have
\begin{eqnarray*}
&&E \biggl[ \int_{\tau_{k-1}}^{\tau_k} A_1 f(V_N(t)) \, dt \Big| {\cal
F}_{t_{k-1}} \biggr] \mathbf{1}_{G_{N, k-1}} \\
&& \qquad= f'(V_N(\tau_{k-1})) V_N(\tau_{k-1}) \theta s (c + 2 \pi
^2 \log\eps)
\mathbf{1}_{G_{N, k-1}} + O(\theta\eps^{1/2}) + o(1).
\end{eqnarray*}
\end{Lemma}

\begin{pf}
If $0 \leq y \leq\eps$, then
\[
f\bigl(V_N(t) + y\bigr) = f(V_N(t)) + y f'(V_N(t)) + \tfrac{1}{2} f''(\xi_y) y^2
\]
for some $\xi_y$ satisfying $V_N(t) \leq\xi_y \leq V_N(t) + \eps$.
Therefore,
%e156 ###
%
%e157 #&#
\begin{eqnarray}\label{sjb}
&&\biggl| \int_{\tau_{k-1}}^{\tau_k} V_N(t) \int_0^{\eps} \bigl( f\bigl(V_N(t) +
y\bigr) - f(V_N(t)) - y f'(V_N(t)) \bigr) y^{-2} \, dy \, dt \biggr|
\nonumber\\
&& \qquad= \biggl| \int_{\tau_{k-1}}^{\tau_k} V_N(t) \biggl( \int_0^{\eps}
\frac
{1}{2} f''(\xi_y) \, dy \biggr) \, dt \biggr| \nonumber
\\[-8pt]
\\[-8pt]
&& \qquad\leq\theta s \sup_{t \in[\tau_{k-1}, \tau_k]} \sup_{z
\in[V_N(t),
V_N(t) + \eps]} \frac{\eps}{2} V_N(t) |f''(z)| \nonumber\\
&& \qquad\leq C \eps\theta s,
\nonumber
\end{eqnarray}
where the last inequality follows from the fact that $\|h\| < \infty$,
where $h$ is defined in (\ref{xfxy}). Equations (\ref{A1def}) and
(\ref
{sjb}) give
%
%e157 ###
%
%e158 #&#
\begin{eqnarray}\label{431}
&&  E \biggl[ \int_{\tau_{k-1}}^{\tau_k} A_1 f(V_N(t)) \, dt \Big| {\cal
F}_{t_{k-1}} \biggr] \mathbf{1}_{G_{N, k-1}}\nonumber\hspace*{-35pt}
\\[-8pt]
\\[-8pt]
&&   \qquad= (c + 2 \pi^2 \log\eps) E \biggl[ \int_{\tau_{k-1}}^{\tau_k} V_N(t)
f'(V_N(t)) \, dt\Big | {\cal F}_{t_{k-1}} \biggr] \mathbf{1}_{G_{N, k-1}}
+ O(\theta\eps).
\nonumber\hspace*{-35pt}
\end{eqnarray}
Recall that $d(x) = x f'(x)$ for $x \geq0$. Therefore,
%
%e158 ###
%
%e159 #&#
\begin{eqnarray}\label{432}
&&E \biggl[ \int_{\tau_{k-1}}^{\tau_k} V_N(t) f'(V_N(t)) \, dt \Big|
{\cal F}_{t_{k-1}} \biggr] \mathbf{1}_{G_{N, k-1}} \nonumber
\\
&& \qquad= f'(V_N(\tau_{k-1})) V_N(\tau_{k-1}) \theta s \mathbf
{1}_{G_{N, k-1}}\\
&& \qquad  \quad {} +
E \biggl[ \int_{\tau_{k-1}}^{\tau_k} d(V_N(t)) - d(V_N(\tau_{k-1})) \,
dt \Big| {\cal F}_{t_{k-1}} \biggr] \mathbf{1}_{G_{N, k-1}}.
\nonumber
\end{eqnarray}
The absolute value of the second term on the right-hand side of (\ref
{432}) is at most
\[
2 \|d\| E \biggl[\int_{\tau_{k-1}}^{\tau_k} \mathbf{1}_{\{|V_N(t) -
V_N(\tau
_{k-1})| > \eps^2\}} \, dt \Big| {\cal F}_{t_{k-1}} \biggr] \mathbf{1}_{G_{N,
k-1}} + \theta s\|d'\| \eps^2,
\]
which is at most $C \theta\eps^2 + o(1)$ by Lemma~\ref{uniformlem}.

Therefore, the result follows from (\ref{431}) and (\ref{432}), since
$\eps|\log\eps| < \eps^{1/2}$ for sufficiently small $\eps$.
\end{pf}

%
%le48 #&#
\begin{Lemma}\label{lem44}
We have
\begin{eqnarray*}
&&E \biggl[ \int_{\tau_{k-1}}^{\tau_k} A_2 f(V_N(t)) \, dt\Big | {\cal
F}_{t_{k-1}} \biggr] \mathbf{1}_{G_{N, k-1}} \\
&& \qquad= 2 \pi^2 \theta s V_N(\tau_{k-1}) \mathbf{1}_{G_{N, k-1}}
\int_{\eps
}^{\infty} \bigl( f\bigl(V_N(\tau_{k-1}) + y\bigr) - f(V_N(\tau_{k-1})) \bigr)
y^{-2} \, dy\\
&& \qquad  \quad {} + O(\theta\eps)+ o(1).
\end{eqnarray*}
\end{Lemma}

\begin{pf}
For $y \geq0$, define the function $g_y(x) = x f(x + y)$. Note that
$\sup_{y \geq0} \|g_y\| < \infty$ and $\sup_{y \geq0} \|g_y'\| <
\infty$. We have
%
%e159 ###
%
%e160 #&#
\begin{eqnarray}\label{44main}
&&E \biggl[ \int_{\tau_{k-1}}^{\tau_k} A_2 f(V_N(t)) \, dt\Big | {\cal
F}_{t_{k-1}} \biggr] \mathbf{1}_{G_{N, k-1}} \nonumber\hspace*{-35pt}\\
&& \qquad= E \biggl[ 2 \pi^2 \int_{\tau_{k-1}}^{\tau_k} V_N(t) \int
_{\eps
}^{\infty}\bigl ( f\bigl(V_N(t) + y\bigr) - f(V_N(t)) \bigr) y^{-2} \, dy \, dt
\Big| {\cal F}_{t_{k-1}} \biggr]\nonumber\hspace*{-35pt}\\
&& \qquad  \quad {}\times \mathbf{1}_{G_{N, k-1}}\nonumber\hspace*{-35pt}
\\[-8pt]
\\[-8pt]
&& \qquad= 2 \pi^2 \theta s V_N(\tau_{k-1}) \mathbf{1}_{G_{N, k-1}}
\int_{\eps
}^{\infty} \bigl( f\bigl(V_N(\tau_{k-1}) + y\bigr) - f(V_N(\tau_{k-1})) \bigr)
y^{-2} \, dy \nonumber\hspace*{-35pt}\\
&& \qquad \quad{} + 2 \pi^2 E\biggl [ \int_{\tau_{k-1}}^{\tau_k} \int
_{\eps
}^{\infty} \bigl( g_y(V_N(t)) - g_0(V_N(t))- g_y(V_N(\tau_{k-1})) \nonumber\hspace*{-35pt}\\
&&\hspace*{47pt}\hphantom{+ 2 \pi^2 E\biggl [ \int_{\tau_{k-1}}^{\tau_k} \int
_{\eps
}^{\infty} \bigl(} \qquad \quad{}  + g_0(V_N(\tau_{k-1})) \bigr)
y^{-2} \, dy \, dt \Big| {\cal F}_{t_{k-1}} \biggr] \mathbf{1}_{G_{N,
k-1}}.
\nonumber\hspace*{-35pt}
\end{eqnarray}
The absolute value of the second term on the right-hand side of (\ref
{44main}) is at most
\begin{eqnarray*}
&&\frac{2 \pi^2 s}{\eps} \Bigl(2 \sup_{y \geq0} \|g_y\| + 2 \|g_0\| \Bigr)
E \biggl[ \int_{\tau_{k-1}}^{\tau_k} \mathbf{1}_{\{|V_N(t) - V_N(\tau
_{k-1})| > \eps^2\}}\,dt \Big| {\cal F}_{t_{k-1}} \biggr] \mathbf{1}_{G_{N,
k-1}} \\
&& \qquad{}+ 2 \pi^2 \theta s \eps\Bigl(\sup_{y \geq0} \|g_y'\| + \|
g_0'\|
\Bigr),
\end{eqnarray*}
using that $\int_{\eps}^{\infty} y^{-2} \, dy = \eps^{-1}$. By Lemma
\ref{uniformlem}, this expression is at most $C \theta\eps+ o(1)$,
which, combined with (\ref{44main}), implies the result.
\end{pf}

\begin{pf*}{Proof of Proposition~\ref{ZNCSBP}}
Recall that we need to establish (\ref{mainVN}). For $1 \leq k \leq
\theta^{-1}$, define
\[
J_k = f(V_N(\tau_k)) - f(V_N(\tau_{k-1})) - \int_{\tau_{k-1}}^{\tau_k}
A f(V_N(t)) \, dt.
\]
Then
%
%e160 ###
%
%e161 #&#
\begin{equation}\label{sumJk}
f\bigl(V_N(u + s)\bigr) - f(V_N(s)) - \int_u^{u+s} A f(V_N(t)) \, dt = \sum
_{k=1}^{\theta^{-1}} J_k.
\end{equation}
Let $B_{N,0} = G_{N,0}^c$, and for $1 \leq k \leq\theta^{-1}$, let
$B_{N,k} = G_{N, k-1} \cap G_{N,k}^c$. Then $G_N(\eps)^c = \bigcup_{k=0}^{\theta^{-1}} B_{N,k}$ and
\[
1 - P(G_N(\eps)) = \sum_{k=0}^{\theta^{-1}} P(B_{N,k}).
\]
Now
%
%e161 ###
%
%e162 #&#
\begin{eqnarray}\label{Jkeq}
&&E \Biggl[ \Biggl( \sum_{k=1}^{\theta^{-1}} J_k \Biggr) \prod_{i=1}^j
h_i(V_N(s_i)) \Biggr]\nonumber\\
 && \qquad = E \Biggl[ \Biggl( \sum_{k=1}^{\theta^{-1}} J_k
\Biggl( \mathbf{1}_{G_{N, k-1}} + \sum_{\ell= 0}^{k-1} \mathbf{1}_{B_{N,
\ell
}} \Biggr) \Biggr) \prod_{i=1}^j h_i(V_N(s_i)) \Biggr] \nonumber
\\[-8pt]
\\[-8pt]
&& \qquad = E \Biggl[ \Biggl( \sum_{k=1}^{\theta^{-1}} J_k \mathbf{1}_{G_{N, k-1}}
\Biggr) \prod_{i=1}^j h_i(V_N(s_i)) \Biggr] \nonumber\\
&&{} \qquad  \quad + \sum_{\ell= 0}^{\theta
^{-1} - 1} E \Biggl[\Biggl ( \sum_{k = \ell+ 1}^{\theta^{-1}} J_k \Biggr)
\mathbf{1}_{B_{N, \ell}} \prod_{i=1}^j h_i(V_N(s_i)) \Biggr].
\nonumber
\end{eqnarray}
For $0 \leq\ell\leq\theta^{-1} - 1$,
\begin{eqnarray*}
\Biggl| \sum_{k = \ell+ 1}^{\theta^{-1}} J_k \Biggr| &=& \biggl| f\bigl(V_N(u +
s)\bigr) - f(V_N(\tau_{\ell})) - \int_{\tau_{\ell}}^{u + s} A f(V_N(t))
\,
dt \biggr|\\
 &\leq&2 \|f\| + s \|Af\|.
\end{eqnarray*}
Therefore, the absolute value of the second term on the right-hand side
of (\ref{Jkeq}) is at most
\[
\Biggl( \prod_{i=1}^j \|h_i\| \Biggr) (2 \|f\| + s \|Af\| ) \sum
_{\ell= 0}^{\theta^{-1} - 1} P(B_{N,\ell}) \leq C \bigl(1 - P(G_N(\eps))\bigr),
\]
using (\ref{Afbound}). To bound the first term on the right-hand side
of (\ref{Jkeq}), note that by conditioning on ${\cal F}_{t_{k-1}}$,
\begin{eqnarray*}
&&E \Biggl[ \Biggl( \sum_{k=1}^{\theta^{-1}} J_k \mathbf{1}_{G_{N, k-1}}
\Biggr) \prod_{i=1}^j h_i(V_N(s_i)) \Biggr]\\
 && \qquad = \sum_{k=1}^{\theta^{-1}} E
\Biggl[ J_k \mathbf{1}_{G_{N, k-1}} \prod_{i=1}^j h_i(V_N(s_i)) \Biggr]
\\
&& \qquad = \sum_{k=1}^{\theta^{-1}} E \Biggl[ \Biggl( \prod_{i=1}^j h_i(V_N(s_i))
\Biggr) E[J_k|{\cal F}_{t_{k-1}}] \mathbf{1}_{G_{N, k-1}} \Biggr].
\end{eqnarray*}
By Lemmas~\ref{lem41},~\ref{lem42},~\ref{lem43} and~\ref{lem44},
\[
| E[J_k|{\cal F}_{t_{k-1}}] \mathbf{1}_{G_{N, k-1}} | \leq C
\theta\eps^{1/2} + o(1)
\]
for all $k$. Therefore,
\begin{eqnarray*}
&&\Biggl| E \Biggl[ \Biggl( \sum_{k=1}^{\theta^{-1}} J_k \mathbf{1}_{G_{N,
k-1}} \Biggr) \prod_{i=1}^j h_i(V_N(s_i)) \Biggr] \Biggr|\\
 && \qquad \leq \Biggl(
\prod_{i=1}^j \|h_i\| \Biggr)\Biggl ( \sum_{k=1}^{\theta^{-1}} E [
|E[J_k|{\cal F}_{t_{k-1}}] \mathbf{1}_{G_{N, k-1}}| ] \Biggr)
\\
&& \qquad \leq C \eps^{1/2} + o(1).
\end{eqnarray*}
It follows that
\[
\Biggl| E \Biggl[ \Biggl( \sum_{k=1}^{\theta^{-1}} J_k \Biggr) \prod
_{i=1}^j h_i(V_N(s_i)) \Biggr] \Biggr| \leq C \eps^{1/2} + C \bigl(1 -
P(G_N(\eps))\bigr) + o(1).
\]
In view of (\ref{sumJk}) and Proposition~\ref{GNProp}, equation (\ref
{mainVN}) now follows by letting $N \rightarrow\infty$ and then
letting $\eps\rightarrow0$.
\end{pf*}

%s6.3 ###
%s6.3 #&#
\subsection{The number of particles}\label{sec63}

Because the value of $Z_N(t)$ approximately determines the number of
particles a short time after time $t$, the fact that the number of
particles converges to a continuous-state branching process follows
rather simply from Proposition~\ref{ZNCSBP}.

\begin{pf*}{Proof of Theorem~\ref{MNCSBP}}
In view of Proposition~\ref{ZNCSBP}, it suffices to show that for any
fixed $t > 0$, we have
%
%e162 ###
%
%e163 #&#
\begin{equation}\label{mvfin}
\biggl| \frac{1}{2 \pi N} M_N((\log N)^3 t) - V_N(t) \biggr| \rightarrow
_p 0.
\end{equation}
Let $\gamma> 0$ be arbitrary. Set $u = 0$ and $s = t$. By Proposition
\ref{GNProp}, we can choose $\eps\in(0, \gamma)$ sufficiently small that
%
%e163 ###
%
%e164 #&#
\begin{equation}\label{new21}
\sup_{\theta}\Bigl ( \limsup_{N \rightarrow\infty} \bigl( 1 - P(G_N(\eps
)) \bigr) \Bigr) < \frac{\gamma}{2},
\end{equation}
where the supremum is taken over all $\theta$ such that $\theta^{-1}
\in\N$. Proposition~\ref{rjumplem} implies that for sufficiently small
$\theta$,
%
%e164 ###
%
%e165 #&#
\begin{equation}\label{new22}
  P\bigl(|V_N(t) - V_N\bigl(t(1 - \theta)\bigr)| \geq\gamma\bigr) \leq C \theta\eps^{-3/2}
+  \bigl(1 - P(G_N(\eps))\bigr) + o(1).\hspace*{-35pt}
\end{equation}
It follows from (\ref{new21}) and (\ref{new22}) that for sufficiently
small $\theta$ and sufficiently large~$N$,
%
%e165 ###
%
%e166 #&#
\begin{equation}\label{mv1}
P\bigl(|V_N(t) - V_N\bigl(t(1 - \theta)\bigr)| < \gamma\bigr) > 1 - \gamma.
\end{equation}

Let $M'_N((\log N)^3 t)$ denote the number of particles at time $(\log
N)^3 t$ whose ancestor at time $u$ is in $(0, L)$ for all $(\log N)^3
(t(1 - \theta)) \leq u \leq(\log N)^3 t$. By Proposition~\ref{EMkProp}
and (\ref{new21}), for sufficiently small $\theta> 0$ and sufficiently
large~$N$,
%
%e166 ###
%
%e167 #&#
\begin{equation}\label{mv2}
P \bigl( M_N((\log N)^3 t) = M'_N((\log N)^3 t) \bigr) > 1 - \gamma.
\end{equation}
By (\ref{Mexp}) and the fact that $1 - \mu^2/2 - \pi^2/2L^2 = 0$,
\begin{eqnarray*}
&&E\bigl[M'_N((\log N)^3 t)|{\cal F}_{(\log N)^3(t(1 - \theta))}\bigr] \\
&& \qquad = \frac{2 N
(\log N)^2 V_N(t(1 - \theta))(1+o(1))}{L} \int_0^L e^{-\mu y} \sin
\biggl( \frac{\pi y}{L} \biggr) \, dy.
\end{eqnarray*}
Now
\begin{eqnarray*}
\int_0^L e^{-\mu y} \sin\biggl( \frac{\pi y}{L} \biggr) \, dy &=& \int
_0^{\infty} \frac{\pi y}{L} e^{-\mu y} \, dy + \int_0^L e^{-\mu y}
\biggl( \sin\biggl( \frac{\pi y}{L} \biggr) - \frac{\pi y}{L} \biggr) \, dy\\
&&{} -
\int_L^{\infty} \frac{\pi y}{L} e^{-\mu y} \, dy \\
&=& \frac{\pi}{L \mu^2} + O \biggl( \int_0^L e^{-\mu y} \frac{y^3}{L^3}
\, dy \biggr) + O( e^{-\mu L}) \\
&=& \frac{\pi}{L \mu^2} + O \biggl( \frac{1}{L^3} \biggr) + O(e^{-\mu L})
= \frac{\pi}{2L} \bigl(1 + o(1)\bigr).
\end{eqnarray*}
It follows that
\begin{eqnarray*}
E\bigl[M'_N((\log N)^3 t)|{\cal F}_{(\log N)^3 (t(1 - \theta))}\bigr] &=& \frac
{\pi N (\log N)^2 V_N(t(1 - \theta))(1+o(1))}{L^2} \\
&=& 2 \pi N V_N\bigl(t(1 - \theta)\bigr) \bigl(1 + o(1)\bigr).
\end{eqnarray*}
Therefore, for sufficiently large $N$,
%
%e167 ###
%
%e168 #&#
\begin{equation}\label{mv3}
 P \biggl( \biggl| E \biggl[ \frac{M'_N((\log N)^3 t)}{2 \pi N} \Big|
{\cal F}_{(\log N)^3 (t(1 - \theta))} \biggr] - V_N\bigl(t(1 - \theta)\bigr)
\biggr| < \gamma\biggr) > 1 - \gamma.\hspace*{-35pt}
\end{equation}
By Proposition~\ref{MNVar}, we have
\[
\operatorname{Var}\bigl(M'_N((\log N)^3 t)|{\cal F}_{(\log N)^3 (t(1 -
\theta))}\bigr)
\leq C \theta\eps^{-1/2} N^2\bigl(1 + o(1)\bigr)
\]
on an event defined in the same manner as $G_{N, k-1}$ but with $(\log
N)^3(t(1 - \theta))$ playing the role of $t_{k-1}$. Combining this
result with (\ref{new21}) and the conditional form of Chebyshev's
inequality, we get for sufficiently small $\theta$ and sufficiently
large $N$,
%
%e168 ###
%
%e169 #&#
\begin{eqnarray}\label{mv4}\quad
 &&  P \biggl( \biggl| \frac{M'_N((\log N)^3 t)}{2 \pi N} - E \biggl[ \frac
{M'_N((\log N)^3 t)}{2 \pi N} \Big| {\cal F}_{(\log N)^3 (t(1 - \theta
))} \biggr] \biggr| < \gamma\biggr)\nonumber
\\[-8pt]
\\[-8pt]
&& \qquad  > 1 - \gamma. \nonumber
\end{eqnarray}
If now follows from (\ref{mv1})--(\ref
{mv4}) that for sufficiently large $N$, we have
\[
P \biggl( \biggl| \frac{1}{2 \pi N} M_N((\log N)^3) - V_N(t) \biggr| < 3
\gamma\biggr) > 1 - 4 \gamma.
\]
Result (\ref{mvfin}) follows.
\end{pf*}

%s7 ###
%s7 #&#
\section{Convergence to the Bolthausen--Sznitman coalescent}\label{sec7}

In this section, we prove Theorem~\ref{bosz}. The strategy will be to
show that a sequence of processes that describe the genealogy of
branching Brownian motion converges to a flow of bridges defined in
\cite{beleg1}, which is known to be dual to the Bolthausen--Sznitman
coalescent.

%s7.1 ###
%s7.1 #&#
\subsection{The flow of bridges}\label{sec71}

Consider the continuous state branching process of Proposition \ref
{ZNCSBP} and Theorem~\ref{MNCSBP} with branching mechanism $\Psi(u)= au
+ 2\pi^2 u \log u$. Recall from~\cite{beleg00} that we can define this
as a two-parameter process $(Z(t,x), t \geq0,  x \geq0)$, where $t$ is
the time parameter, and $x$ is the initial population size. Also recall
from~\cite{beleg00} that we can associate with this continuous-state
branching process a flow of subordinators. On some probability space,
there exists a process
$(S^{(s,t)}(x), 0 \leq s \leq t, x \geq0)$ such that:
\begin{itemize}
\item For every $0 \le s \le t$, the process $S^{(s,t)}=(S^{(s,t)}(x),
x\geq0)$ is a subordinator with Laplace exponent $u_{t-s}$.
\item For every integer $k \ge2$ and every $0 \leq t_1 \leq\cdots
\leq t_k$, the subordinators $S^{(t_1,t_2)}, \ldots, S^{(t_{k-1},t_k)}$
are independent, and
\[
S^{(t_1,t_k)} = S^{(t_{k-1},t_k)}\circ\cdots\circ S^{(t_1,t_2)}.
\]
\item The processes $(Z(t,x), t \geq0, x \geq0)$ and $(S^{(0,t)}(x),
t \geq0, x \geq0)$ have the same finite-dimensional marginals.
\end{itemize}
Here $S^{(s,t)}(x)$ can be understood as the descendants in the
population at time $t$ of the first $x$ individuals in the population
at time $s$.

Suppose that we start with the initial population $Z(0)=z.$ For each $s
\leq t$, we can define the renormalized process $(B_{s,t}(x), 0 \leq x
\leq1)$ by
\[
B_{s,t}(x) = S^{(s,t)} \bigl( x S^{(0,s)} (z) \bigr)/ S^{(0,t)}(z).
\]
It is easily seen that $B_{s,t}$ is a \textit{bridge}, which we define
as in~\cite{beleg1} to be a nondecreasing, $[0, 1]$-valued stochastic
process $(B(r), 0 \leq r \leq1)$ with exchangeable increments and
right-continuous paths such that $B(0)=0$ and $B(1)=1$.

It follows from (\ref{utlambda}) that when $\Psi(u)= au + 2\pi^2 u
\log
u$, the subordinator $S^{(s,t)}$ is a stable subordinator with index
$e^{-2 \pi(t-s)}$. Consequently, letting $R_{s,t}$ denote the range of
$B_{s,t}$, the lengths of the disjoint open intervals whose union is
$[0, 1] \setminus R_{s,t}$ are independent of $S^{(0,s)}(z)$ and have
the Poisson--Dirichlet distribution with parameters $(e^{-2 \pi(t-s)},
0)$. See~\cite{py97} for a definition and further discussion of the
two-parameter Poisson--Dirichlet distribution and its connections with
stable subordinators. It now follows (see Example 2 in~\cite{beleg1})
that $(B_{s,t}(x), 0 \leq s \leq t, 0 \leq x \leq1)$ is a \textit{flow
of bridges}, which is a collection $(B_{s,t}, 0 \leq s \leq t)$ of
bridges such that if $\operatorname{Id}$ denotes the identity function from
$[0,1]$ to itself, then:
\begin{itemize}
\item For every $s < t< u$, we have $B_{s,u} = B_{t,u} \circ B_{s,t}$.
\item The law of $B_{s,t}$ only depends on $t-s$.\vadjust{\goodbreak}
\item If $s_1 < s_2 <\cdots< s_n$, then the bridges $B_{s_1,s_2},
\ldots, B_{s_{n-1}, s_n}$ are independent.
\item$B_{0,0}=\operatorname{Id}$ and $B_{0,t} \to\operatorname{Id}$
as $t \to0$ in
probability, in the sense of Skorohod's topology.
\end{itemize}
Note that we are using a different convention for the time parameters
than in~\cite{beleg1}. The bridge $B_{s,t}$ defined here would be
called $B_{-t, -s}$ in~\cite{beleg1}.

If $B$ is a bridge, define, for $u \in[0, 1]$,
%
%e169 ###
%
%e170 #&#
\begin{equation}\label{bridgeinv}
B^{-1}(u) = \inf\{s \in[0,1]\dvtx  B(s) \geq u\}.
\end{equation}
If $s < t < u$, then $B_{s,u}^{-1} = B_{s,t}^{-1} \circ B_{t,u}^{-1}$.
Given independent random variables $U_1,\ldots, U_n$ with the uniform
distribution on $[0,1]$, we can define $\pi(B)$ to be the partition of
$\{1,\ldots, n\}$ such that $i$ and $j$ are in the same block of $\pi
(B)$ if and only if $B^{-1}(U_i) = B^{-1}(U_j)$. Now, given a flow of
bridges $(B_{s,t}, 0 \leq s \leq t)$ and independent uniform random
variables $U_1,\ldots, U_n$, we can fix a time $t > 0$ and consider the
partition-valued process $(\Pi(s), 0 \leq s \leq t)$ defined by $\Pi(s)
= \pi(B_{t-s, t})$. The main result of Bertoin and Le Gall \cite
{beleg1} establishes that this process is a so-called exchangeable
coalescent process and that there is in fact a one-to-one
correspondence between flows of bridges and exchangeable coalescent
processes. In the example above, in which the flow of bridges is
defined from a continuous-state branching process with $\Psi(u) = au +
2 \pi^2 u \log u$, the process $(\pi(B_{t - s/2\pi, t}), 0 \leq s
\leq
2 \pi t)$ is the Bolthausen--Sznitman coalescent run for time
$2 \pi t$ (see, e.g., Example 2 in~\cite{beleg1}).

%s7.2 ###
%s7.2 #&#
\subsection{Flows describing the genealogy of branching Brownian motion}\label{sec72}

To represent the genealogy of branching Brownian motion, we now
introduce a sequence of discrete versions of these flows of bridges. We
fix $K \in\N$ and the times $0 = t_0 < t_1 <\cdots< t_K.$ For $0 \le
i < j \le K$ we will define a process $(B_{t_i,t_j}^N(s), 0 \le s \le1)$.

We consider the branching Brownian motion $X_N$ at the successive times
$t_j (\log N)^3$. We assign labels to the particles at these times, and
denote by $u_{i,j}$ the label of the $i$th largest particle at time
$t_j(\log N)^3$, that is, the particle in position $X_{i, N}(t_j(\log
N)^3)$. We first define a collection of independent random variables
$(v_{i,j}, i \geq0, 0 \leq j \leq K)$ having the uniform distribution
on $[0,1]$. For $i \leq M_N(0)$, we define $u_{i,0} = v_{i,0}$. That
is, the individuals at time zero are labeled by independent uniform
random variables. For $j \geq1$, the $u_{i,j}$ are sequences of length
$j+1$ which are defined inductively by saying that $u_{i,j} = (u_{p(i),
j-1}, v_{i,j})$, where $u_{p(i), j-1}$ is the label of the particle at
time $t_{j-1}(\log N)^3$ from which the $i$th particle at time
$t_j(\log N)^3$ has descended. That is, we concatenate $v_{i,j}$ with
the label of the ancestor of the $i$th particle to obtain the label of
the $i$th particle. The particles at time $t_j(\log N)^3$ can now be
ordered using the lexicographical order of their labels. We denote by
$x_{i,j}$ the position of the $i$th individual in this lexicographical
order at time $t_j(\log N)^3$.

We now assign weights to the individuals. For $0 \leq j \leq K$ and $1
\leq i \leq M_N(t_j(\log N)^3)$, define
\[
w(i,j) =
\cases{\displaystyle \displaystyle\frac{1}{Z_N(t_j(\log N)^3)}
e^{\mu
x_{i,j}} \sin\biggl( \frac{\pi x_{i,j}}{L} \biggr) \mathbf{1}_{\{x_{i,j}
\leq L\}} ,&\quad if $0 \leq j \leq K-1 $,\vspace*{3pt}\cr\displaystyle
{\displaystyle\frac{1}{M_N(t_j(\log N)^3)}} ,&\quad if $j = K$.
}
\]
That is, the particles are weighted proportional to their contribution
to the sum in (\ref{ZNdef}), except at time $t_K(\log N)^3$ when
all particles are weighted equally. We use these weights because we
will later sample particles uniformly at time $t_K(\log N)^3$, but the
number of descendants that a particle at time $t_i(\log N)^3$ has at
time $t_K(\log N)^3$ will be roughly proportional to the weight that it
has been assigned. Also define $A_i(j,k)$ to be the set of descendants
at time $t_k(\log N)^3$ of the $i$th individual at time $t_j(\log
N)^3$. More precisely, $A_i(j,k)$ is the set of indices $\ell$ such
that the individual at position $x_{\ell, k}$ at time $t_k(\log N)^3$
is descended from the individual in position $x_{i,j}$ at time
$t_j(\log N)^3$. We are now ready to define the discrete bridges.
First, for $0 \leq y \leq1$, and $0 \leq j \leq K$, define
%
%e170 ###
%
%e171 #&#
\begin{equation}\label{Lidef}
L_j(y) = \max\Biggl\{ I \in\N\dvtx  \sum_{i=1}^{I} w(i,j) \le y \Biggr\},
\end{equation}
with the convention that the maximum of the empty set is $0$. We think
of $L_j(y)$ as being approximately the $y$th quantile of the population
at time $t_j(\log N)^3$ when individuals are weighted as above and
ordered according to their labels. Then for $0 \le y \le1$ and $0 \leq
j < k \leq K$, let
\[
B_{t_j, t_k}^N(y) = \sum_{i=1}^{L_j(y)} \sum_{m \in A_i(j,k)} w(m,k).
\]
Note that these discrete bridges $B_{t_j, t_k}^N$ are not exactly
bridges in the sense defined above; for example, their increments are
not exactly exchangeable because there are only finitely many particles
at time $t_j$. However, we will show in Lemmas~\ref{L1ptcv} and \ref
{Lmulticv} below that these discrete bridges converge to the bridges
$B_{t_j, t_k}$.

%
%le49 #&#
\begin{Lemma} \label{cocycle}
If $0 \leq i<j<k \leq K$, then $B_{t_i,t_k}^N=B_{t_j,t_k}^N \circ
B^N_{t_i,t_j}$ and\break $(B_{t_i, t_k}^N)^{-1} = (B_{t_i, t_j}^N)^{-1} \circ
(B_{t_j, t_k}^N)^{-1}$, where the inverse functions are defined as in
(\ref{bridgeinv}).
\end{Lemma}

\begin{pf}
For $0 \le y \le1$,
\[
B_{t_i, t_k}^N(y) = \sum_{\ell=1}^{L_i(y)} \sum_{m \in A_{\ell
}(i,k)} w(m,k).
\]
Note that $m \in A_{\ell}(i,k)$ for some ${\ell} \le L_i(y)$ if and
only if $m \in A_{\ell}(j,k)$ for some ${\ell} \le
L_j(B_{t_i,t_j}^N(y)).$ Therefore,
\[
B_{t_i, t_k}^N(y) = \sum_{\ell=1}^{L_j(B_{t_i,t_j}^N(y))} \sum_{m
\in
A_{\ell}(j,k)} w(m,k) = B_{t_j, t_k}^N(B_{t_i, t_j}^N(y)).
\]
That is, $B_{t_i,t_k}^N=B_{t_j,t_k}^N \circ B^N_{t_i,t_j}$. Also,
\begin{eqnarray*}
(B_{t_i, t_k}^N)^{-1}(y) &=& \inf\{s\dvtx  B_{t_i, t_k}^N(s) \geq y\}
\\
&=& \inf\{s\dvtx  B_{t_j, t_k}^N(B_{t_i, t_j}^N(s)) \geq y\} \\
&=& \inf\{s\dvtx  B_{t_i, t_j}^N(s) \geq(B_{t_j, t_k}^N)^{-1}(y)\}
\\
&=& (B_{t_i, t_j}^N)^{-1}((B_{t_j, t_k}^N)^{-1}(y)),
\end{eqnarray*}
which implies that $(B_{t_i, t_k}^N)^{-1} = (B_{t_i, t_j}^N)^{-1} \circ
(B_{t_j, t_k}^N)^{-1}$.
\end{pf}

%s7.3 ###
%s7.3 #&#
\subsection{Convergence of one bridge}\label{sec73}

Let $(B_{s,t}, 0 \leq s \leq t)$ be the flow of bridges defined above
from the continuous-state branching process with branching mechanism
$\Psi(u) = au + 2 \pi^2 u \log u$. We will now show that for $1 \leq i
\leq K$, the sequence of discrete bridges $(B^N_{0, t_i}(u), 0\le u \le
1)$ converges to $(B_{0, t_i}(u), 0\le u \le1)$ in the sense of
finite-dimensional distributions. The first step is the following
extension of Proposition~\ref{ZNCSBP}.

%
%le50 #&#
\begin{Lemma}\label{ZNCSBPprime}
Assume that the initial population is subdivided into $m$ possibly
random subgroups $S_1, \ldots, S_m$, and that given the initial
positions of the particles, they evolve according to branching Brownian
motion killed at 0. Assume that $Y_N(0)/(N (\log N)^3)$ converges to
zero in probability.
Let $Z_{i,N}(t)$ denote the contribution to the sum in (\ref{ZNdef})
from particles descended from one of the particles that is in $S_i$ at
time zero, and let $M_{i,N}(t)$ denote the number of particles at time
$t$ descended from one of the particles that is in $S_i$ at time zero.
Assume that the initial joint distribution of
%
%e171 ###
%
%e172 #&#
\begin{equation}\label{Ejointin}
\biggl( \frac{Z_{i,N}(0)}{N(\log N)^2} \biggr)_{i=1}^m
\end{equation}
converges as $N \rightarrow\infty$ to some probability measure $\rho$
on $[0, \infty)^m$.
Then the finite-dimensional distributions of the $m$-dimensional
vector-valued processes
\[
\biggl\{\biggl (\frac{Z_{i,N}(t(\log N)^3)}{N(\log N)^2} \biggr)_{i=1}^m, t
\geq0 \biggr\} \quad\mbox{and} \quad \biggl\{ \biggl(\frac
{M_{i,N}(t(\log N)^3)}{2 \pi N} \biggr)_{i=1}^m, t > 0 \biggr\}
\]
each converge as $N \rightarrow\infty$ to the finite-dimensional
distributions of $ \{(Z_i(t))_{i=1}^m,\break t\ge0 \},$ where
$(Z_i(0))_{i=1}^m$ has distribution $\rho$, and conditional on
$(Z_i(0))_{i=1}^m$, each $Z_i$ evolves independently as a
continuous-state branching process with branching mechanism $\Psi(u) =
au + 2 \pi u \log u$.
\end{Lemma}

\begin{pf}
While this is in principle a simple extension of Proposition \ref
{ZNCSBP}, some care is needed in the proof because the components of
the process are not independent but only conditionally independent
given the initial configuration. To ease notation, we only show here
the proof of the one-dimensional marginal convergence (which is all
that is needed later), as the general result is conceptually identical
but more cumbersome. Thus, let $t > 0$, and fix arbitrary bounded and
continuous test functions $f_1, \ldots, f_m\dvtx  [0, \infty) \rightarrow
\R
$. By Skorohod's Representation Theorem, we may assume that all the
branching Brownian motions $X_N$ are constructed on the same
probability space in such a way that the expression in (\ref{Ejointin})
converges almost surely to $(Z_i(0))_{i=1}^m$ having joint distribution
$\rho$.

For $i = 1,\ldots, m$, let $X_{i,N}$ denote the branching Brownian
motion obtained by considering only the descendants of the particles in
$S_i$. Let ${\cal F} = \sigma(X_{i,N}(0), i=1, \ldots, m, N=1, 2,
\dots
)$ be the filtration generated by all the processes at time zero for
all subgroups. Let also ${\cal G} = \sigma(Z_1(0), \ldots, Z_m(0))$.
Note that the random variables $Z_{1,N}(t(\log N)^3),\ldots,
Z_{m,N}(t(\log N)^3)$ are conditionally independent given~${\cal F}$. Therefore,
%
%e172 ###
%
%e173 #&#
\begin{eqnarray}\label{Ejoint2}\qquad
E \Biggl[\prod_{i=1}^m f_i \biggl( \frac{Z_{i,N}(t(\log N)^3)}{N (\log
N)^2} \biggr) \Biggr] & =& E \Biggl[E \Biggl[ \prod_{i=1}^m f_i\biggl (
\frac{Z_{i,N}(t(\log N)^3)}{N (\log N)^2} \biggr) \Big|{\cal F}
\Biggr] \Biggr] \nonumber
\\[-8pt]
\\[-8pt]
& =& E \Biggl[\prod_{i=1}^m E \biggl[ f_i \biggl( \frac{Z_{i,N}(t(\log
N)^3)}{N (\log N)^2} \biggr) \Big|{\cal F} \biggr] \Biggr].
\nonumber
\end{eqnarray}
By Proposition~\ref{ZNCSBP}, for $1 \leq i \leq m$ we have that
$E[f_i(Z_{i,N}(t(\log N)^3)/\break(N(\log N)^2))| {\cal F}]$ converges almost surely
to the random variable $E_{Z_i(0)}[f_i(Z(t))]$, where $E_x[f_i(Z(t))]$
denotes the expected value for the continuous-state branching process started
from the value $Z(0) =x$.
The application of Proposition~\ref{ZNCSBP} is justified here because
the condition that $Y_N(0)/(N(\log N)^3)$ converges in probability to
zero is satisfied for the entire process, and hence the analogous
condition is satisfied for each of the $m$ components.
We may rewrite this random variable as
\[
E_{Z_i(0)}[f_i(Z(t))] = E[f_i(Z_i(t))|{\cal G}].
\]
Since all random variables on the right-hand side of (\ref{Ejoint2})
are bounded, we deduce by the dominated convergence theorem that
\begin{eqnarray*}
\lim_{N\to\infty} E \Biggl[\prod_{i=1}^m f_i \biggl( \frac{Z_{i,N}(t(\log
N)^3)}{N (\log N)^2} \biggr) \Biggr] &=& E \Biggl[\prod_{i=1}^m
E[f_i(Z_i(t))|{\cal G}] \Biggr]\\
 &=& E \Biggl[\prod_{i=1}^m f_i(Z_i(t)) \Biggr],
\end{eqnarray*}
since the random variables $(Z_i(t))_{i=1}^m$ are conditionally
independent given ${\cal G}$. This completes the proof of convergence
for the processes $Z_{i,N}$. The proof for the processes $M_{i,N}$ is
identical, except that we invoke Theorem~\ref{MNCSBP} instead of
Proposition~\ref{ZNCSBP}.
\end{pf}

Before proving the convergence of bridges, we establish the following
lemma, which states that at a typical time $t$, no single particle
makes too large a contribution to $Z_N(t)$.

%
%le51 #&#
\begin{Lemma}\label{maxcont}
Let
\[
m_N(s) = \max_{1 \leq i \leq M_N(s(\log N)^3)} e^{\mu X_{i,N}(s(\log
N)^3)} \sin\biggl( \frac{\pi X_{i,N}(s (\log N)^3)}{L} \biggr).
\]
Then for all $s \geq0$, we have $m_N(s)/(N(\log N)^2) \rightarrow0$
in probability as \mbox{$N \rightarrow\infty$}.
\end{Lemma}

\begin{pf}
Suppose $(x_N)_{N=1}^{\infty}$ is a sequence such that $e^{\mu x_N}/(N
(\log N)^3) \rightarrow0$ as $N \rightarrow\infty$. Letting $w_N = L
- x_N$, we have $w_N \rightarrow\infty$ as $N \rightarrow\infty$. Therefore,
%
%e173 ###
%
%e174 #&#
\begin{eqnarray}\label{maxxn}
e^{\mu x_N} \sin\biggl( \frac{\pi x_N}{L} \biggr) &=& e^{\mu(L - w_N)}
\sin\biggl( \frac{\pi w_N}{L} \biggr) \leq\frac{\pi e^{\mu L}}{L} \cdot
w_N e^{-\mu w_N}\nonumber
\\[-8pt]
\\[-8pt] &=& o(N (\log N)^2).
\nonumber
\end{eqnarray}
Observe that $Y_N(s(\log N)^3)/(N (\log N)^3)$ converges in probability
to zero, which is true by assumption when $s = 0$ and by Proposition
\ref{GNProp} when $s > 0$. Therefore, given any subsequence
$(N_j)_{j=1}^{\infty}$, there is a further subsequence
$(N_{j_k})_{k=1}^{\infty}$ such that $Y_{N_{j_k}}(s(\log
N_{j_k})^3)/(N_{j_k} (\log N_{j_k})^3)$ converges to zero almost
surely. It follows from (\ref{maxxn}) that
$m_{N_{j_k}}(s)/(N_{j_k}(\log N_{j_k})^2)$ converges to zero almost
surely, which implies the result.
\end{pf}

%
%le52 #&#
\begin{Lemma}\label{L1ptcv}
Assume the hypotheses of Theorem~\ref{bosz} hold.
Recall that $0 = t_0 < t_1 <\cdots< t_K$. Let $m \geq1$ and let $0 =
u_0 < u_1 <\cdots< u_m = 1$.
Then for each fixed $i$, with $1 \leq i \leq K$, we have
\[
(B^N_{0, t_i}(u_j))_{j=1}^m \Rightarrow(B_{0, t_i}(u_j))_{j=1}^m,
\]
where $\Rightarrow$ denotes convergence in distribution as $N
\rightarrow\infty$.
\end{Lemma}

\begin{pf}
It suffices to prove the joint convergence of the increments  $(B^N_{0,
t_i}(u_j) - B^N_{0, t_i}(u_{j-1}))_{j=1}^m$.
Define $L_0$ as in (\ref{Lidef}), and for $1 \leq j \leq m$, let
\[
S_j = \{L_0(u_{j-1}) + 1, L_0(u_{j-1})+2, \ldots, L_0(u_j)\}
\]
be the subset of particles in the population at time zero associated
with the quantiles in $[u_{j-1}, u_j)$. Note that the $S_j$ are
disjoint, and divide the population at time zero into $m$ subgroups.
We treat the positions of the particles in these $m$ subgroups as $m$
random starting configurations, to which we will apply
Lemma~\ref{ZNCSBPprime}.\vadjust{\goodbreak}

For $1 \leq j \leq m$, define the process $(Z_{j,N}(t), t \geq0)$ as
in Lemma~\ref{ZNCSBPprime}. We claim that the distribution of
%
%e174 ###
%
%e175 #&#
\begin{equation}
\label{Ejointin3}
\biggl(\frac{Z_{j,N}(0)}{N(\log N)^2} \biggr)_{j=1}^m
\end{equation}
converges as $N \rightarrow\infty$ to some probability measure $\rho$
on $[0, \infty)^m$. Here $\rho$ has the distribution\vspace*{1pt} of $(\delta_j
X)_{j=1}^m$, where $\delta_j = u_j - u_{j-1}$ for $1 \le j \le m$ and
$X$ has distribution $\nu$. To check that this convergence holds, note that
\[
|{Z_{j,N}(0)} - \delta_j Z_N(0) | \le2 m_N(0),
\]
where $m_N(0)$ is defined as in Lemma~\ref{maxcont} and the error term
$2m_N(0)$ comes from the fact that $\sum_{k=1}^j e^{\mu x_{k,0}}\sin
(\pi x_{k,0}/L)$ increases discontinuously with~$j$. In view of Lemma
\ref{maxcont}, the convergence of the distribution of (\ref{Ejointin3})
to $\rho$ follows by Slutsky's theorem (see Corollary 3.3 in Chapter 3
of~\cite{ek86}) and Proposition~\ref{ZNCSBP}. Therefore, the hypotheses
of Lemma~\ref{ZNCSBPprime} are satisfied.

Assume for now that $i \leq K-1$. By Lemma~\ref{ZNCSBPprime},
\[
\biggl(\frac{Z_{j,N}(t_i(\log N)^3)}{N(\log N)^2} \biggr)_{j=1}^m
\Rightarrow(Z_j(t_i))_{j=1}^m,
\]
where $\{(Z_j(t))_{j=1}^m, t \geq0\}$ is defined as in Lemma \ref
{ZNCSBPprime}. Thus for any $\alpha>0$,
%
%e175 ###
%
%e176 #&#
\begin{equation}\label{Eprevent}\quad
\biggl( \frac{Z_{j,N}(t_i(\log N)^3)}{Z_N(t_i( \log N)^3) \vee\alpha N
(\log N)^2} \biggr)_{j=1}^m \Rightarrow\biggl( \frac{Z_j(t_i)}{\alpha
\vee\sum_{k=1}^m Z_k(t_i)} \biggr)_{j=1}^m.
\end{equation}
Choose $\gamma>0$, and let $\alpha$ be such that $P(Z(t_i) < \alpha)
\le\gamma$, where $(Z(t), t \geq0)$ is a continuous-state branching
process with branching mechanism $\Psi$ and initial distribution $\nu$,
which is possible because $\nu(\{0\}) = 0$ and $(Z(t), t \geq0)$ never
goes extinct. Thus, by Proposition~\ref{ZNCSBP} we have
for $N$ large enough,
\[
P\bigl(Z_N(t_i (\log N)^3) < \alpha N (\log N)^2\bigr) \le2\gamma.
\]
Now fix $f_1, \ldots, f_m$, some arbitrary bounded and continuous test
functions on $[0,1]$ and let $M= \|f_1\|\cdots\|f_m\|$. Thus, we have
\begin{eqnarray*}
&& \Biggl|E \Biggl[\prod_{j=1}^m f_j \bigl( B^N_{0, t_i}(u_j) - B^N_{0,
t_i}(u_{j-1}) \bigr) \Biggr] - E \Biggl[ \prod_{j=1}^m f_j \biggl( \frac
{Z_j(t_i)}{\sum_{k=1}^m Z_k(t_i)} \biggr) \Biggr] \Biggr| \\
&& \qquad \le\Biggl|E \Biggl[ \prod_{j=1}^m f_j \biggl(\frac
{Z_{j,N}(t_i(\log N)^3)}{\alpha N (\log N)^2 \vee Z_N(t_i (\log
N)^3)} \biggr) \Biggr]\\
&&\hspace*{48pt} \qquad \hphantom{\le\Biggl|}{} -
E \Biggl[\prod_{j=1}^m f_j \biggl( \frac{Z_j(t_i)}{\alpha\vee\sum
_{k=1}^m Z_k (t_i)} \biggr) \Biggr] \Biggr| \\
&& \qquad \quad{} + M P\bigl(Z_N(t_i (\log N)^3) < \alpha N (\log N)^2\bigr) + M
P\bigl(Z(t_i) < \alpha\bigr).
\end{eqnarray*}
Taking the limsup of both sides, we find that the first term in the
right-hand side of the above inequality converges to 0 by (\ref
{Eprevent}), and the second and third terms are respectively smaller
than $2M \gamma$ and $M \gamma$. Since $\gamma> 0$ is arbitrary, and since
\[
\biggl(\frac{Z_j(t_i)}{\sum_{k=1}^m Z_k(t_i)} \biggr)_{j=1}^m
\]
has the same distribution as $(B_{0, t_i}(u_j) - B_{0,
t_i}(u_{j-1}))_{j=1}^m$, this finishes the proof when $1 \leq i \leq K-1$.

The proof when $i = K$ is the same, except $Z_{j,N}(t_i(\log
N)^3)/(N(\log N)^2)$ needs to be replaced throughout the argument by
$M_{j,N}(t_i(\log N)^3)/(2 \pi N)$, where the processes $(M_{j,N}(t), t
\geq0)$ are defined as in Lemma~\ref{ZNCSBPprime}.
\end{pf}

%s7.4 ###
%s7.4 #&#
\subsection{Joint convergence of bridges}\label{sec74}

In this subsection we extend the convergence obtained in Lemma \ref
{L1ptcv} to the joint convergence of the finite-dimen\-sional
distributions of several bridges. We begin by establishing a result
about the convergence of the distribution of a single bridge,
conditional on the branching Brownian motion up to the starting point
of the bridge.

%
%le53 #&#
\begin{Lemma}\label{condconv}
Assume the hypotheses of Theorem~\ref{bosz} hold.
Recall that $0 = t_0 < t_1 <\cdots< t_K$. Let $m \geq1$, and let $0 =
u_0 < u_1 <\cdots< u_m$. Let $f\dvtx \break [0,1]^{m+1} \rightarrow\R$ be
bounded and continuous. For $0 \leq i \leq K-1$, we have
\begin{eqnarray*}
&&E\bigl[f(B^N_{t_i, t_{i+1}}(u_0),\ldots, B^N_{t_i, t_{i+1}}(u_m))|{\cal
F}_{t_i (\log N)^3}\bigr] \\
&& \qquad \rightarrow_p E[f(B_{t_i, t_{i+1}}(u_0),\ldots,
B_{t_i, t_{i+1}}(u_m))],
\end{eqnarray*}
where $\rightarrow_p$ denotes convergence in probability as $N
\rightarrow\infty$.
\end{Lemma}

\begin{pf}
Let $(Z(t), t \geq0)$ be a continuous-state branching process with
branching mechanism $\Psi$ and initial distribution $\nu$. By
Proposition~\ref{ZNCSBP},
\[
\frac{Z_N(t_i(\log N)^3)}{N(\log N)^2} \Rightarrow Z(t_i).
\]
Also, we have $Y_N(t_i (\log N)^3)/(N (\log N)^3) \rightarrow_p 0$ by
assumption if $i = 0$ and by Proposition~\ref{GNProp} if $i \geq1$.
Therefore, by Skorohod's Representation Theorem, the branching Brownian
motion processes $(X_N, N \geq1)$ can be constructed on a single
probability space in such a way that $Z_N(t_i(\log N)^3)/\break(N(\log N)^2)
\rightarrow Z(t_i)$ a.s. and $Y_N(t_i (\log N)^3)/(N(\log N)^3)
\rightarrow0$ a.s.
Furthermore, it can be arranged that the processes $X_N$ evolve
independently of one another after time $t_i (\log N)^3$.

Let ${\tilde{\cal F}}_t = \sigma(X_N(s), N \geq1, 0 \leq s \leq t)$
be the $\sigma$-field generated by all the information up to time $t$
by all processes.
By the Markov property, conditional on ${\tilde{\cal F}}_{t_i (\log
N)^3}$, the process $X_N$ evolves after time $t_i (\log N)^3$ like a
branching Brownian motion with absorption whose initial configuration
is that of $X_N(t_i (\log N)^3)$.\vadjust{\goodbreak}
Therefore, we can apply Lemma~\ref{L1ptcv}, with $t_{i+1} - t_i$
playing the role of $t_i$ in Lemma~\ref{L1ptcv}, to get that on this
probability space
\begin{eqnarray*}
&&E\bigl[f(B^N_{t_i, t_{i+1}}(u_0),\ldots, B^N_{t_i, t_{i+1}}(u_m))|{\tilde
{\cal F}}_{t_i (\log N)^3}\bigr] \\
&& \qquad \rightarrow E[f(B_{0, t_{i+1} -
t_i}(u_0),\ldots, B_{0, t_{i+1} - t_i}(u_m))] \qquad\mbox{a.s.}
\end{eqnarray*}
The result follows because the bridges $B_{0, t_{i+1} - t_i}$ and
$B_{t_i, t_{i+1}}$ have the same law.
\end{pf}

%
%le54 #&#
\begin{Lemma} \label{Lmulticv}
Assume the hypotheses of Theorem~\ref{bosz} hold.
Recall that $0 = t_0 \leq t_1 <\cdots< t_K$ and let $0 = u_0 < u_1
\leq\cdots< u_m \leq1$. Then
%
%e176 ###
%
%e177 #&#
\begin{equation}\label{jointconv}
(B^N_{t_i, t_{i+1}}(u_j))\mathop{\mathop{}_{0 \leq i \le K-1 }}_{ 1
\leq j \leq
m} \Rightarrow(B_{t_i, t_{i+1}}(u_j))\mathop{\mathop{}_{0 \leq i
\leq K-1 }}_{
1 \leq j \leq m},
\end{equation}
where the bridges $B_{t_i, t_{i+1}}$, $0\le i \le K-1$, are independent.
\end{Lemma}

\begin{pf}
We proceed by induction. The convergence of $(B^N_{t_0, t_1}(u_j))_{1
\leq j \leq m}$ to $(B_{t_0, t_1}(u_j))_{1 \leq j \leq m}$ is a
consequence of Lemma~\ref{L1ptcv}. Thus assume\vspace*{1pt} that the convergence
(\ref{jointconv}) holds for $0 \leq i \leq k-1$ with $2\le k \le K - 1$.
Let $f_1, \ldots, f_k: [0,1]^{m+1} \to\R$ be bounded continuous functions.
By Proposition~\ref{ZNCSBP}, we know that
%
%e177 ###
%
%e178 #&#
\begin{equation}
(Z_{N}(t_1(\log N)^3), \ldots, Z_N(t_k(\log N)^3)) \Rightarrow(Z(t_1),
\ldots, Z(t_k)),
\end{equation}
where $(Z(t), t \geq0)$ is a continuous-state branching process with
branching mechanism $\Psi$ and initial distribution $\nu$.
Let ${\tilde{\cal F}}_t = \sigma(X_N(s), N \geq1, 0 \leq s \leq t)$
be the $\sigma$-field generated by the information up to time $t$. To
simplify notation, we write $\beta_i^N = (B^N_{t_i, t_{i+1}}(u_j))_{1
\le j \le m}$ and $\beta_i = (B_{t_i, t_{i+1}}(u_j))_{j=1}^m$. Since
$\beta_k^N$ is conditionally independent of
$\beta_1^N , \ldots, \beta_{k-1}^N$ given ${\cal F}_{t_{k-1}(\log N)^3}$,
%
%e178 ###
%
%e179 #&#
\begin{eqnarray}\label{Eindep0}
E \Biggl[\prod_{i=1}^k f_i(\beta_{i}^N) \Biggr] = E\Biggl [ \Biggl( \prod
_{i=1}^{k-1} f_i(\beta_{i}^N) \Biggr) E \bigl[ f_k(\beta^N_k)
|{\tilde{\cal F}}_{t_{k-1}(\log N)^3} \bigr] \Biggr].
\end{eqnarray}
Lemma~\ref{condconv} states that
%
%e179 ###
%
%e180 #&#
\begin{equation}
\label{Eindep1}E\bigl[f_k(\beta^N_k)|{\tilde{\cal F}}_{t_{k-1}(\log N)^3}\bigr]
\rightarrow_p E[f_k(\beta_k)],
\end{equation}
where $\rightarrow_p$ denotes convergence in probability as $N
\rightarrow\infty$.
Using the identity of real numbers
\[
x' y' - xy = x'(y' - y) +y(x'-x)
\]
in (\ref{Eindep0}) with $x' = \prod_{i=1}^{k-1} f_i(\beta_{i}^N)$, $y'
= E[f_k(\beta_{k}^N)|{\tilde{\cal F}}_{t_{k-1}(\log N)^3}]$, $x=\break
\prod_{i=1}^{k-1} E[f_i(\beta_i)]$ and $y = E[f_k(\beta_k)]$, and then
taking the expectation, we get
\begin{eqnarray*}
&&E \Biggl[\prod_{i=1}^k f_i(\beta_{i}^N) \Biggr] - \prod_{i=1}^k E[f_i(\beta
_i)]\\
&& \qquad  = E \Biggl[\prod_{i=1}^{k-1} f_i(\beta_{i}^N)
\bigl( E \bigl[ f_k(\beta^N_k) |{\tilde{\cal
F}}_{t_{k-1}(\log N)^3} \bigr] - E[f_k(\beta_k)] \bigr) \Biggr] \\
&& \qquad  \quad {} + E[f_k(\beta_k)] \Biggl(E \Biggl[\prod_{i=1}^{k-1} f_i(\beta
_{i}^N) \Biggr] - \prod_{i=1}^{k-1}E[f_i(\beta_{i})] \Biggr).
\end{eqnarray*}
The first term on the right-hand side converges to 0 by the dominated
convergence theorem and (\ref{Eindep1}) since $f_1, \ldots, f_k$ are
bounded, and the second term converges to 0 by the induction
hypothesis. This finishes the proof of Lemma~\ref{Lmulticv}.
\end{pf}

%s7.5 ###
%s7.5 #&#
\subsection{Tightness}\label{sec75}

Our goal in this subsection is to prove the following tightness result.

%
%le55 #&#
\begin{Lemma}\label{Ltight}
Assume the hypotheses of Theorem~\ref{bosz} hold.
For $0 \leq i \leq\break K-1$, the sequence of random discrete bridges
$(B^N_{t_i, t_{i+1}}(u), 0 \le u \le1)$ is a tight sequence with
respect to the Skorohod topology.
\end{Lemma}

\begin{pf}
For $\delta> 0$ and a function $B\dvtx [0,1] \to[0,1]$, define
\[
w'(B, \delta) = \inf_{\{x_j\}} \max_j \sup_{x,y \in[x_j,x_{j+1})}
|B(x) - B(y)|,
\]
where the infimum is taken over all subdivisions $\{x_j\}$ of $[0,1]$ with
$0 = x_0 < x_1 <\cdots< x_m = 1$ and
$\min(x_{j+1} - x_j) \ge\delta$. It suffices to show (see Chapter~13
of~\cite{bill}) that for all $\eps> 0$,
there exists $\delta>0$ such that
%
%e180 ###
%
%e181 #&#
\begin{equation}\label{w}
\limsup_{N \rightarrow\infty} P\bigl(w'(B^N_{t_i,t_{i+1}}, \delta) \ge
\eps
\bigr) \le\eps.
\end{equation}

Assume for now that $i \leq K-2$.
To prove (\ref{w}), we need to show that two jumps do not occur very
close to one another. Let $\eps>0$. Let $(Z(t), t \geq0)$ be a
continuous-state branching process with branching mechanism $\Psi$ and
initial distribution $\nu$. Since $\nu(\{0\}) = 0$, the
continuous-state branching process does not explode or go extinct, so
we can fix $0 < a < 1$ such that
\[
P\bigl(a < Z(t_{i+1}) < 1/a\bigr) \geq1 - \eps/4.
\]
Let $A(a,N)$ be the event that $a N (\log N)^2 < Z_N(t_{i+1}(\log N)^3)
< a^{-1} N (\log N)^2$. By Proposition~\ref{ZNCSBP}, we can choose
$N_0$ so that for all $N \ge N_0$, we have\break $P(A(a,N)) \geq1 - \eps/2$.

For $0 \leq x \leq1$, let
\begin{eqnarray*}
Z^N_{t_i, t_{i+1}}(x) &=& Z_N(t_{i+1} (\log N)^3) B^N_{t_i,
t_{i+1}}(x)\\
&=&\sum_{\ell= 1}^{L_i(x)} \sum_{m \in A_{\ell}(i,i+1)} e^{\mu x_{m,
i+1}} \sin\biggl( \frac{\pi x_{m, i+1}}{L} \biggr) .
\end{eqnarray*}
We now define our subdivision $\{x_j\}$. Let $x_0 = 0$, and for $j \geq
1$ such that $x_{j-1}<1$, let
\[
x_j = 1 \wedge\min\{x \ge0\dvtx  Z^N_{t_i,t_{i+1}}(x) -
Z^N_{t_i,t_{i+1}}(x_{j-1}) \geq a \eps N (\log N)^2\}.
\]
Since $P(A(a,N)) \geq1 - \eps/2$, and since this subdivision ensures
that\break $|B^N_{t_i, t_{i+1}}(x) - B^N_{t_i, t_{i+1}}(y)| < \eps$ for all
$x,y \in[x_j, x_{j+1})$ on the event $A(a,N)$, it remains only to show
that there is a $\delta>0$ such that
\[
\limsup_{N \rightarrow\infty} P \Bigl( A(a,N) \cap\Bigl\{ \min_{j}(x_j -
x_{j-1}) < \delta\Bigr\} \Bigr) \le\eps/2.
\]
Let $D_j$ be the event that $x_j \leq1 - \delta$. On the event
$A(a,N)$, there can be at most $1/\eps a^2$ values of $x_j$ less than
$1$. Also, on the event $D_j$, we have $x_j - x_{j-1} \leq\delta$ if
and only if $Z^N_{t_i, t_{i+1}}(x_{j-1} + \delta) - Z^N_{t_i,
t_{i+1}}(x_{j-1}) \geq a \eps N (\log N)^2$. Therefore, it suffices to
show that there exists $\delta> 0$ such that
\[
\limsup_{N \to\infty}P\bigl(Z^N_{t_i, t_{i+1}}(1) - Z^N_{t_i, t_{i+1}}(1 -
\delta) \geq a \eps N (\log N)^2\bigr) \leq\eps/4
\]
and for all $0 \le i \le(1/\eps a^2) - 1$,
\[
\limsup_{N \to\infty}P\bigl(D_i \cap\{Z^N_{t_i, t_{i+1}}(x_{j-1} +
\delta)
- Z^N_{t_i, t_{i+1}}(x_{j-1}) \geq a \eps N (\log N)^2\}\bigr) \le\eps^2a^2/4.
\]
In view of Lemma~\ref{maxcont}, both of these statements follow from an
application of Proposition~\ref{ZNCSBP}, in which the distribution of
$\delta Z(t_i)$ plays the role of $\nu$.

If $i = K-1$, the proof proceeds in the same way, except that instead
of working with $Z^N_{t_i, t_{i+1}}$, we define $M^N_{t_{K-1}, t_K}(x)
= M_N(t_K(\log N)^3) B_{t_{K-1}, t_K}(x)$.
The subdivision is defined by $x_0 = 0$ and, for $j \geq1$,
\[
x_j = 1 \wedge\min\{x \ge0\dvtx  M^N_{t_{K-1},t_K}(x) -
M^N_{t_{K-1},t_K}(x_{j-1}) \geq2 \pi a \eps N\}.
\]
The proof concludes with an application of Theorem~\ref{MNCSBP} rather
than Proposition~\ref{ZNCSBP}.
\end{pf}

Because the tightness of each sequence $(B^N_{t_i, t_{i+1}}(u), 0 \le u
\le1)$ implies the joint tightness of the $K$ sequences of bridges,
Lemmas~\ref{Lmulticv} and~\ref{Ltight} combine to yield the following
corollary.

%
%co56 #&#
\begin{Cor}\label{skbridge}
The sequence of processes $((B^N_{t_0, t_1}(u), B^N_{t_1,
t_2}(u),\ldots,\break  B^N_{t_{K-1}, t_K}(u)), 0 \le u \le1)$ converges in
the Skorohod\vspace*{1pt}
topology to $((B_{t_0, t_1}(u),\break  B_{t_1, t_2}(u),\ldots, B_{t_{K-1},
t_K}(u)), 0 \le u \le1)$.
\end{Cor}

%s7.6 ###
%s7.6 #&#
\subsection{Coalescence}\label{sec76}

Let $D$ be the set of functions $f\dvtx  [0,1] \rightarrow\R$ that are
right continuous and have left limits. Let $\rho$ denote the Skorohod
metric on $D$. Let $\Lambda$ denote the set of functions $\lambda:
[0,1] \rightarrow[0,1]$ that are continuous and strictly increasing
and satisfy $\lambda(0) = 0$ and $\lambda(1) = 1$. Recall (see Chapter
12 of~\cite{bill}) that\vadjust{\goodbreak} if $f, f_1, f_2, \dots$ are functions in $D$,
then $\lim_{n \rightarrow\infty} \rho(f_n, f) = 0$ if and only if
there exists a sequence of functions $(\lambda_n)_{n=1}^{\infty}$ in
$\Lambda$ such that
%
%e181 ###
%
%e182 #&#
\begin{equation}\label{sk1}
\lim_{n \rightarrow\infty} \sup_{0 \leq t \leq1} |f_n(\lambda
_n(t)) - f(t)| = 0
\end{equation}
and
%
%e182 ###
%
%e183 #&#
\begin{equation}\label{sk2}
\lim_{n \rightarrow\infty} \sup_{0 \leq t \leq1} |\lambda_n(t) -
t| = 0.
\end{equation}

The lemma below is similar to Lemma 1 of~\cite{beleg1} but differs in
that we do not require the processes $B_N$ to have exchangeable increments.

%
%le57 #&#
\begin{Lemma}\label{onebr}
Suppose $b, b_1, b_2, \dots$ are functions from $[0,1]$ to $[0,1]$ that
are nondecreasing and right continuous and have left limits at every
point other than~$0$. Suppose $\lim_{N \rightarrow\infty} \rho(b_N, b)
= 0$, where $\rho$ denotes the Skorohod metric. Suppose
$(x_N)_{N=1}^{\infty}$ and $(y_N)_{N=1}^{\infty}$ are sequences in
$[0,1]$ such that $x_N \rightarrow x$ and $y_N \rightarrow y$ as $N
\rightarrow\infty$. Suppose $x$ and $y$ are not in the closure of the
range of $b$. Then for sufficiently large $N$ we have $b_N^{-1}(x_N) =
b_N^{-1}(y_N)$ if and only if $b^{-1}(x) = b^{-1}(y)$. Furthermore,
%
%e183 ###
%
%e184 #&#
\begin{equation}\label{skconv}
\lim_{N \rightarrow\infty} b_N^{-1}(x_N) = b^{-1}(x).
\end{equation}
\end{Lemma}

\begin{pf}
Because $x$ is not in the closure of the range of $b$, there exists
some maximal open interval $(u, v)$ with $u < x < v$ such that $(u, v)$
does not intersect the range of $b$. For sufficiently small $\delta$,
we have $u + 2 \delta< x < v + 2 \delta$, which implies $u + \delta<
x_N < v - \delta$ for sufficiently large $N$. By condition (\ref{sk1})
applied to $b_N$ and $b$, for sufficiently large $N$ the interval $(u +
\delta, v - \delta)$ does not intersect the range of $b_N$. Therefore,
there exists $\gamma_N$ such that $b_N(\gamma_N) \geq v - \delta$ and
$b_N(\gamma_N-) \leq u + \delta$. Then $b_N^{-1}(x_N) = \gamma_N$ for
sufficiently large $N$. Also, there is a sequence of functions
$(\lambda
_N)_{N=1}^{\infty}$ in $\Lambda$ such that $\lambda_N(b^{-1}(x)) =
\gamma_N$ for sufficiently large $N$ by (\ref{sk1}) and therefore
$\lim
_{N \rightarrow\infty} \gamma_N = b^{-1}(x)$ by (\ref{sk2}). Result
(\ref{skconv}) follows.

Suppose $b^{-1}(x) = b^{-1}(y)$. Because $b$ is right continuous with
left limits, we have $u < y < v$. Arguing as above, we have
$b_N^{-1}(y_N) = \gamma_N$ for sufficiently large $N$, and thus
$b_N^{-1}(x_N) = b_N^{-1}(y_N)$ for sufficiently large $N$.
Alternatively, suppose $b^{-1}(x) \neq b^{-1}(y)$. We may assume
without loss of generality that $x < y$. Then $y > v$, and there is
some open interval $(r, s)$ with $v < r < y < s$ such that $(r, s)$
does not intersect the range of $b$. As above, there exists $\delta>
0$ and $\xi_N$ such that for sufficiently large $N$, we have $b_N(\xi
_N) \geq s - \delta$, $b_N(\xi_N-) \leq r + \delta$, and $b_N^{-1}(y_N)
= \xi_N$. Therefore, $b_N^{-1}(x_N) \neq b_N^{-1}(y_N)$ for
sufficiently large $N$, and the lemma follows.
\end{pf}

\begin{pf*}{Proof of Theorem~\ref{bosz}}
Fix times $0 = t_0 < t_1 <\cdots< t_K = t$. By Corollary~\ref
{skbridge} and Skorohod's representation theorem, we may work on a
probability space on which the sequence of discrete bridges\vspace*{-1pt}
$((B^N_{t_0, t_1}(u), B^N_{t_1, t_2}(u),\ldots,\break  B^N_{t_{K-1}, t_K}(u)),
0 \le u \le1)$ converges almost surely to $((B_{t_0, t_1}(u), B_{t_1,
t_2}(u),\ldots,\break  B_{t_{K-1}, t_K}(u)), 0 \le u \le1)$. Note that in
this setting, almost sure convergence means that $\rho(B^N_{t_i,
t_{i+1}},  B_{t_i, t_{i+1}}) \rightarrow0$ as $N \rightarrow\infty$
for $i = 0, 1,\ldots, K-1$, where $\rho$ denotes the Skorohod metric.

Fix a positive integer $n$, and let $U_1,\ldots, U_n$ be independent
random variables having the uniform distribution on $[0, 1]$. For $0
\leq i \leq K-1$, define the partition $\pi(B^N_{t_i, t}) = \pi
(B^N_{t_i, t_K})$ to be the partition of $\{1,\ldots, n\}$ such that
$i$ and $j$ are in the same block of the partition if and only if
$(B^N_{t_i, t_K})^{-1}(U_i) = (B^N_{t_i, t_K})^{-1}(U_j)$. Likewise,
define $\pi(B_{t_i, t}) = \pi(B_{t_i, t_K})$ to be the partition of
$\{
1,\ldots, n\}$ such that $i$ and $j$ are in the same block of the
partition if and only if $B_{t_i, t_K}^{-1}(U_i) = B_{t_i, t_K}^{-1}(U_j)$.
It follows from the definition of the processes $B^N_{t_i, t_K}$ that
$i$ and $j$ are in the same block of the partition if and only if the
individuals who are in positions $\lceil U_i M_N(t_K(\log N)^3) \rceil$
and $\lceil U_j M_N(t_K(\log N)^3) \rceil$ in the lexicographical order
at time $t_K (\log N)^3$ are descended from the same ancestor at time
$t_i (\log N)^3$. As a result, we have the equality in distribution
%
%e184 ###
%
%e185 #&#
\begin{eqnarray}\label{eqdist}
&&( \pi(B^N_{t_{K-1}, t_K}),\ldots, \pi(B^N_{t_0, t_K}) ) \nonumber
\\[-8pt]
\\[-8pt]&& \qquad =_d
\bigl( \Pi_N\bigl(2 \pi(t - t_{K-1})\bigr),\ldots, \Pi_N\bigl(2 \pi(t - t_0)\bigr) \bigr),
\nonumber
\end{eqnarray}
where $\Pi_N$ is the process defined in Theorem~\ref{bosz}. We note
that the sampling scheme here using the random variables $U_1,\ldots,
U_n$ corresponds to sampling with replacement from the individuals at
time $t_K(\log N)^3$, but the difference between sampling with and
without replacement is unimportant because the probability of sampling
the same individual twice tends to zero as $N \rightarrow\infty$.

We claim that for $0 \leq i \leq K-1$, almost surely
%
%e185 ###
%
%e186 #&#
\begin{equation}\label{claim167}
\pi(B^N_{t_i, t_K}) = \pi(B_{t_i, t_K})
\end{equation}
for sufficiently large $N$. Because we know the process $(\pi(B_{t -
s/2\pi, t}), 0 \leq s \leq2 \pi t)$ is the Bolthausen--Sznitman
coalescent run for time $t$, this claim in combination with (\ref
{eqdist}) will imply Theorem~\ref{bosz}.

We now prove (\ref{claim167}) by backward induction. Since the lengths
of the intervals of the complement of the range of $B_{t_{K-1}, t_K}$
have a Poisson--Dirichlet distribution and thus sum to 1 (see, e.g.,
Proposition 2 in~\cite{py97}), the closure of the range of $B_{t_{K-1},
t_K}$ has Lebesgue measure zero almost surely.
Therefore, almost surely $U_1,\ldots, U_n$ are not in the closure of
the range of $B_{t_{K-1}, t_K}$. It follows from Lemma~\ref{onebr} that
$\pi(B^N_{t_{K-1}, t_K}) = \pi(B_{t_{K-1}, t_K})$ for sufficiently
large $N$ almost surely. Furthermore,
\[
\lim_{N \rightarrow\infty}(B^N_{t_{K-1}, t_K})^{-1}(U_j) = B_{t_{K-1},
t_K}^{-1}(U_j)
\]
almost surely for $j = 1,\ldots, n$. Also, by Lemma 2 of~\cite{beleg1},
the random variables $B_{t_{K-1}, t_K}^{-1}(U_j)$ each have the uniform
distribution on $[0, 1]$.

For the induction step, suppose that for some $i = 2,\ldots, K-1$, the
following hold:
\begin{itemize}
\item We have ${ \lim_{N \rightarrow\infty}} (B^N_{t_i,
t_K})^{-1}(U_j) = B_{t_i, t_K}^{-1}(U_j)$ almost surely for $j =
1,\ldots, n$.

\item The random variables $B_{t_i, t_K}^{-1}(U_j)$ each have the
uniform distribution on $[0, 1]$.
\end{itemize}
Now $(B^N_{t_{i-1}, t_K})^{-1}(U_j) = (B^N_{t_{i-1},
t_i})^{-1}((B^N_{t_i, t_K})^{-1}(U_j))$ by Lemma~\ref{cocycle}, and
likewise $B_{t_{i-1}, t_K}^{-1}(U_j) = B_{t_{i-1}, t_i}^{-1}(B_{t_i,
t_K}^{-1}(U_j))$. Because the random variables\break $B_{t_i, t_K}^{-1}(U_j)$
each have the uniform distribution on $[0, 1]$ and are independent of
$B_{t_{i-1}, t_i}$, almost surely none of these random variables is in
the closure of the range of $B_{t_{i-1}, t_i}$. Since also $(B^N_{t_i,
t_K})^{-1}(U_j) \rightarrow B_{t_i, t_K}^{-1}(U_j)$ almost surely for
$j = 1,\ldots, n$, Lemma~\ref{onebr} implies that $\pi(B^N_{t_{i-1},
t_K}) = \pi(B_{t_{i-1}, t_K})$ for sufficiently large $N$ almost
surely. Furthermore, $(B^N_{t_{i-1}, t_K})^{-1}(U_j) \rightarrow
B_{t_{i-1}, t_K}^{-1}(U_j)$ almost surely for $j = 1,\ldots, n$. By
Lemma~2 of~\cite{beleg1}, the random variables $B_{t_{i-1},
t_K}^{-1}(U_j)$ each have the uniform distribution on $[0, 1]$. The
claim (\ref{claim167}) now follows by induction.
\end{pf*}

\section*{Acknowledgments}
The authors thank two referees for carefully reading the paper and
making a number of helpful comments.

%suskaldyti doi

% imsref loaded by smiklovaite, 2012-03-16 13:33:33
%

\printaddresses


\begin{thebibliography}{66}
% BibTex style file: ims.bst, 2011-05-30
% Default style options (sort=0,type=number).
% Used options (sort=1,type=number).

%b1 ###
%b1 #&#
\bibitem{aija11}
%
\begin{barticle}[mr]
\bauthor{\bsnm{A{\"{\i}}d{\'e}kon},~\bfnm{Elie}\binits{E.}} \AND
\bauthor{\bsnm{Jaffuel},~\bfnm{Bruno}\binits{B.}}
(\byear{2011}).
\btitle{Survival of branching random walks with absorption}.
\bjournal{Stochastic Process. Appl.}
\bvolume{121}
\bpages{1901--1937}.
\bid{doi={10.1016/j.spa.2011.04.006}, issn={0304-4149}, mr={2819234}}
\bptok{imsref}%
\end{barticle}
%
\endbibitem

%b2 ###
%b2 #&#
\bibitem{bago07}
%
\begin{barticle}[mr]
\bauthor{\bsnm{Basdevant},~\bfnm{Anne-Laure}\binits{A.-L.}} \AND
\bauthor{\bsnm{Goldschmidt},~\bfnm{Christina}\binits{C.}}
(\byear{2008}).
\btitle{Asymptotics of the allele frequency spectrum associated with the
{B}olthausen--{S}znitman coalescent}.
\bjournal{Electron. J. Probab.}
\bvolume{13}
\bpages{486--512}.
\bid{doi={10.1214/EJP.v13-494}, issn={1083-6489}, mr={2386740}}
\bptok{imsref}%
\end{barticle}
%
\endbibitem

%b3 ###
%b3 #&#
\bibitem{bego08}
%
\begin{barticle}[mr]
\bauthor{\bsnm{B{\'e}rard},~\bfnm{Jean}\binits{J.}} \AND
\bauthor{\bsnm{Gou{\'e}r{\'e}},~\bfnm{Jean-Baptiste}\binits{J.-B.}}
(\byear{2010}).
\btitle{Brunet--{D}errida behavior of branching-selection particle
systems on
the line}.
\bjournal{Comm. Math. Phys.}
\bvolume{298}
\bpages{323--342}.
\bid{doi={10.1007/s00220-010-1067-y}, issn={0010-3616}, mr={2669438}}
\bptok{imsref}%
\end{barticle}
%
\endbibitem

%b4 ###
%b4 #&#
\bibitem{bego09}
%
\begin{barticle}[mr]
\bauthor{\bsnm{B{\'e}rard},~\bfnm{Jean}\binits{J.}} \AND
\bauthor{\bsnm{Gou{\'e}r{\'e}},~\bfnm{Jean-Baptiste}\binits{J.-B.}}
(\byear{2011}).
\btitle{Survival probability of the branching random walk killed below
a linear
boundary}.
\bjournal{Electron. J. Probab.}
\bvolume{16}
\bpages{396--418}.
\bid{doi={10.1214/EJP.v16-861}, issn={1083-6489}, mr={2774095}}
\bptok{imsref}%
\end{barticle}
%
\endbibitem

%b5 ###
%b5 #&#
\bibitem{bbs}
%
\begin{barticle}[mr]
\bauthor{\bsnm{Berestycki},~\bfnm{Julien}\binits{J.}},
\bauthor{\bsnm{Berestycki},~\bfnm{Nathana{\"e}l}\binits{N.}} \AND
\bauthor{\bsnm{Schweinsberg},~\bfnm{Jason}\binits{J.}}
(\byear{2011}).
\btitle{Survival of near-critical branching {B}rownian motion}.
\bjournal{J. Stat. Phys.}
\bvolume{143}
\bpages{833--854}.
\bid{doi={10.1007/s10955-011-0224-9}, issn={0022-4715}, mr={2811463}}
\bptok{imsref}%
\end{barticle}
%
\endbibitem

%b6 ###
%b6 #&#
\bibitem{ensaios}
%
\begin{bbook}[mr]
\bauthor{\bsnm{Berestycki},~\bfnm{Nathana{\"e}l}\binits{N.}}
(\byear{2009}).
\btitle{Recent Progress in Coalescent Theory}.
\bseries{Ensaios Matem\'aticos [Mathematical Surveys]}
\bvolume{16}.
\bpublisher{Sociedade Brasileira de Matem\'atica}, \baddress{Rio de Janeiro}.
\bid{mr={2574323}}
\bptok{imsref}%
\end{bbook}
%
\endbibitem

%b7 ###
%b7 #&#
\bibitem{beleg00}
%
\begin{barticle}[mr]
\bauthor{\bsnm{Bertoin},~\bfnm{Jean}\binits{J.}} \AND
\bauthor{\bsnm{Le~Gall},~\bfnm{Jean-Fran{\c{c}}ois}\binits{J.-F.}}
(\byear{2000}).
\btitle{The {B}olthausen--{S}znitman coalescent and the genealogy of
continuous-state branching processes}.
\bjournal{Probab. Theory Related Fields}
\bvolume{117}
\bpages{249--266}.
\bid{doi={10.1007/s004400050008}, issn={0178-8051}, mr={1771663}}
\bptok{imsref}%
\end{barticle}
%
\endbibitem

%b8 ###
%b8 #&#
\bibitem{beleg1}
%
\begin{barticle}[mr]
\bauthor{\bsnm{Bertoin},~\bfnm{Jean}\binits{J.}} \AND
\bauthor{\bsnm{Le~Gall},~\bfnm{Jean-Fran{\c{c}}ois}\binits{J.-F.}}
(\byear{2003}).
\btitle{Stochastic flows associated to coalescent processes}.
\bjournal{Probab. Theory Related Fields}
\bvolume{126}
\bpages{261--288}.
\bid{doi={10.1007/s00440-003-0264-4}, issn={0178-8051}, mr={1990057}}
\bptok{imsref}%
\end{barticle}
%
\endbibitem

%b9 ###
%b9 #&#
\bibitem{bepi00}
%
\begin{barticle}[mr]
\bauthor{\bsnm{Bertoin},~\bfnm{Jean}\binits{J.}} \AND
\bauthor{\bsnm{Pitman},~\bfnm{Jim}\binits{J.}}
(\byear{2000}).
\btitle{Two coalescents derived from the ranges of stable subordinators}.
\bjournal{Electron. J. Probab.}
\bvolume{5}
\bpages{1--17 (electronic)}.
\bid{issn={1083-6489}, mr={1768841}}
\bptok{imsref}%
\end{barticle}
%
\endbibitem

%b10 ###
%b10 #&#
\bibitem{bigg}
%
\begin{barticle}[mr]
\bauthor{\bsnm{Biggins},~\bfnm{J.~D.}\binits{J.~D.}}
(\byear{1976}).
\btitle{The first- and last-birth problems for a multitype age-dependent
branching process}.
\bjournal{Adv. in Appl. Probab.}
\bvolume{8}
\bpages{446--459}.
\bid{issn={0001-8678}, mr={0420890}}
\bptok{imsref}%
\end{barticle}
%
\endbibitem

%b11 ###
%b11 #&#
\bibitem{bill}
%
\begin{bbook}[mr]
\bauthor{\bsnm{Billingsley},~\bfnm{Patrick}\binits{P.}}
(\byear{1999}).
\btitle{Convergence of Probability Measures},
\bedition{2nd} ed.
% Statistics}.
\bpublisher{Wiley}, \baddress{New York}.
\bid{doi={10.1002/9780470316962}, mr={1700749}}
\bptok{imsref}%
\end{bbook}
%
\endbibitem

%b12 ###
%b12 #&#
\bibitem{regvar}
%
\begin{bbook}[mr]
\bauthor{\bsnm{Bingham},~\bfnm{N.~H.}\binits{N.~H.}},
\bauthor{\bsnm{Goldie},~\bfnm{C.~M.}\binits{C.~M.}} \AND
\bauthor{\bsnm{Teugels},~\bfnm{J.~L.}\binits{J.~L.}}
(\byear{1987}).
\btitle{Regular Variation}.
\bseries{Encyclopedia of Mathematics and Its Applications}
\bvolume{27}.
\bpublisher{Cambridge Univ. Press}, \baddress{Cambridge}.
\bid{mr={0898871}}
\bptok{imsref}%
\end{bbook}
%
\endbibitem

%b13 ###
%b13 #&#
\bibitem{bbcemsw}
%
\begin{barticle}[mr]
\bauthor{\bsnm{Birkner},~\bfnm{Matthias}\binits{M.}},
\bauthor{\bsnm{Blath},~\bfnm{Jochen}\binits{J.}},
\bauthor{\bsnm{Capaldo},~\bfnm{Marcella}\binits{M.}},
\bauthor{\bsnm{Etheridge},~\bfnm{Alison}\binits{A.}},
\bauthor{\bsnm{M{\"o}hle},~\bfnm{Martin}\binits{M.}},
\bauthor{\bsnm{Schweinsberg},~\bfnm{Jason}\binits{J.}} \AND
\bauthor{\bsnm{Wakolbinger},~\bfnm{Anton}\binits{A.}}
(\byear{2005}).
\btitle{Alpha-stable branching and beta-coalescents}.
\bjournal{Electron. J. Probab.}
\bvolume{10}
\bpages{303--325 (electronic)}.
\bid{doi={10.1214/EJP.v10-241}, issn={1083-6489}, mr={2120246}}
\bptok{imsref}%
\end{barticle}
%
\endbibitem

%b14 ###
%b14 #&#
\bibitem{bosz98}
%
\begin{barticle}[mr]
\bauthor{\bsnm{Bolthausen},~\bfnm{E.}\binits{E.}} \AND
\bauthor{\bsnm{Sznitman},~\bfnm{A.~S.}\binits{A.~S.}}
(\byear{1998}).
\btitle{On {R}uelle's probability cascades and an abstract cavity method}.
\bjournal{Comm. Math. Phys.}
\bvolume{197}
\bpages{247--276}.
\bid{doi={10.1007/s002200050450}, issn={0010-3616}, mr={1652734}}
\bptok{imsref}%
\end{barticle}
%
\endbibitem

%b15 ###
%b15 #&#
\bibitem{bovkur}
%
\begin{bincollection}[mr]
\bauthor{\bsnm{Bovier},~\bfnm{Anton}\binits{A.}} \AND
\bauthor{\bsnm{Kurkova},~\bfnm{Irina}\binits{I.}}
(\byear{2007}).
\btitle{Much ado about {D}errida's {GREM}}.
In \bbooktitle{Spin Glasses}.
\bseries{Lecture Notes in Math.}
\bvolume{1900}
\bpages{81--115}.
\bpublisher{Springer}, \baddress{Berlin}.
\bid{doi={10.1007/978-3-540-40908-3_4}, mr={2309599}}
\bptok{imsref}%
\end{bincollection}
%
\endbibitem

%b16 ###
%b16 #&#
\bibitem{bram2}
%
\begin{barticle}[mr]
\bauthor{\bsnm{Bramson},~\bfnm{Maury}\binits{M.}}
(\byear{1983}).
\btitle{Convergence of solutions of the {K}olmogorov equation to travelling
waves}.
\bjournal{Mem. Amer. Math. Soc.}
\bvolume{44}
\bpages{iv+190}.
\bid{issn={0065-9266}, mr={0705746}}
\bptok{imsref}%
\end{barticle}
%
\endbibitem

%b17 ###
%b17 #&#
\bibitem{bram1}
%
\begin{barticle}[mr]
\bauthor{\bsnm{Bramson},~\bfnm{Maury~D.}\binits{M.~D.}}
(\byear{1978}).
\btitle{Maximal displacement of branching {B}rownian motion}.
\bjournal{Comm. Pure Appl. Math.}
\bvolume{31}
\bpages{531--581}.
\bid{issn={0010-3640}, mr={0494541}}
\bptok{imsref}%
\end{barticle}
%
\endbibitem

%b18 ###
%b18 #&#
\bibitem{bd2}
%
\begin{barticle}[auto:STB|2012/03/12|15:33:09]
\bauthor{\bsnm{Brunet},~\bfnm{E.}\binits{E.}} \AND
\bauthor{\bsnm{Derrida},~\bfnm{B.}\binits{B.}}
(\byear{1999}).
\btitle{Microscopic models of traveling wave equations}.
\bjournal{Comput. Phys. Comm.}
\bvolume{121-122}
\bpages{376--381}.
\bptok{imsref}%
\end{barticle}
%
\endbibitem

%b19 ###
%b19 #&#
\bibitem{bd1}
%
\begin{barticle}[mr]
\bauthor{\bsnm{Brunet},~\bfnm{Eric}\binits{E.}} \AND
\bauthor{\bsnm{Derrida},~\bfnm{Bernard}\binits{B.}}
(\byear{1997}).
\btitle{Shift in the velocity of a front due to a cutoff}.
\bjournal{Phys. Rev. E (3)}
\bvolume{56}
\bpages{2597--2604}.
\bid{doi={10.1103/PhysRevE.56.2597}, issn={1539-3755}, mr={1473413}}
\bptok{imsref}%
\end{barticle}
%
\endbibitem

%b20 ###
%b20 #&#
\bibitem{bd3}
%
\begin{barticle}[mr]
\bauthor{\bsnm{Brunet},~\bfnm{{\'E}ric}\binits{{\'E}.}} \AND
\bauthor{\bsnm{Derrida},~\bfnm{Bernard}\binits{B.}}
(\byear{2001}).
\btitle{Effect of microscopic noise on front propagation}.
\bjournal{J.~Stat. Phys.}
\bvolume{103}
\bpages{269--282}.
\bid{doi={10.1023/A:1004875804376}, issn={0022-4715}, mr={1828730}}
\bptok{imsref}%
\end{barticle}
%
\endbibitem

%b21 ###
%b21 #&#
\bibitem{bdmm1}
%
\begin{barticle}[mr]
\bauthor{\bsnm{Brunet},~\bfnm{E.}\binits{E.}},
\bauthor{\bsnm{Derrida},~\bfnm{B.}\binits{B.}},
\bauthor{\bsnm{Mueller},~\bfnm{A.~H.}\binits{A.~H.}} \AND
\bauthor{\bsnm{Munier},~\bfnm{S.}\binits{S.}}
(\byear{2006}).
\btitle{Noisy traveling waves: Effect of selection on genealogies}.
\bjournal{Europhys. Lett.}
\bvolume{76}
\bpages{1--7}.
\bid{doi={10.1209/epl/i2006-10224-4}, issn={0295-5075}, mr={2299937}}
\bptok{imsref}%
\end{barticle}
%
\endbibitem

%b22 ###
%b22 #&#
\bibitem{bdmm2}
%
\begin{barticle}[mr]
\bauthor{\bsnm{Brunet},~\bfnm{{\'E}.}\binits{{\'E}.}},
\bauthor{\bsnm{Derrida},~\bfnm{B.}\binits{B.}},
\bauthor{\bsnm{Mueller},~\bfnm{A.~H.}\binits{A.~H.}} \AND
\bauthor{\bsnm{Munier},~\bfnm{S.}\binits{S.}}
(\byear{2007}).
\btitle{Effect of selection on ancestry: An exactly soluble case and its
phenomenological generalization}.
\bjournal{Phys. Rev. E (3)}
\bvolume{76}
\bpages{041104}.
\bid{doi={10.1103/PhysRevE.76.041104}, issn={1539-3755}, mr={2365627}}
\bptok{imsref}%
\end{barticle}
%
\endbibitem

%b23 ###
%b23 #&#
\bibitem{cablam}
%
\begin{barticle}[mr]
\bauthor{\bsnm{Caballero},~\bfnm{Ma.~Emilia}\binits{M.~E.}},
\bauthor{\bsnm{Lambert},~\bfnm{Amaury}\binits{A.}} \AND
\bauthor{\bsnm{Uribe~Bravo},~\bfnm{Ger{\'o}nimo}\binits{G.}}
(\byear{2009}).
\btitle{Proof(s) of the {L}amperti representation of continuous-state branching
processes}.
\bjournal{Probab. Surv.}
\bvolume{6}
\bpages{62--89}.
\bid{doi={10.1214/09-PS154}, issn={1549-5787}, mr={2592395}}
\bptok{imsref}%
\end{barticle}
%
\endbibitem

%b24 ###
%b24 #&#
\bibitem{cw}
%
\begin{bbook}[mr]
\bauthor{\bsnm{Chung},~\bfnm{K.~L.}\binits{K.~L.}} \AND
\bauthor{\bsnm{Williams},~\bfnm{R.~J.}\binits{R.~J.}}
(\byear{1990}).
\btitle{Introduction to Stochastic Integration},
\bedition{2nd} ed.
\bpublisher{Birkh\"auser}, \baddress{Boston, MA}.
\bid{mr={1102676}}
\bptok{imsref}%
\end{bbook}
%
\endbibitem

%b25 ###
%b25 #&#
\bibitem{ds1}
%
\begin{barticle}[mr]
\bauthor{\bsnm{Derrida},~\bfnm{B.}\binits{B.}} \AND
\bauthor{\bsnm{Simon},~\bfnm{D.}\binits{D.}}
(\byear{2007}).
\btitle{The survival probability of a branching random walk in presence
of an
absorbing wall}.
\bjournal{Europhys. Lett. EPL}
\bvolume{78}
\bpages{Art. 60006}.
\bid{doi={10.1209/0295-5075/78/60006}, issn={0295-5075}, mr={2366713}}
\bptok{imsref}%
\end{barticle}
%
\endbibitem

%b26 ###
%b26 #&#
\bibitem{dk99}
%
\begin{barticle}[mr]
\bauthor{\bsnm{Donnelly},~\bfnm{Peter}\binits{P.}} \AND
\bauthor{\bsnm{Kurtz},~\bfnm{Thomas~G.}\binits{T.~G.}}
(\byear{1999}).
\btitle{Particle representations for measure-valued population models}.
\bjournal{Ann. Probab.}
\bvolume{27}
\bpages{166--205}.
\bid{doi={10.1214/aop/1022677258}, issn={0091-1798}, mr={1681126}}
\bptok{imsref}%
\end{barticle}
%
\endbibitem

%b27 ###
%b27 #&#
\bibitem{dimr07}
%
\begin{barticle}[mr]
\bauthor{\bsnm{Drmota},~\bfnm{Michael}\binits{M.}},
\bauthor{\bsnm{Iksanov},~\bfnm{Alex}\binits{A.}},
\bauthor{\bsnm{Moehle},~\bfnm{Martin}\binits{M.}} \AND
\bauthor{\bsnm{Roesler},~\bfnm{Uwe}\binits{U.}}
(\byear{2007}).
\btitle{Asymptotic results concerning the total branch length of the
{B}olthausen--{S}znitman coalescent}.
\bjournal{Stochastic Process. Appl.}
\bvolume{117}
\bpages{1404--1421}.
\bid{doi={10.1016/j.spa.2007.01.011}, issn={0304-4149}, mr={2353033}}
\bptok{imsref}%
\end{barticle}
%
\endbibitem

%b28 ###
%b28 #&#
\bibitem{durr96}
%
\begin{bbook}[mr]
\bauthor{\bsnm{Durrett},~\bfnm{Richard}\binits{R.}}
(\byear{1996}).
\btitle{Stochastic Calculus: A Practical Introduction}.
\bpublisher{CRC Press}, \baddress{Boca Raton, FL}.
\bid{mr={1398879}}
\bptok{imsref}%
\end{bbook}
%
\endbibitem

%b29 ###
%b29 #&#
\bibitem{dumay}
%
\begin{barticle}[mr]
\bauthor{\bsnm{Durrett},~\bfnm{Rick}\binits{R.}} \AND
\bauthor{\bsnm{Mayberry},~\bfnm{John}\binits{J.}}
(\byear{2010}).
\btitle{Evolution in predator-prey systems}.
\bjournal{Stochastic Process. Appl.}
\bvolume{120}
\bpages{1364--1392}.
\bid{doi={10.1016/j.spa.2010.03.011}, issn={0304-4149}, mr={2639750}}
\bptnote{check year}%
\bptok{imsref}%
\end{barticle}
%
\endbibitem

%b30 ###
%b30 #&#
\bibitem{durem}
%
\begin{barticle}[auto:STB|2012/03/12|15:33:09]
\bauthor{\bsnm{Durrett},~\bfnm{R.}\binits{R.}} \AND
\bauthor{\bsnm{Remenik},~\bfnm{D.}\binits{D.}}
(\byear{2011}).
\btitle{Brunet--Derrida particle systems, free boundary problems and Wiener--Hopf
equations}.
\bjournal{Ann. Probab.}
\bvolume{39}
\bpages{2043--2078}.
\bptok{imsref}%
\end{barticle}
%
\endbibitem

%b31 ###
%b31 #&#
\bibitem{ek86}
%
\begin{bbook}[mr]
\bauthor{\bsnm{Ethier},~\bfnm{Stewart~N.}\binits{S.~N.}} \AND
\bauthor{\bsnm{Kurtz},~\bfnm{Thomas~G.}\binits{T.~G.}}
(\byear{1986}).
\btitle{Markov Processes: Characterization and Convergence}.
%Probability
% and Mathematical Statistics}.
\bpublisher{Wiley}, \baddress{New York}.
\bid{doi={10.1002/9780470316658}, mr={0838085}}
\bptok{imsref}%
\end{bbook}
%
\endbibitem

%b32 ###
%b32 #&#
\bibitem{fz}
%
\begin{barticle}[mr]
\bauthor{\bsnm{Fang},~\bfnm{Ming}\binits{M.}} \AND
\bauthor{\bsnm{Zeitouni},~\bfnm{Ofer}\binits{O.}}
(\byear{2010}).
\btitle{Consistent minimal displacement of branching random walks}.
\bjournal{Electron. Commun. Probab.}
\bvolume{15}
\bpages{106--118}.
\bid{doi={10.1214/ECP.v15-1533}, issn={1083-589X}, mr={2606508}}
\bptok{imsref}%
\end{barticle}
%
\endbibitem

%b33 ###
%b33 #&#
\bibitem{fisher}
%
\begin{barticle}[auto:STB|2012/03/12|15:33:09]
\bauthor{\bsnm{Fisher},~\bfnm{R.~A.}\binits{R.~A.}}
(\byear{1937}).
\btitle{The wave of advance of advantageous genes}.
\bjournal{Ann. Eugenics}
\bvolume{7}
\bpages{355--369}.
\bptok{imsref}%
\end{barticle}
%
\endbibitem

%b34 ###
%b34 #&#
\bibitem{ghs}
%
\begin{barticle}[mr]
\bauthor{\bsnm{Gantert},~\bfnm{Nina}\binits{N.}},
\bauthor{\bsnm{Hu},~\bfnm{Yueyun}\binits{Y.}} \AND
\bauthor{\bsnm{Shi},~\bfnm{Zhan}\binits{Z.}}
(\byear{2011}).
\btitle{Asymptotics for the survival probability in a killed branching random
walk}.
\bjournal{Ann. Inst. Henri Poincar\'e Probab. Stat.}
\bvolume{47}
\bpages{111--129}.
\bid{doi={10.1214/10-AIHP362}, issn={0246-0203}, mr={2779399}}
\bptok{imsref}%
\end{barticle}
%
\endbibitem

%b35 ###
%b35 #&#
\bibitem{goma05}
%
\begin{barticle}[mr]
\bauthor{\bsnm{Goldschmidt},~\bfnm{Christina}\binits{C.}} \AND
\bauthor{\bsnm{Martin},~\bfnm{James~B.}\binits{J.~B.}}
(\byear{2005}).
\btitle{Random recursive trees and the {B}olthausen--{S}znitman coalescent}.
\bjournal{Electron. J. Probab.}
\bvolume{10}
\bpages{718--745 (electronic)}.
\bid{doi={10.1214/EJP.v10-265}, issn={1083-6489}, mr={2164028}}
\bptok{imsref}%
\end{barticle}
%
\endbibitem

%b36 ###
%b36 #&#
\bibitem{grey}
%
\begin{barticle}[mr]
\bauthor{\bsnm{Grey},~\bfnm{D.~R.}\binits{D.~R.}}
(\byear{1974}).
\btitle{Asymptotic behaviour of continuous time, continuous state-space
branching processes}.
\bjournal{J. Appl. Probab.}
\bvolume{11}
\bpages{669--677}.
\bid{issn={0021-9002}, mr={0408016}}
\bptok{imsref}%
\end{barticle}
%
\endbibitem

%b37 ###
%b37 #&#
\bibitem{hamm}
%
\begin{barticle}[mr]
\bauthor{\bsnm{Hammersley},~\bfnm{J.~M.}\binits{J.~M.}}
(\byear{1974}).
\btitle{Postulates for subadditive processes}.
\bjournal{Ann. Probab.}
\bvolume{2}
\bpages{652--680}.
\bid{mr={0370721}}
\bptok{imsref}%
\end{barticle}
%
\endbibitem

%b38 ###
%b38 #&#
\bibitem{hh07}
%
\begin{barticle}[mr]
\bauthor{\bsnm{Harris},~\bfnm{J.~W.}\binits{J.~W.}} \AND
\bauthor{\bsnm{Harris},~\bfnm{S.~C.}\binits{S.~C.}}
(\byear{2007}).
\btitle{Survival probabilities for branching {B}rownian motion with
absorption}.
\bjournal{Electron. Commun. Probab.}
\bvolume{12}
\bpages{81--92 (electronic)}.
\bid{doi={10.1214/ECP.v12-1259}, issn={1083-589X}, mr={2300218}}
\bptok{imsref}%
\end{barticle}
%
\endbibitem

%b39 ###
%b39 #&#
\bibitem{hhk06}
%
\begin{barticle}[mr]
\bauthor{\bsnm{Harris},~\bfnm{J.~W.}\binits{J.~W.}},
\bauthor{\bsnm{Harris},~\bfnm{S.~C.}\binits{S.~C.}} \AND
\bauthor{\bsnm{Kyprianou},~\bfnm{A.~E.}\binits{A.~E.}}
(\byear{2006}).
\btitle{Further probabilistic analysis of the
{F}isher--{K}olmogorov--{P}etrovskii--{P}iscounov equation: One sided
travelling-waves}.
\bjournal{Ann. Inst. Henri Poincar\'e Probab. Stat.}
\bvolume{42}
\bpages{125--145}.
\bid{doi={10.1016/j.anihpb.2005.02.005}, issn={0246-0203}, mr={2196975}}
\bptok{imsref}%
\end{barticle}
%
\endbibitem

%b40 ###
%b40 #&#
\bibitem{har99}
%
\begin{barticle}[mr]
\bauthor{\bsnm{Harris},~\bfnm{Simon~C.}\binits{S.~C.}}
(\byear{1999}).
\btitle{Travelling-waves for the {FKPP} equation via probabilistic arguments}.
\bjournal{Proc. Roy. Soc. Edinburgh Sect. A}
\bvolume{129}
\bpages{503--517}.
\bid{doi={10.1017/S030821050002148X}, issn={0308-2105}, mr={1693633}}
\bptok{imsref}%
\end{barticle}
%
\endbibitem

%b41 ###
%b41 #&#
\bibitem{jaffuel}
%
\begin{bmisc}[auto:STB|2012/03/12|15:33:09]
\bauthor{\bsnm{Jaffuel},~\bfnm{B.}\binits{B.}}
(\byear{2009}).
\bhowpublished{The critical barrier for the survival of the branching
random walk with
absorption. Preprint. Available at arXiv:\arxivurl{0911.2227}}.
\bptok{imsref}%
\end{bmisc}
%
\endbibitem

%b42 ###
%b42 #&#
\bibitem{jirina}
%
\begin{barticle}[mr]
\bauthor{\bsnm{Ji{\v{r}}ina},~\bfnm{Miloslav}\binits{M.}}
(\byear{1958}).
\btitle{Stochastic branching processes with continuous state space}.
\bjournal{Czechoslovak Math. J.}
\bvolume{8}({83})
\bpages{292--313}.
\bid{issn={0011-4642}, mr={0101554}}
\bptok{imsref}%
\end{barticle}
%
\endbibitem

%b43 ###
%b43 #&#
\bibitem{kesten}
%
\begin{barticle}[mr]
\bauthor{\bsnm{Kesten},~\bfnm{Harry}\binits{H.}}
(\byear{1978}).
\btitle{Branching {B}rownian motion with absorption}.
\bjournal{Stochastic Process. Appl.}
\bvolume{7}
\bpages{9--47}.
\bid{issn={0304-4149}, mr={0494543}}
\bptok{imsref}%
\end{barticle}
%
\endbibitem

%b44 ###
%b44 #&#
\bibitem{king}
%
\begin{barticle}[mr]
\bauthor{\bsnm{Kingman},~\bfnm{J.~F.~C.}\binits{J.~F.~C.}}
(\byear{1975}).
\btitle{The first birth problem for an age-dependent branching process}.
\bjournal{Ann. Probab.}
\bvolume{3}
\bpages{790--801}.
\bid{mr={0400438}}
\bptok{imsref}%
\end{barticle}
%
\endbibitem

%b45 ###
%b45 #&#
\bibitem{king82}
%
\begin{barticle}[mr]
\bauthor{\bsnm{Kingman},~\bfnm{J.~F.~C.}\binits{J.~F.~C.}}
(\byear{1982}).
\btitle{The coalescent}.
\bjournal{Stochastic Process. Appl.}
\bvolume{13}
\bpages{235--248}.
\bid{doi={10.1016/0304-4149(82)90011-4}, issn={0304-4149}, mr={0671034}}
\bptok{imsref}%
\end{barticle}
%
\endbibitem

%b46 ###
%b46 #&#
\bibitem{kpp}
%
\begin{barticle}[auto:STB|2012/03/12|15:33:09]
\bauthor{\bsnm{Kolmogorov},~\bfnm{A.}\binits{A.}},
\bauthor{\bsnm{Petrovsky},~\bfnm{I.}\binits{I.}} \AND
\bauthor{\bsnm{Piscounov},~\bfnm{N.}\binits{N.}}
(\byear{1937}).
\btitle{\'Etude de l'equation de la diffusion avec croissance de la
quantit\'e
de mati\`ere et son application \`a un probl\`eme biologique}.
\bjournal{Moscow Univ. Math. Bull.}
\bvolume{1}
\bpages{1--25}.
\bptok{imsref}%
\end{barticle}
%
\endbibitem

%b47 ###
%b47 #&#
\bibitem{ls1}
%
\begin{barticle}[mr]
\bauthor{\bsnm{Lalley},~\bfnm{S.~P.}\binits{S.~P.}} \AND
\bauthor{\bsnm{Sellke},~\bfnm{T.}\binits{T.}}
(\byear{1987}).
\btitle{A conditional limit theorem for the frontier of a branching {B}rownian
motion}.
\bjournal{Ann. Probab.}
\bvolume{15}
\bpages{1052--1061}.
\bid{issn={0091-1798}, mr={0893913}}
\bptok{imsref}%
\end{barticle}
%
\endbibitem

%b48 ###
%b48 #&#
\bibitem{lamp1}
%
\begin{barticle}[mr]
\bauthor{\bsnm{Lamperti},~\bfnm{John}\binits{J.}}
(\byear{1967}).
\btitle{The limit of a sequence of branching processes}.
\bjournal{Z. Wahrsch. Verw. Gebiete}
\bvolume{7}
\bpages{271--288}.
\bid{mr={0217893}}
\bptok{imsref}%
\end{barticle}
%
\endbibitem

%b49 ###
%b49 #&#
\bibitem{lamp2}
%
\begin{barticle}[mr]
\bauthor{\bsnm{Lamperti},~\bfnm{John}\binits{J.}}
(\byear{1967}).
\btitle{Continuous state branching processes}.
\bjournal{Bull. Amer. Math. Soc.}
\bvolume{73}
\bpages{382--386}.
\bid{issn={0002-9904}, mr={0208685}}
\bptok{imsref}%
\end{barticle}
%
\endbibitem

%b50 ###
%b50 #&#
\bibitem{lawler}
%
\begin{bbook}[mr]
\bauthor{\bsnm{Lawler},~\bfnm{Gregory~F.}\binits{G.~F.}}
(\byear{2006}).
\btitle{Introduction to Stochastic Processes}, \bedition{2nd} ed.
\bpublisher{Chapman \& Hall/CRC}, \baddress{Boca Raton, FL}.
\bid{mr={2255511}}
\bptok{imsref}%
\end{bbook}
%
\endbibitem

%b51 ###
%b51 #&#
\bibitem{li06}
%
\begin{barticle}[mr]
\bauthor{\bsnm{Li},~\bfnm{Zenghu}\binits{Z.}}
(\byear{2006}).
\btitle{A limit theorem for discrete {G}alton--{W}atson branching
processes with
immigration}.
\bjournal{J. Appl. Probab.}
\bvolume{43}
\bpages{289--295}.
\bid{doi={10.1239/jap/1143936261}, issn={0021-9002}, mr={2225068}}
\bptok{imsref}%
\end{barticle}
%
\endbibitem

%b52 ###
%b52 #&#
\bibitem{ma09}
%
\begin{barticle}[mr]
\bauthor{\bsnm{Ma},~\bfnm{Chunhua}\binits{C.}}
(\byear{2009}).
\btitle{A limit theorem of two-type {G}alton--{W}atson branching
processes with
immigration}.
\bjournal{Statist. Probab. Lett.}
\bvolume{79}
\bpages{1710--1716}.
\bid{doi={10.1016/j.spl.2009.04.008}, issn={0167-7152}, mr={2547941}}
\bptok{imsref}%
\end{barticle}
%
\endbibitem

%b53 ###
%b53 #&#
\bibitem{maillard}
%
\begin{bmisc}[auto:STB|2012/03/12|15:33:09]
\bauthor{\bsnm{Maillard},~\bfnm{P.}\binits{P.}}
(\byear{2011}).
\bhowpublished{The number of absorbed individuals in branching Brownian
motion with a
barrier. Preprint. Available at arXiv:\arxivurl{1004.1426}}.
\bptok{imsref}%
\end{bmisc}
%
\endbibitem

%b54 ###
%b54 #&#
\bibitem{mckean}
%
\begin{barticle}[mr]
\bauthor{\bsnm{McKean},~\bfnm{H.~P.}\binits{H.~P.}}
(\byear{1975}).
\btitle{Application of {B}rownian motion to the equation of
{K}olmogorov--{P}etrovskii--{P}iskunov}.
\bjournal{Comm. Pure Appl. Math.}
\bvolume{28}
\bpages{323--331}.
\bid{issn={0010-3640}, mr={0400428}}
\bptok{imsref}%
\end{barticle}
%
\endbibitem

%b55 ###
%b55 #&#
\bibitem{mmq}
%
\begin{barticle}[mr]
\bauthor{\bsnm{Mueller},~\bfnm{Carl}\binits{C.}},
\bauthor{\bsnm{Mytnik},~\bfnm{Leonid}\binits{L.}} \AND
\bauthor{\bsnm{Quastel},~\bfnm{Jeremy}\binits{J.}}
(\byear{2011}).
\btitle{Effect of noise on front propagation in reaction-diffusion
equations of
{KPP} type}.
\bjournal{Invent. Math.}
\bvolume{184}
\bpages{405--453}.
\bid{doi={10.1007/s00222-010-0292-5}, issn={0020-9910}, mr={2793860}}
\bptok{imsref}%
\end{barticle}
%
\endbibitem

%b56 ###
%b56 #&#
\bibitem{nev87}
%
\begin{bincollection}[mr]
\bauthor{\bsnm{Neveu},~\bfnm{J.}\binits{J.}}
(\byear{1988}).
\btitle{Multiplicative martingales for spatial branching processes}.
In \bbooktitle{Seminar on {S}tochastic {P}rocesses, 1987 ({P}rinceton, {NJ},
1987)}.
\bseries{Progr. Probab. Statist.}
\bvolume{15}
\bpages{223--242}.
\bpublisher{Birkh\"auser}, \baddress{Boston, MA}.
\bid{mr={1046418}}
\bptok{imsref}%
\end{bincollection}
%
\endbibitem

%b57 ###
%b57 #&#
\bibitem{neveu}
%
\begin{bmisc}[auto:STB|2012/03/12|15:33:09]
\bauthor{\bsnm{Neveu},~\bfnm{J.}\binits{J.}}
(\byear{1992}).
\bhowpublished{A continuous-state branching process in relation with
the GREM model of
spin glass theory. Rapport interne  267, \'Ecole polytechnique}.
\bptok{imsref}%
\end{bmisc}
%
\endbibitem

%b58 ###
%b58 #&#
\bibitem{pem}
%
\begin{barticle}[mr]
\bauthor{\bsnm{Pemantle},~\bfnm{Robin}\binits{R.}}
(\byear{2009}).
\btitle{Search cost for a nearly optimal path in a binary tree}.
\bjournal{Ann. Appl. Probab.}
\bvolume{19}
\bpages{1273--1291}.
\bid{doi={10.1214/08-AAP585}, issn={1050-5164}, mr={2538070}}
\bptnote{check year}%
\bptok{imsref}%
\end{barticle}
%
\endbibitem

%b59 ###
%b59 #&#
\bibitem{pit99}
%
\begin{barticle}[mr]
\bauthor{\bsnm{Pitman},~\bfnm{Jim}\binits{J.}}
(\byear{1999}).
\btitle{Coalescents with multiple collisions}.
\bjournal{Ann. Probab.}
\bvolume{27}
\bpages{1870--1902}.
\bid{doi={10.1214/aop/1022677552}, issn={0091-1798}, mr={1742892}}
\bptok{imsref}%
\end{barticle}
%
\endbibitem

%b60 ###
%b60 #&#
\bibitem{py97}
%
\begin{barticle}[mr]
\bauthor{\bsnm{Pitman},~\bfnm{Jim}\binits{J.}} \AND
\bauthor{\bsnm{Yor},~\bfnm{Marc}\binits{M.}}
(\byear{1997}).
\btitle{The two-parameter {P}oisson--{D}irichlet distribution derived
from a
stable subordinator}.
\bjournal{Ann. Probab.}
\bvolume{25}
\bpages{855--900}.
\bid{doi={10.1214/aop/1024404422}, issn={0091-1798}, mr={1434129}}
\bptok{imsref}%
\end{barticle}
%
\endbibitem

%b61 ###
%b61 #&#
\bibitem{revyor}
%
\begin{bbook}[mr]
\bauthor{\bsnm{Revuz},~\bfnm{Daniel}\binits{D.}} \AND
\bauthor{\bsnm{Yor},~\bfnm{Marc}\binits{M.}}
(\byear{1999}).
\btitle{Continuous Martingales and {B}rownian Motion},
\bedition{3rd} ed.
\bseries{Grundlehren der Mathematischen Wissenschaften [Fundamental Principles
of Mathematical Sciences]}
\bvolume{293}.
\bpublisher{Springer}, \baddress{Berlin}.
\bid{mr={1725357}}
\bptok{imsref}%
\end{bbook}
%
\endbibitem

%b62 ###
%b62 #&#
\bibitem{sag99}
%
\begin{barticle}[mr]
\bauthor{\bsnm{Sagitov},~\bfnm{Serik}\binits{S.}}
(\byear{1999}).
\btitle{The general coalescent with asynchronous mergers of ancestral lines}.
\bjournal{J.~Appl. Probab.}
\bvolume{36}
\bpages{1116--1125}.
\bid{issn={0021-9002}, mr={1742154}}
\bptok{imsref}%
\end{barticle}
%
\endbibitem

%b63 ###
%b63 #&#
\bibitem{saw76}
%
\begin{barticle}[mr]
\bauthor{\bsnm{Sawyer},~\bfnm{Stanley}\binits{S.}}
(\byear{1976}).
\btitle{Branching diffusion processes in population genetics}.
\bjournal{Adv. in Appl. Probab.}
\bvolume{8}
\bpages{659--689}.
\bid{issn={0001-8678}, mr={0432250}}
\bptok{imsref}%
\end{barticle}
%
\endbibitem

%b64 ###
%b64 #&#
\bibitem{silv}
%
\begin{barticle}[mr]
\bauthor{\bsnm{Silverstein},~\bfnm{M.~L.}\binits{M.~L.}}
(\byear{1967/1968}).
\btitle{A new approach to local times}.
\bjournal{J. Math. Mech.}
\bvolume{17}
\bpages{1023--1054}.
\bid{mr={0226734}}
\bptnote{check year}%
\bptok{imsref}%
\end{barticle}
%
\endbibitem

%b65 ###
%b65 #&#
\bibitem{ds2}
%
\begin{barticle}[mr]
\bauthor{\bsnm{Simon},~\bfnm{Damien}\binits{D.}} \AND
\bauthor{\bsnm{Derrida},~\bfnm{Bernard}\binits{B.}}
(\byear{2008}).
\btitle{Quasi-stationary regime of a branching random walk in presence
of an
absorbing wall}.
\bjournal{J. Stat. Phys.}
\bvolume{131}
\bpages{203--233}.
\bid{doi={10.1007/s10955-008-9504-4}, issn={0022-4715}, mr={2386578}}
\bptok{imsref}%
\end{barticle}
%
\endbibitem

\end{thebibliography}
\end{document}